\documentclass[a4paper,twoside,10pt]{article}

\usepackage{amsmath}
\usepackage{amsthm}
\usepackage{amsfonts}
\usepackage{subcaption}
\usepackage[margin=0.85in]{geometry}

\usepackage{MnSymbol}%
\usepackage{wasysym}%

\newtheorem{proposition}{Proposition}
\newtheorem{lemma}{Lemma}

\usepackage{graphicx}
\usepackage[utf8]{inputenc}

\usepackage{booktabs}
\usepackage[export]{adjustbox}

\newcommand{\ka}[1]{\textcolor{black}{#1}}
\newcommand{\nm}[1]{\textcolor{black}{#1}}
\newcommand{\an}[1]{\textcolor{black}{#1}}
\newcommand{\ak}[1]{\textcolor{black}{#1}}

\usepackage[prependcaption,colorinlistoftodos]{todonotes}

\DeclareMathOperator*{\argmin}{arg\,min}
\allowdisplaybreaks

\makeatletter
\def\thanks#1{\protected@xdef\@thanks{\@thanks
        \protect\footnotetext{#1}}}
\makeatother

\usepackage[font=small,labelfont=bf]{caption}
\usepackage{float}

\usepackage[ruled,vlined]{algorithm2e}
\usepackage{subcaption}
\usepackage[section]{placeins}
\usepackage{xparse}
\usepackage{authblk}

\title{Optimal intervention strategies to mitigate the COVID-19 pandemic  effects}

\author{Andreas Kasis$^{\ast}$\thanks{Andreas Kasis, Stelios Timotheou and 
Marios Polycarpou are with the KIOS Research and Innovation Center of Excellence and the Department of Electrical and Computer Engineering, University of Cyprus, Cyprus; e-mails: kasis.andreas@ucy.ac.cy, timotheou.stelios@ucy.ac.cy, mpolycar@ucy.ac.cy},  Stelios Timotheou, Nima Monshizadeh\thanks{Nima Monshizadeh is with the Engineering and Technology Institute, University of Groningen, Nijenborgh 4, 9747AG, Groningen, The Netherlands. email: n.monshizadeh@rug.nl} 
 and Marios Polycarpou}

\begin{document}

\date{}

\maketitle

\begin{abstract}
Governments across the world are currently facing the task of selecting suitable intervention strategies to cope with the effects of the COVID-19 pandemic.
This is a highly challenging task, since harsh measures may result in economic collapse while a relaxed strategy might lead to a high death toll.
\ka{Motivated by this, 
we consider the problem of forming intervention strategies to mitigate the impact of the COVID-19 pandemic that optimize the trade-off between the number of deceases and the socio-economic costs. 
We demonstrate that the healthcare capacity and the testing rate highly affect the optimal intervention strategies. 
Moreover, we propose an approach that enables practical strategies, with a small number of  policies and policy changes, that are close to optimal. In particular, we provide tools to decide which policies should be implemented and when should a government change to a different policy. 
Finally, we consider how the presented results are affected by uncertainty in the initial reproduction number and infection fatality rate and 
 demonstrate that parametric uncertainty has a more substantial effect when stricter strategies
are adopted.
}
\end{abstract}


\setlength{\belowdisplayskip}{1.25pt}
\setlength{\abovedisplayskip}{1.25pt}

\section{Introduction}

A novel coronavirus {was first} reported in Wuhan, China {in} December 2019 \cite{li2020substantial}.
The virus, now known as SARS-COV-2 \cite{of2020species},  spread rapidly through China and the rest of the world causing the COVID-19 disease, being officially declared a pandemic by the World Health Organization (WHO) on March 17th, 2020.
Since the outbreak of the COVID-19 pandemic, the world has been facing an unprecedented human tragedy along {with} fears of economic devastation.
As a result, more than \ka{$90$} million infected cases and \ka{$2$} million deaths have been reported to this date {(\today)}.
To cope with the effects of the virus, governments across the world have implemented a range of non-pharmaceutical  interventions such as closing schools, banning public {events} and imposing  social distancing, self-isolation and lockdown policies.
Although such interventions may curtail the infection rate of the disease and hence the spread of the virus \cite{maier2020effective}, \cite{maharaj2012controlling}, they impose an enormous economic effect.
According to the International Monetary Fund  \cite{imf_outlook}, the economic impact  of the pandemic is expected to cause the steepest worldwide recession in over $40$ years
 and result in a loss of more than \ka{$5\%$} of the GDP in the developed world. 
Hence, although a combination of social distancing and lockdown policies may be effective in containing the virus, it might be highly costly in terms of economical impact, which naturally makes government decision making a multi-objective problem.


Mathematical models are fundamental to describe the dynamic evolution of pandemics and to form effective policies to mitigate their impact.
A seminal study in this area is  \cite{kermack1927contribution}, which describes the widely used susceptible - infected - recovered (SIR) model.
 A comprehensive review of epidemiology models can be found in \cite{hethcote2000mathematics}.
Such models enable the study of the progression of various diseases over time, and {facilitate} the characterization of their asymptotic behaviour and dependence on model parameters. 
Recently, there have been various approaches to model the progression of the COVID-19 outbreak. 
A common approach is to apply different extensions to the SIR model, \an{e.g.} \cite{dehning2020inferring}.
In addition, a more involved compartmental model has been developed in \cite{giordano2020modelling}, offering larger modelling flexibility compared to simpler models.
\ka{Furthermore, \cite{della2020network} developed an extended model which took into account the regional heterogeneity of the pandemic.}

Two important parameters in the study of  epidemic progression are the basic reproduction number and the {infection fatality} rate.
The former is interpreted as the  number of new people that the average person  transmits the disease to while the latter enables an estimate of the fatalities resulting from the disease.
The initial reproduction number, i.e. the basic reproduction number at the onset of the disease,  which we denote by $\overline{R}_0$, is also of particular importance to accurately model the disease progression and in deciding the extend of government policies.
However,  there is significant uncertainty in estimating these  parameters,
as demonstrated via numerous studies that estimate $\overline{R}_0$ using statistical data from different countries  \cite{alimohamadi2020estimate}, \cite{yuan2020monitoring}, {\cite{gatto2020spread}},
and various studies that have reported different {infection fatality} rates \cite{verity2020estimates}, \cite{mallapaty2020deadly}, \cite{salje2020estimating}.
Hence, 
 it is important to consider the effect of parametric uncertainty {in forming} effective government {strategies}.


\subsection*{Contribution}

This study {uses tools from optimal control theory to address} the problem of forming a practical and efficient government intervention {strategy} that limits the number of fatalities due to the COVID-19 pandemic with a low social and economic cost until a vaccine is {fully deployed}.

In particular, we consider a {controlled} SIDARE (Susceptible, Infected undetected, infected Detected, Acutely symptomatic - threatened, Recovered, deceased - Extinct) model
\ka{that takes into account the effect of government intervention policies.}
The considered model 
\ka{enables} {the integration of} features such as the impact of the available healthcare capacity and testing rate.
The contribution of this study is summarized {as follows}:
\begin{itemize}
\item[(i)] \ka{\textbf{Fatalities vs economics cost.}}
We \ka{present} the relation between the number of fatalities and cost of optimal government intervention,
and study how this relation is affected by the amount of testing and the capacity of the healthcare system to treat patients.
  \ka{We demonstrate} the effect of these parameters in the decease rate of the pandemic and the resulting cost associated with the optimal intervention strategy.  In addition, for a range of adopted decease tolerance levels,  we provide insights  on the shape of the optimal intervention {strategy} and its dependence on the \ka{adopted test policy}. 
  \item[(ii)]  \ka{\textbf{Which policies and when.}} We consider the fact that a government can only implement a limited number  of policies  and  policy changes  over the time span of the pandemic, due to practicality and implementability reasons
and to avoid the  social fatigue resulting from frequent  changes in policy.
Our approach provides tools to decide which policies should be implemented and when should a government change to a different policy. 
\ka{We demonstrate that  a small number of policies and policy changes yields a close to optimal government strategy.
In particular, our results suggest that the additional cost incurred from implementing $4$ policies and $6$ policy changes  is  less than $1\%$ compared to {the optimal} continuously changing strategy.} 
 \item[(iii)] \ka{\textbf{Impact of uncertainty.}} We consider the impact of uncertainty in {the value} of the initial basic reproduction number $\bar{R}_0$ and the infection fatality rate on  the  decease rates resulting from optimal government {strategies} associated with particular decease tolerance levels.
We demonstrate that parametric uncertainty has a \ka{larger} impact when stricter government policies, {associated with lower \an{decease} tolerances,} are adopted.
\end{itemize}

\section{\ka{Results}}\label{sec_results}

\subsection{\ka{Problem description}}
\begin{figure}[h!]
\begin{center}
\includegraphics[scale=0.40]{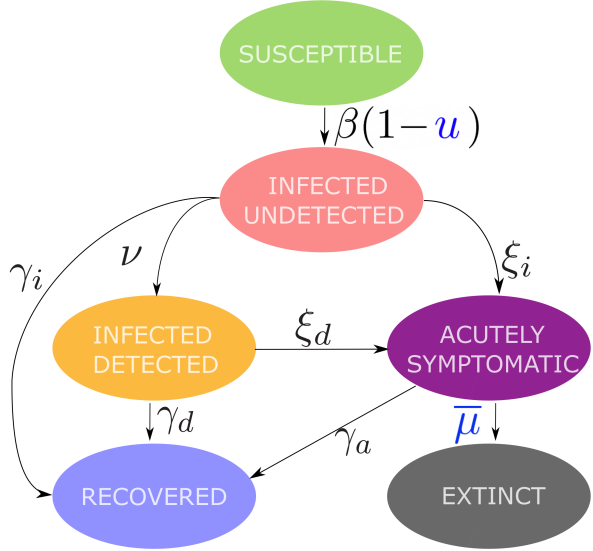}
\end{center}
\vspace{-4mm}
\caption{\textbf{\ka{The controlled SIDARE model.}} Schematic representation of the controlled SIDARE model, used to describe the evolution of the COVID-19 pandemic. \ka{The model splits the population into Susceptible, Infected undetected, infected Detected, Acutely symptomatic - threatened, Recovered and deceased - Extinct.
Model parameters $\beta, \xi_i, \xi_d, \nu, \gamma_i, \gamma_d$ 
and $\gamma_a$ describe the transition rates  between the states.
The effect of government interventions is described by $u$ which limits the rate of infection.
The rate at which the acutely symptomatic population deceases is described by $\overline{\mu}$, which depends on the healthcare system capacity.
}}
\label{SIDARE_figure}
\vspace{-5mm}
\end{figure}

\ka{To study the progression of the pandemic, we consider a controlled SIDARE (Susceptible, Infected undetected, infected Detected, Acutely symptomatic - threatened, Recovered, deceased - Extinct) model (see Fig. \ref{SIDARE_figure}), where the effects of the healthcare capacity limit  and  non-pharmaceutical government interventions on the mortality and infection rates are taken into account.
A full mathematical description of the controlled SIDARE model and explanations on its components are provided in Section \ref{Sec_methods}.
Note that the terms threatened and acutely symptomatic, as well as deceased and extinct are used interchangeably.
In addition, we form a multi-objective optimization problem whose cost function consists of three components: (i) the socio-economic cost of government intervention, (ii) the cost associated with hospitalization and medical care of the acutely symptomatic population and (iii) a cost proportional to the portion of  deceased population. \nm{We seek intervention strategies that minimize the aforementioned cost function. We then investigate how such strategies can be obtained with
a limited} number of distinct policies and policy changes in the described optimization problem.
Further motivation and details concerning the mathematical formulation of the optimization problem is provided in Section \ref{Sec_methods}, while our approach to solve it is detailed in the Appendix.
}

\begin{figure*}[t!]
{
\textbf{a}\begin{subfigure}{0.14\textwidth}
\centering
 \includegraphics[trim = 0mm 0mm 0mm 0mm, scale=0.7]{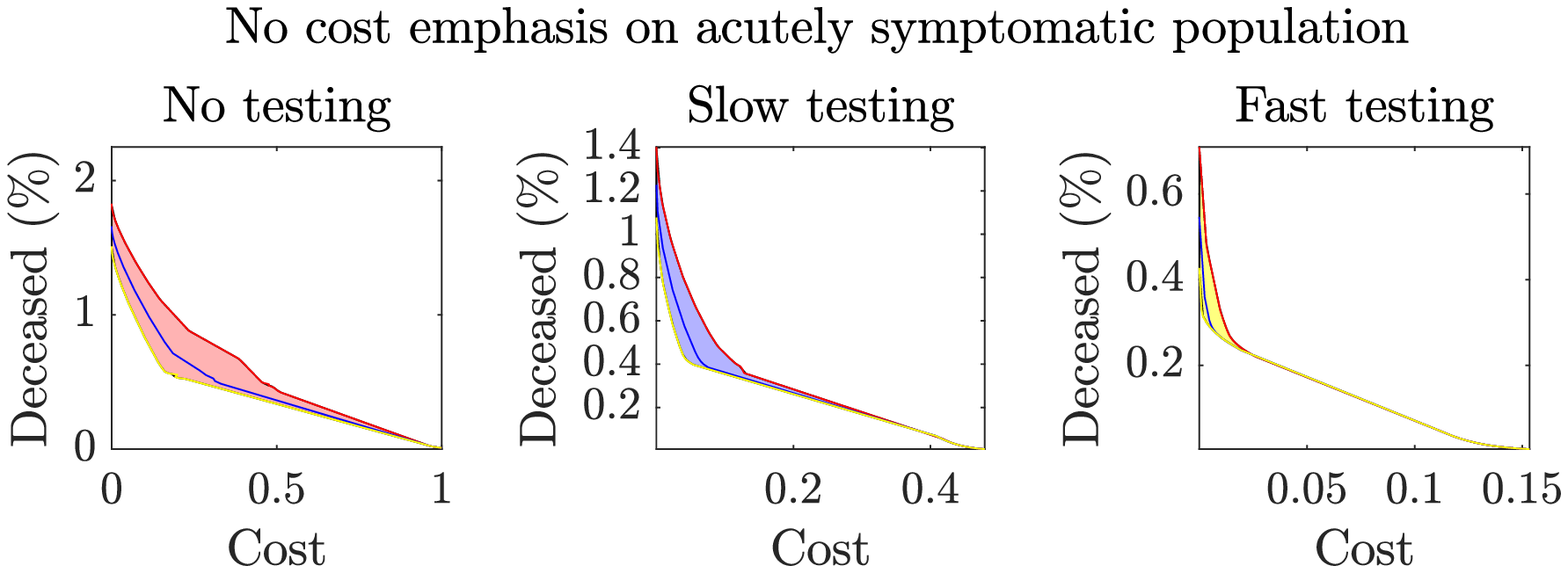}
\label{Cost_Dth_0}
\end{subfigure}
  }
\\[-3mm]
\adjustbox{valign=c}{
 \textbf{b}\begin{subfigure}{0.2\textwidth}
\includegraphics[trim = 0mm 0mm 7mm 0mm,clip, scale=0.7]{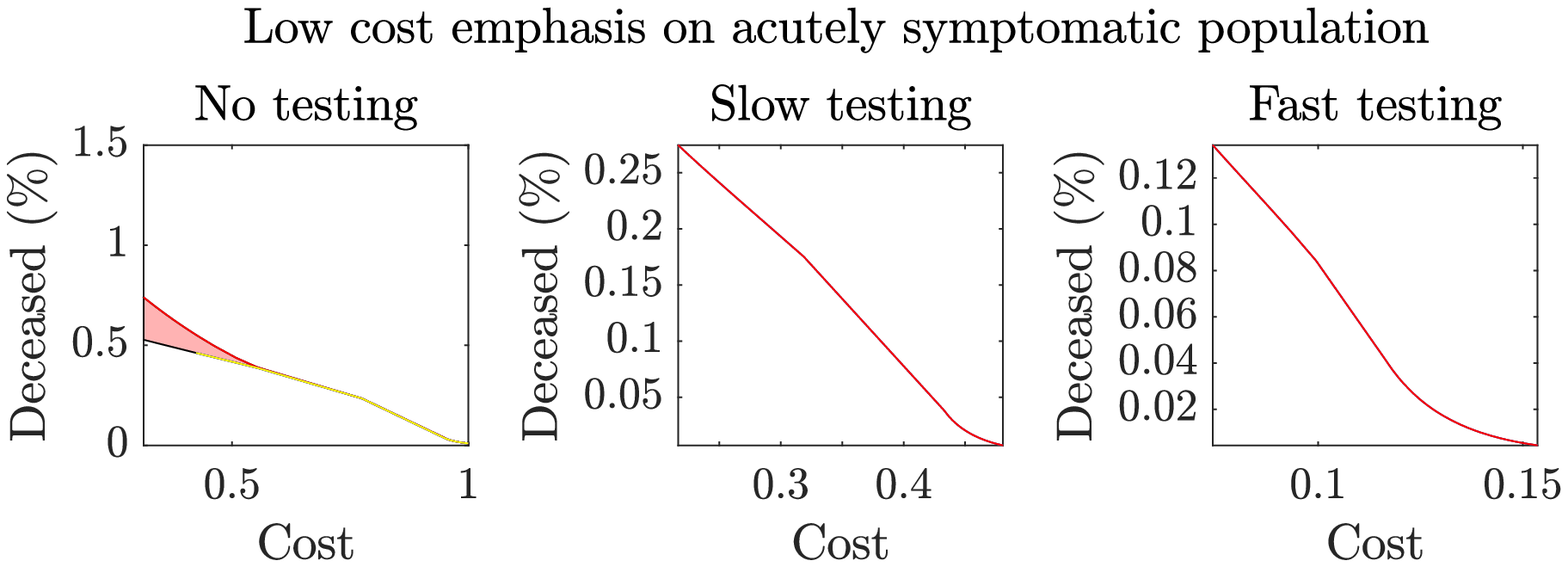}
\label{Cost_Dth_50}
\end{subfigure}
\hspace{0.57\textwidth}
\begin{subfigure}{0.2\textwidth}
\centering
 \includegraphics[trim = 133mm 0mm 0mm 0mm, clip, scale=0.7]{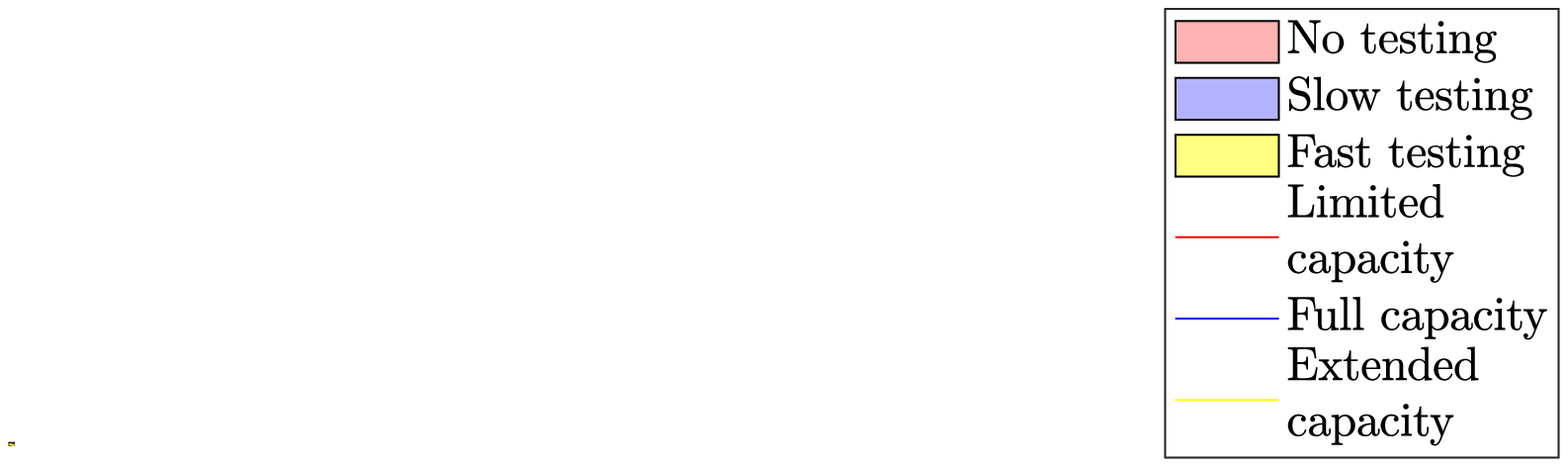}
\end{subfigure}
}
\\[-3mm]
\adjustbox{valign=b}{
\textbf{c}\begin{subfigure}{0.2\textwidth}
\includegraphics[trim = 0mm 0mm 0mm 0mm, scale=0.7]{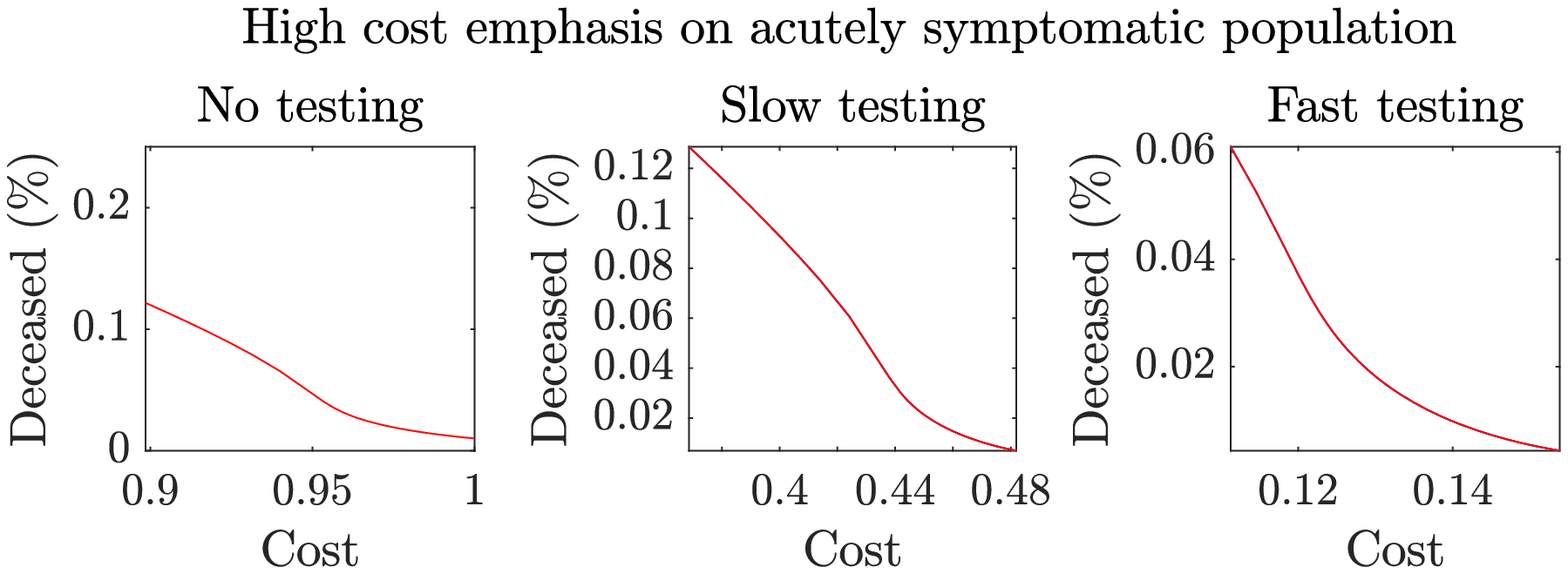}
\label{Cost_Dth_100}
\end{subfigure}
}
\vspace{-6mm}
\caption{\textbf{\ka{Deceased population vs. cost.}} Proportion of deceased population versus cost of optimal government intervention \an{when the available healthcare capacity for COVID-19 patients  is limited (red), full (blue), and extended  (yellow)} and
 when no testing \ka{(left)}, slow testing \ka{(center)} and fast testing \ka{(right)}  policies are adopted.
\ka{In addition, we present the cases where} no \ka{emphasis} 
(\ka{\textbf{a}}), low  \ka{emphasis} (\ka{\textbf{b}}) 
 and high \ka{emphasis} 
 (\ka{\textbf{c}})  is given to the cost associated with the acutely symptomatic population.
\an{When identical relations are obtained for different healthcare capacity levels, as for example in (\textbf{c}), then only the lowest capacity is presented.}
Shaded regions show the ranges of the relations between deceased population percentage and cost of government interaction when the healthcare capacity is between \an{the limited and extended levels}.
\ka{All presented costs are normalised using as basis the cost of the optimal government strategy  with no testing resulting to $0.01\%$ deceases.}
}
\label{Cost_Dth_All}
\vspace{-5mm}
\end{figure*}

\subsection{Deceased population vs. cost of government intervention}

%
\ka{In this section we study the problem of forming an optimal government strategy and its dependence on selected parameters.
In particular, we consider the impact of (i) the healthcare capacity limit, (ii) the testing rate  and (iii) the cost emphasis on the acutely symptomatic and deceased population, on the optimal  intervention strategy  and the resulting portion of deceased population.
 Figure \ref{Cost_Dth_All} depicts the relation between the \an{portion of deceases} and the optimal cost of government intervention, resulting from solving the considered optimization problem (described in equation \eqref{Problem_to_min} in Section \ref{Sec_methods}), for a range of cases for testing rate, healthcare capacity and emphasis on the acutely symptomatic and deceased population.
 It should be noted that all costs presented in Fig. \ref{Cost_Dth_All} are normalised using as basis the cost of the optimal government strategy  with no testing resulting to $0.01\%$ deceases.} 
 


In particular, we considered the cases 
\an{of
(i) limited capacity, where two-thirds of the current healthcare capacity is used for COVID-19 patients, 
(ii) full capacity, where the total capacity is used and 
(iii) extended capacity, where the total capacity is increased by one third due to government investment and is available for COVID-19 patients.}
%
 \ka{In addition, we considered the cases where 
  (I) no testing, 
  (II) slow testing and 
  (III) fast testing  policies are implemented. 
 Finally, we consider three different cases for the cost emphasis on acutely symptomatic population  corresponding to no emphasis (\ka{\textbf{a}}),  low emphasis (\ka{\textbf{b}}) and high emphasis (\ka{\textbf{c}}). 
In addition,   a broad range of cost weights associated with the deceased population
was considered in each case, with aim to provide a rich set of policy options.
 Note  that a zero cost policy is only demonstrated in Fig. \ref{Cost_Dth_All} (\ka{\textbf{a}}) since having a non-zero emphasis on the acutely symptomatic population necessarily results in an optimal intervention strategy with a non-zero cost.
 The exact values used to produce the results presented in Fig. \ref{Cost_Dth_All} are provided in Section \ref{Sec_methods}.}
 
From Fig. \ref{Cost_Dth_All} we deduce the following: \\[1mm]
$(i)$ The healthcare system capacity significantly affects the portion of deceased population, particularly when a low cost (Cost $< 50\%$) strategy with  no testing  is implemented. 
This is particularly reflected in Fig. \ref{Cost_Dth_All} (\ka{\textbf{a}}) \ka{which} demonstrates that 
increasing the available healthcare capacity \ka{from the limited level to the extended level} 
 results in up to a $50\%$ decrease in deceases. 
\\[1mm]
 $(ii)$ When high cost government intervention strategies are adopted (Cost $> 60\%$), then the amount of threatened population never exceeds the healthcare capacity limit and hence its value does not affect the decease rate.
The latter is demonstrated from the fact that there is no shaded regions in most scenarios  depicted in Fig. \ref{Cost_Dth_All} (\ka{\textbf{b,c}}). \\[1mm]
$(iii)$ Increasing the amount of testing enables significantly fewer deaths for the same government intervention cost. 
This is reflected in Fig. \ref{Cost_Dth_All} (\ka{\textbf{a}}), which demonstrates that slow testing approximately halves the \an{portion of deceased} when compared to no testing, when a low intensity government strategy (Cost $< 50\%$) is adopted.
In addition, fast testing results in \ka{approximately} half the deceases compared to slow testing and enables low cost strategies (Cost $< 15\%$).
It should be noted though that, although fast testing policies enable a reduction in costs and decease rates, they may not  always be feasible since they require sufficient resources  in terms of testing equipment and trained personnel.
 \\[1mm]
$(iv)$ When a decease tolerance is set,  a {faster testing} policy enables a less intense government strategy, and hence a lower government intervention cost.
For example, when a $0.1\%$ decease tolerance is considered, a no testing policy requires a cost of more that $90\%$, while slow and fast testing policies \an{yield} the same \an{amount of deceases} with costs of less than $40\%$ and $10\%$ respectively.  


 \begin{figure}[ht!]
\textbf{a}
\begin{subfigure}{0.4\columnwidth}
\centering
\hspace{-15mm}
\includegraphics[scale=0.6]{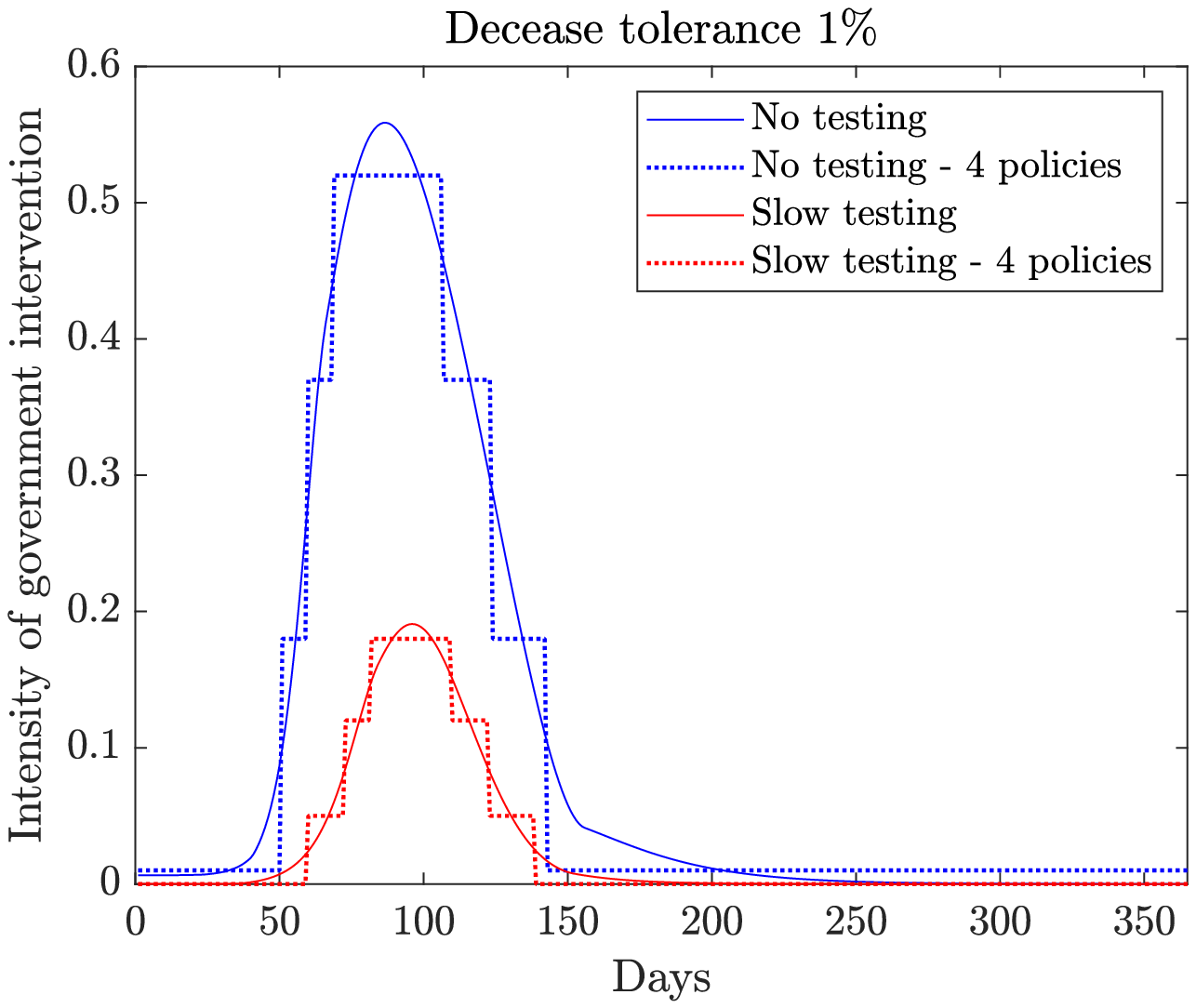}
\label{Policy_1}
\end{subfigure}
\hspace{15mm} \vspace{-3mm}
\begin{subfigure}{0.4\columnwidth}
\centering
\vspace{-3mm}
\includegraphics[scale=0.6]{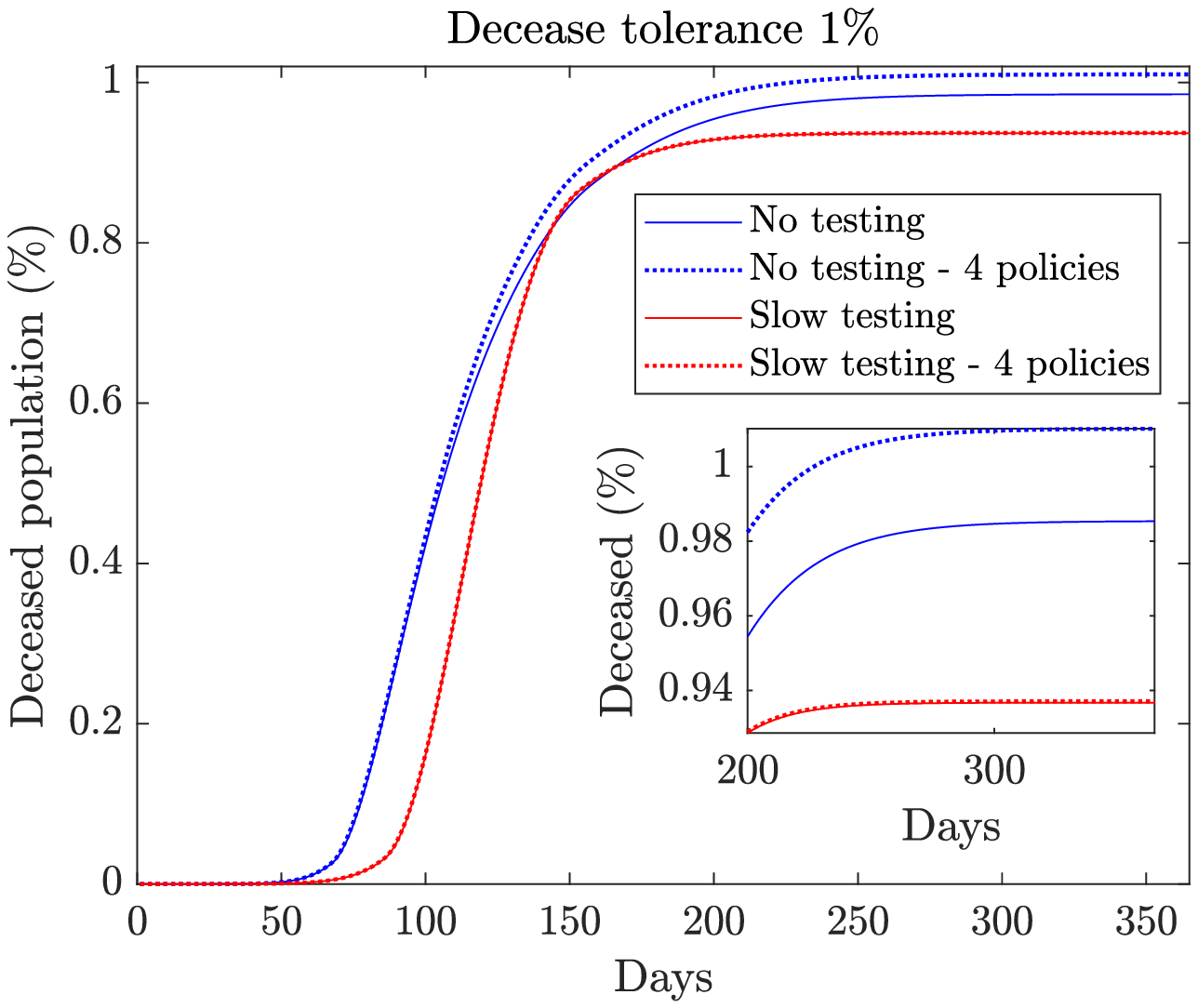}
\end{subfigure}\\
\textbf{b}\begin{subfigure}{0.4\columnwidth}
\centering
\hspace{-15mm}
\includegraphics[scale=0.6]{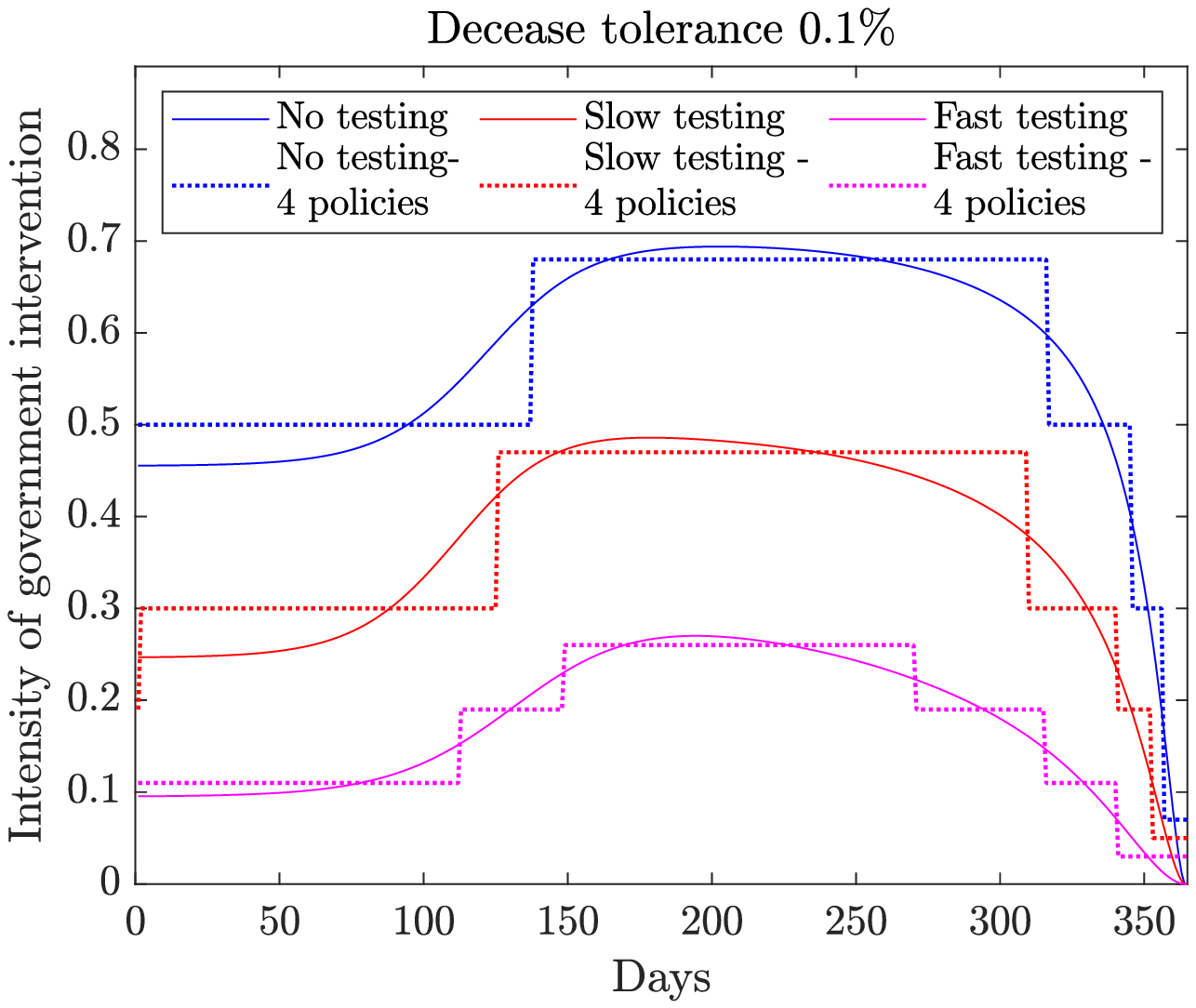}
\label{Policy_01}
\end{subfigure}
\hspace{15mm} \vspace{-3mm}
\begin{subfigure}{0.4\columnwidth}
\centering
\vspace{-3mm}
\includegraphics[scale=0.6]{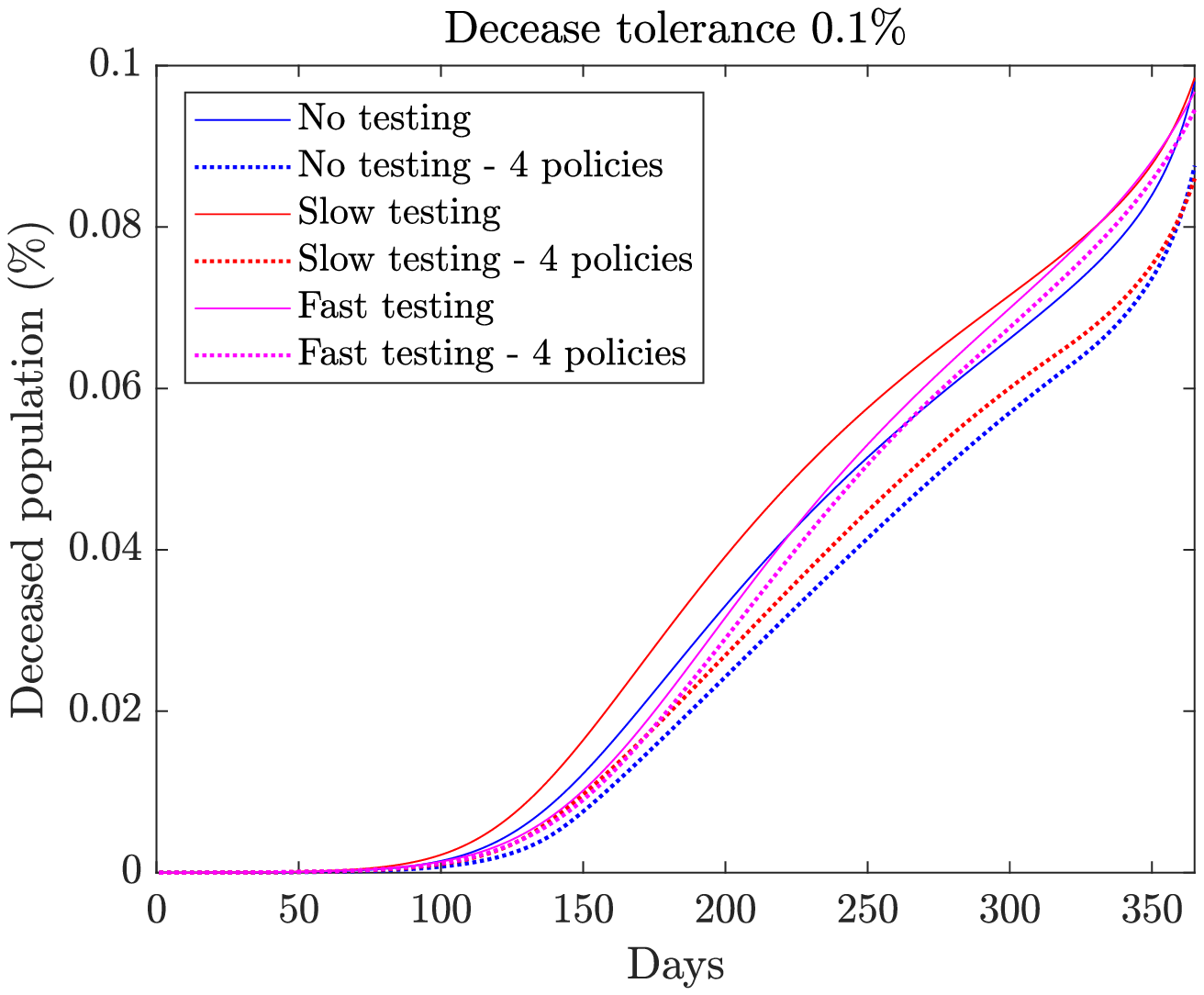}
\end{subfigure}\\
\textbf{c}\begin{subfigure}{0.4\columnwidth}
\centering
\hspace{-15mm} 
\includegraphics[scale=0.6]{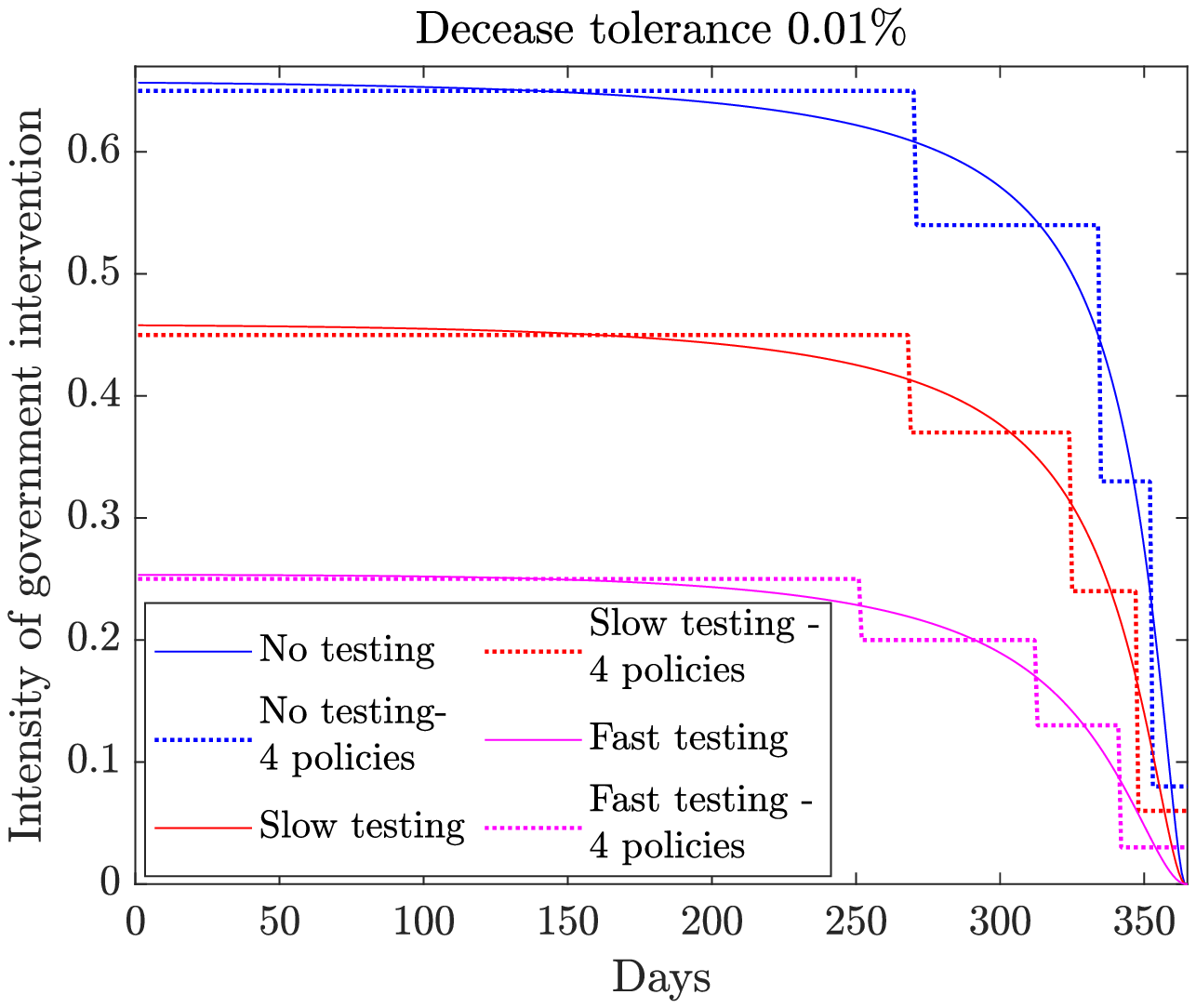}
\label{Policy_001}
\end{subfigure}
\hspace{15mm} \vspace{-7mm}
\begin{subfigure}{0.4\columnwidth}
\centering
\vspace{-3mm}
\includegraphics[scale=0.6]{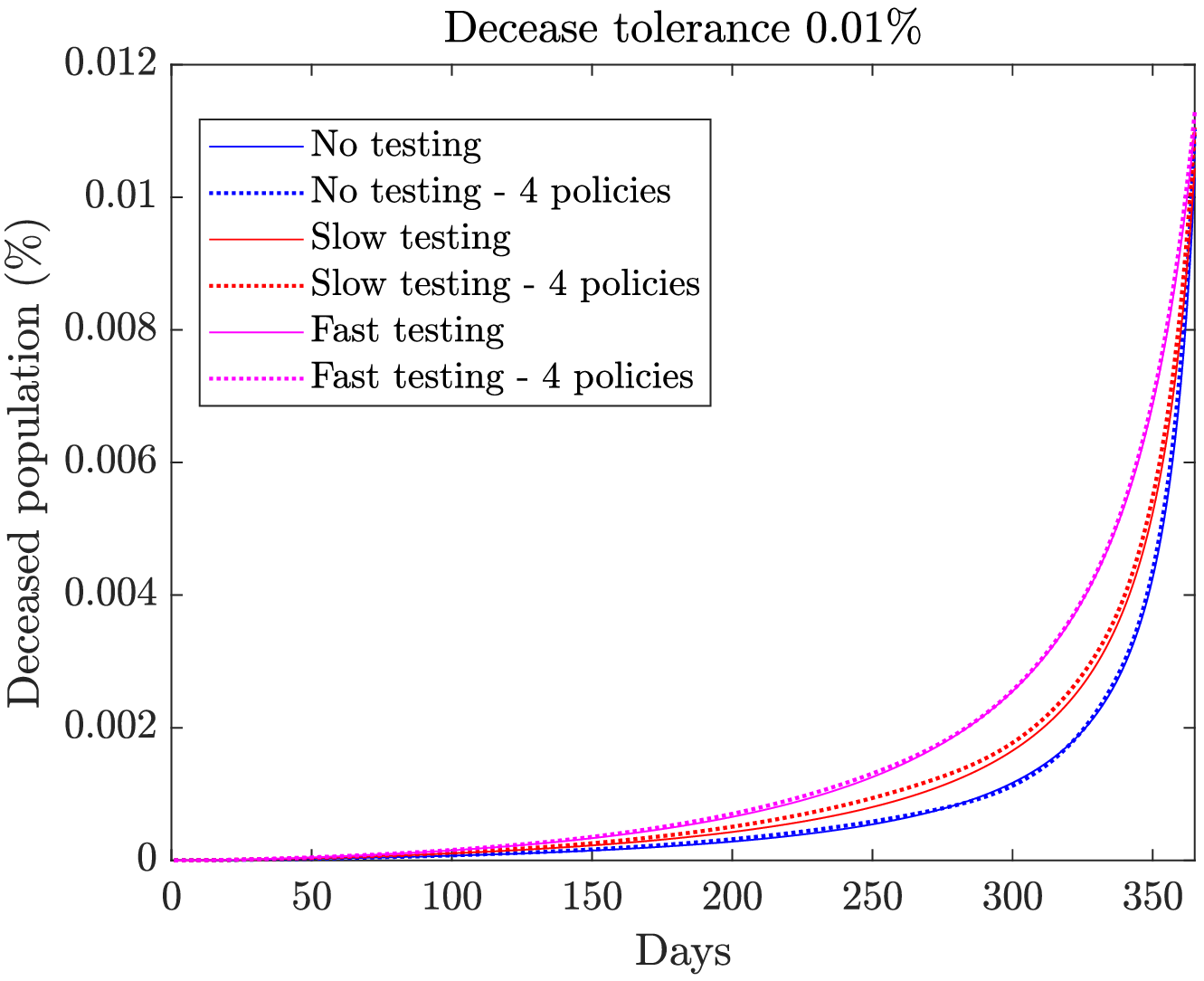}
\end{subfigure}
\caption{\ka{\textbf{Optimal intervention strategies and deceased population.}}
Optimal government intervention strategy (left) and proportion of deceased population (right) versus time for  \ka{government strategies with} (\ka{\textbf{a}}) $1\%$ decease tolerance, (\ka{\textbf{b}})  $0.1\%$ decease tolerance  and (\ka{\textbf{c}}) $0.01\%$ decease tolerance  when 
(i) no testing (blue),
 (ii) slow testing (red) and
  (iii) fast testing (magenta)  policies are adopted.
  \ka{The intensity of government intervention strategies correspond to the strictness of government policies, where $0$ corresponds to no interventions and $1$ to the strictest possible intervention (e.g. a full scale lockdown).
The approach to obtain the optimal continuous strategy is explained in the appendix.}
   Dotted plots correspond to \ka{optimized} discrete implementations of the selected strategies by allowing a maximum of $4$ policy levels and $6$ policy changes.
   \ka{The distinct policy levels and the times where the policy was changed where selected in an optimized way, as described in Section \ref{Sec_methods} and the Appendix.}
Implementations with $7$ and $10$   policy levels and $12$ and $18$ policy changes are \an{presented} in the Appendix.}
\label{Policy_All}
\vspace{-5mm}
\end{figure}

\subsection{Government intervention strategies}



Using the findings depicted in Fig. \ref{Cost_Dth_All}, we aimed to draw \ka{efficient} intervention {strategies} that restrict the portion of the deceased population to specific tolerated amounts with the minimum cost.
The selected portions of decease tolerances where $1\%$, $0.1\%$ and $0.01\%$. 
The approach to obtain the optimal intervention {strategies} is  described in the appendix.
The corresponding \an{intervention strategies} for each decease tolerance level and
 (i) no testing, 
 (ii) slow testing, 
 and 
 (iii) fast testing
 policy levels are depicted in Fig. \ref{Policy_All} (\textbf{a,b,c}, left). 
 Note that a fast testing policy \an{yielded} less than $1\%$ \an{deceases} for any intervention strategy.
 Figure \ref{Policy_All} depicts the intensity of optimal government intervention strategies \ka{(left)} and the resulting \an{portion of deceased} (right) for each decease tolerance level.
 It demonstrates the impact of testing availability in designing intervention strategies associated with selected \an{decease} tolerances.


\ka{The intensity of government intervention is modelled in Section \ref{Sec_methods} with a parameter $u$ (see equation  \eqref{SIDARE_u} in Section \ref{Sec_methods} and Fig. \ref{SIDARE_figure}), where a value of $u = 0$ corresponds to no government interventions and $u=1$ to the strictest intervention policy (e.g. a full scale lockdown).}  
Note that  forming actual government strategies using the provided values of $u$, is a \ka{nontrivial} task.
This problem is equivalent to obtaining  the basic reproduction number resulting from  different intervention policies, which has been considered in \cite{flaxman2020estimating}.
For example, from \cite{flaxman2020estimating} it can be deduced that for Italy, a school closure policy results to $u = 0.02$ while a lockdown policy to $u = 0.8$.
Motivated by this, we consider any policy with $u > 0.6$ as a {\textit{very high}} intensity policy.
In addition, policies with $u \in [0, 0.2]$, $u \in [0.2, 0.4]$ and $u \in [0.4, 0.6]$ are {referred to} as \textit{low}, \textit{medium} and \textit{high} intensity policies respectively.

From Fig. \ref{Policy_All} (\ka{\textbf{a}}), it follows that adopting a high intensity intervention strategy ($u > 0.5$) for a period of approximately $50$ days is required to \an{limit the deceases to}  $1\%$ when no testing is performed. Interestingly, a slow testing policy allows a similar \an{portion of deceased} with a government strategy of about $3$ times lower intensity.  
In addition, Fig. \ref{Policy_All} (\ka{\textbf{b}}) demonstrates that a $0.1\%$ decease tolerance requires intervention strategies ranging from high to very high $(u \in [0.5, 0.7])$, medium to high $(u \in [0.25, 0.50])$, and low to medium $(u \in [0.10, 0.25])$ intensities when no, slow  and fast  testing policies are adopted. 
Furthermore, Fig. \ref{Policy_All} (\ka{\textbf{c}}) shows that a decease tolerance of $0.01\%$ requires very slowly changing strategies of very high ($u \approx 0.65$), high ($u \approx 0.45$) and medium ($u \approx 0.25$) intensity when no, slow and fast  testing policies are respectively adopted.

\subsection{Implementing a limited {number} of policies and policy changes}



 \begin{figure}
\centering
\includegraphics[trim = 0mm 0mm 0 0, scale=0.65]{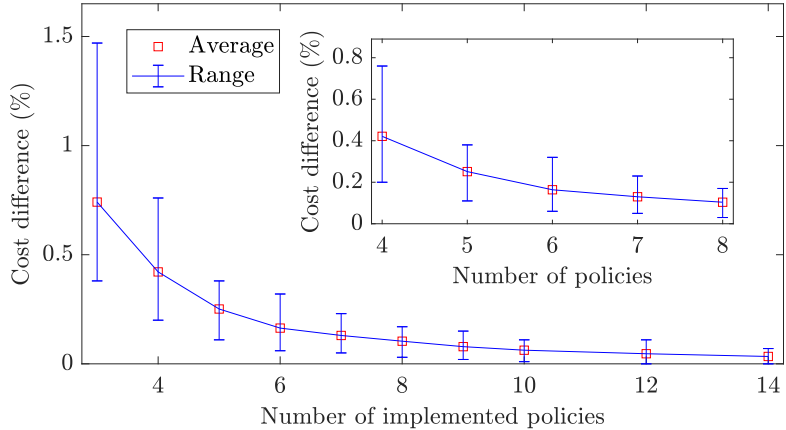}
\caption{\textbf{\ka{Additional cost from implementing a limited number of policies.}} 
\ka{Implementing a limited number of policies results in increased costs compared to the optimal continuously changing intervention strategy. 
This figure depicts the average and range of the percentage differences between the costs of the continuous strategies presented in Fig. \ref{Policy_All} and strategies with   a small number of policies, for different  numbers of  allowed distinct implemented policies.}
%
 In each case, the number of allowed policy changes was twice the number of the implemented policies minus two.
\ka{The boxed plot focuses on implementing between $4$ and $8$ distinct policies and demonstrates that $4$ policies result in a cost difference of less than $1\%$ in all cases.}}
\label{cost_vs_cardinality}
\vspace{-4mm}
\end{figure}

\ka{An implementable government strategy should only have a limited number of distinct policies. 
In addition, frequent changes in the intervention strategy may result in social fatigue {and confusion}, decreasing the  {receptiveness} of the population to the policy instructions.}
 The policies that follow by discretizing  the continuous strategies described in the previous section are  depicted in Fig. \ref{Policy_All} with dotted plots.
 \an{The approach to obtain optimized strategies with a small number of policies and policy changes is described in the appendix.}
Figure \ref{Policy_All} demonstrates that implementing $4$ {policies} 
 and allowing a maximum of $6$ changes among them results in {almost identical levels of decease rates compared to the continuously changing strategies}.
Hence, a close to optimal government response may be obtained with a {relatively} small number of distinct policies.


{Implementing an optimal strategy with a limited number of policies and policy changes results in an increased cost compared to the optimal continuously changing strategy.
The cost differences for discrete implementations of the strategies presented in Fig.  \ref{Policy_All} are depicted in Fig. \ref{cost_vs_cardinality}, where a broad range of allowed number of policies  is considered.} 
{For all cases the number of allowed policy changes \ka{was twice the number of implemented policies minus two.}
 From Fig. \ref{cost_vs_cardinality}, it follows that as the number of policies grows, the percentage cost difference decreases.
 Furthermore, as follows from the boxed plot within Fig. \ref{cost_vs_cardinality}, it can be seen that a low cost difference can be obtained with a small number of policies.
 In particular, adopting  $4$ or more intervention policies   allowed a cost difference of less than $1\%$. 
 The latter demonstrates the effectiveness of implementing a small number of policies and policy changes. 
 
 Hence, a small number of policies  and policy changes suffices for a close to optimal government response, while at the same time resolves issues of   implementability and  social fatigue.

\begin{figure*}[t!]
\vspace{0mm}
\hspace{0.2mm}
  \begin{subfigure}[c]{0.01\textwidth}
   \hspace{-2mm} \textbf{a}
  \end{subfigure}
\begin{subfigure}{0.14\textwidth}
\centering
\includegraphics[scale=0.497]{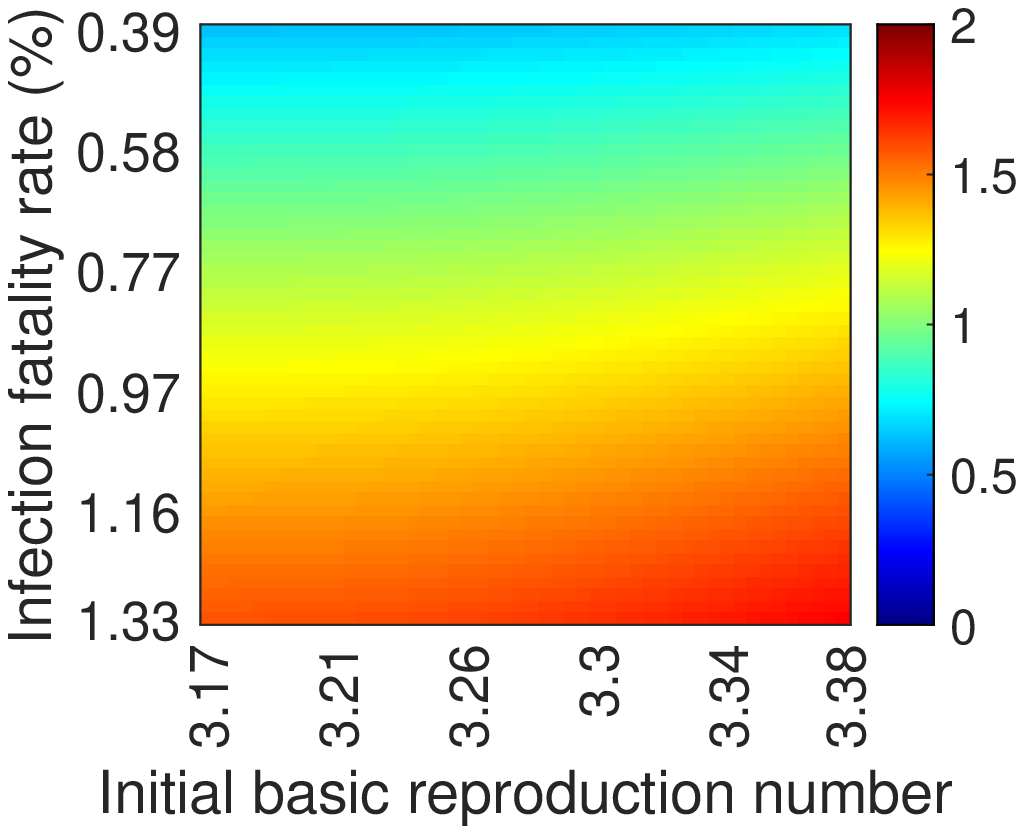}
\label{Heat_1}
\end{subfigure}
\hspace{0.189\textwidth}
\begin{subfigure}{0.20\textwidth}
\centering
\includegraphics[scale=0.497]{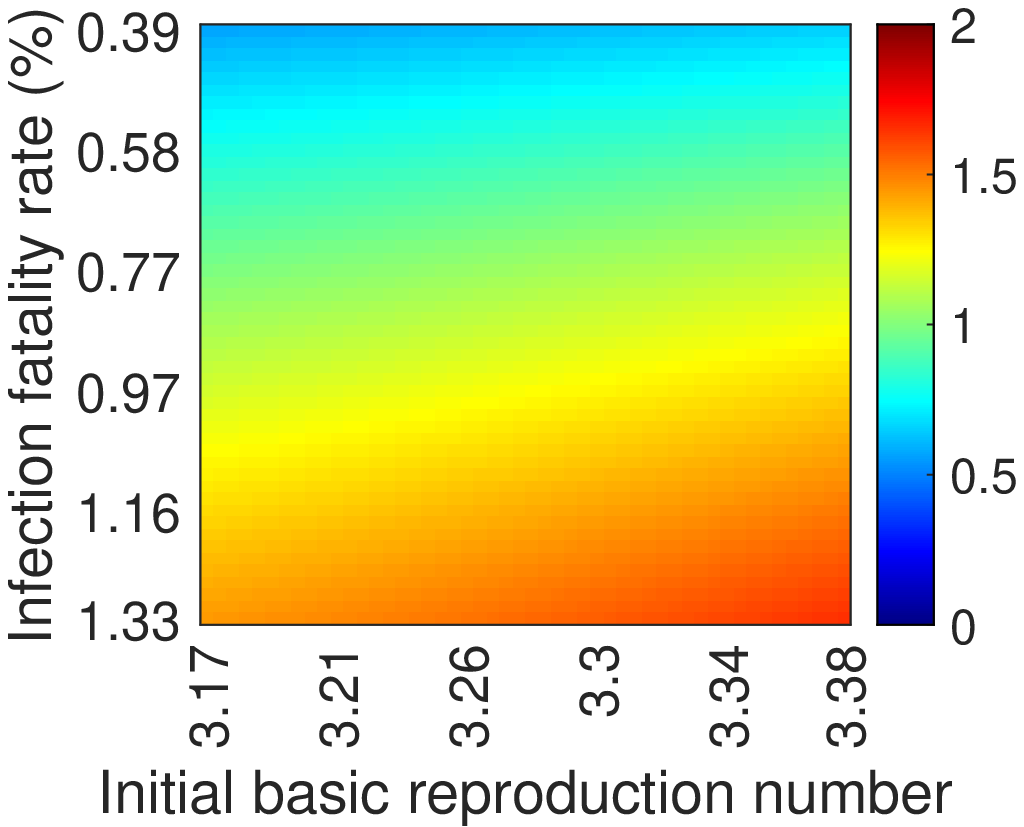}
\label{Heat_2}
\end{subfigure}
\hspace{0.11\textwidth} 
  \footnotesize{
  \begin{tabular}{p{2.68cm} p{1.75cm}}
  \toprule
      Description & Deceased  \\
    \midrule
    1$\%$ tol. - no test & $1.74\%$  \\
    1$\%$ tol. - slow test & $1.65\%$  \\
    0.1$\%$ tol. - no test & $0.29\%$  \\
    0.1$\%$ tol. - slow test & $0.29\%$  \\
    0.1$\%$ tol. - fast test & $0.28\%$  \\
    0.01$\%$ tol. - no test & $0.057\%$  \\
    0.01$\%$ tol. - slow test & $0.085\%$  \\
    0.01$\%$ tol. - fast test & $0.116\%$  \\
    \bottomrule
    \vspace{6mm}
    \end{tabular}
  }
\\
{
  \begin{subfigure}[c]{0.01\textwidth}
   \hspace{-2mm}  \normalsize{\textbf{b}}
  \end{subfigure}
\begin{subfigure}{0.2\textwidth}
 \includegraphics[scale=0.497]{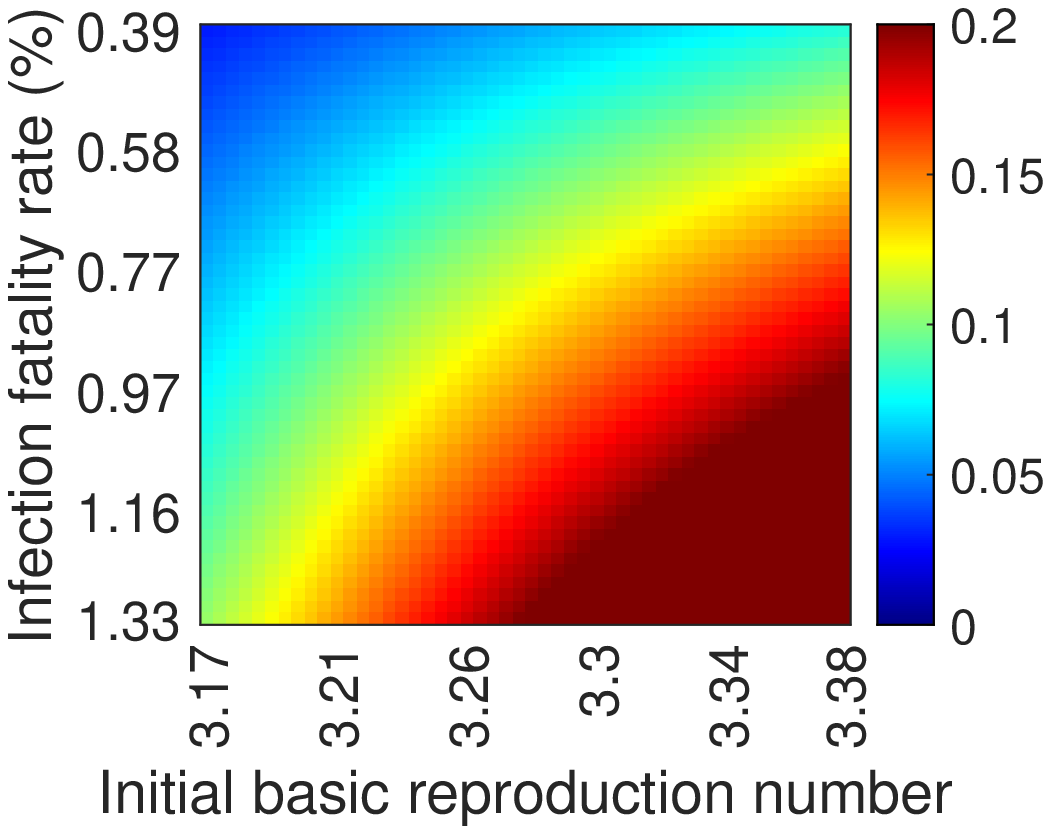}  
\label{Heat_3}
  \end{subfigure} 
\hspace{0.125\textwidth}
\begin{subfigure}{0.20\textwidth}
\centering
\includegraphics[scale=0.497]{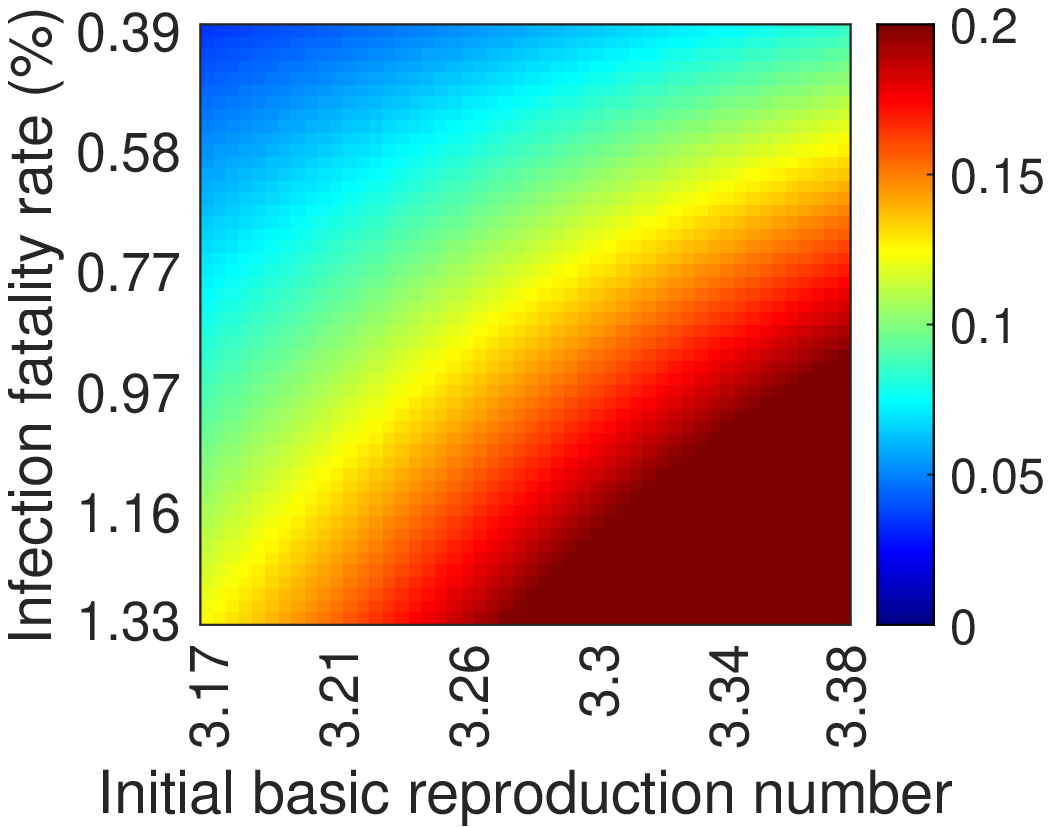}
\label{Heat_4}
\end{subfigure}
\hspace{0.125\textwidth}
\begin{subfigure}{0.20\textwidth}
\centering
\includegraphics[scale=0.497]{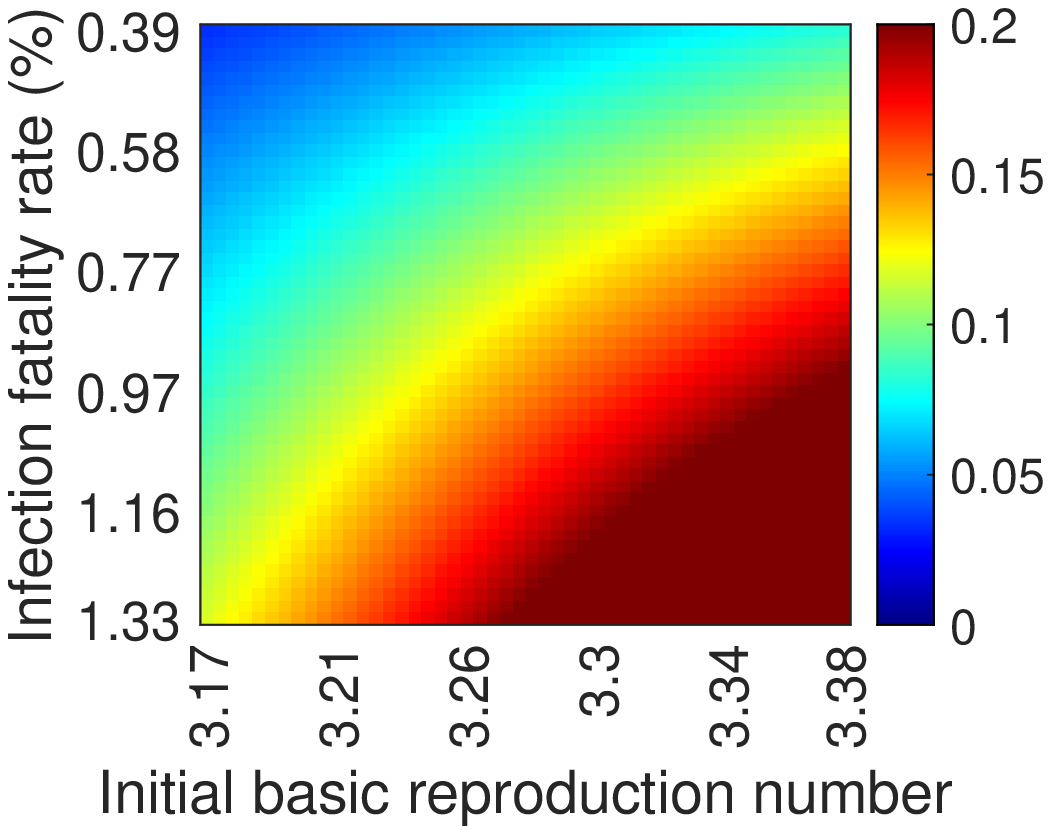}
\label{Heat_5}
\end{subfigure}
}
\\[0.0mm]
{
  \begin{subfigure}[c]{0.01\textwidth}
   \hspace{-2mm} \normalsize{\textbf{c}}
  \end{subfigure}
\begin{subfigure}{0.2\textwidth}
\includegraphics[scale=0.497]{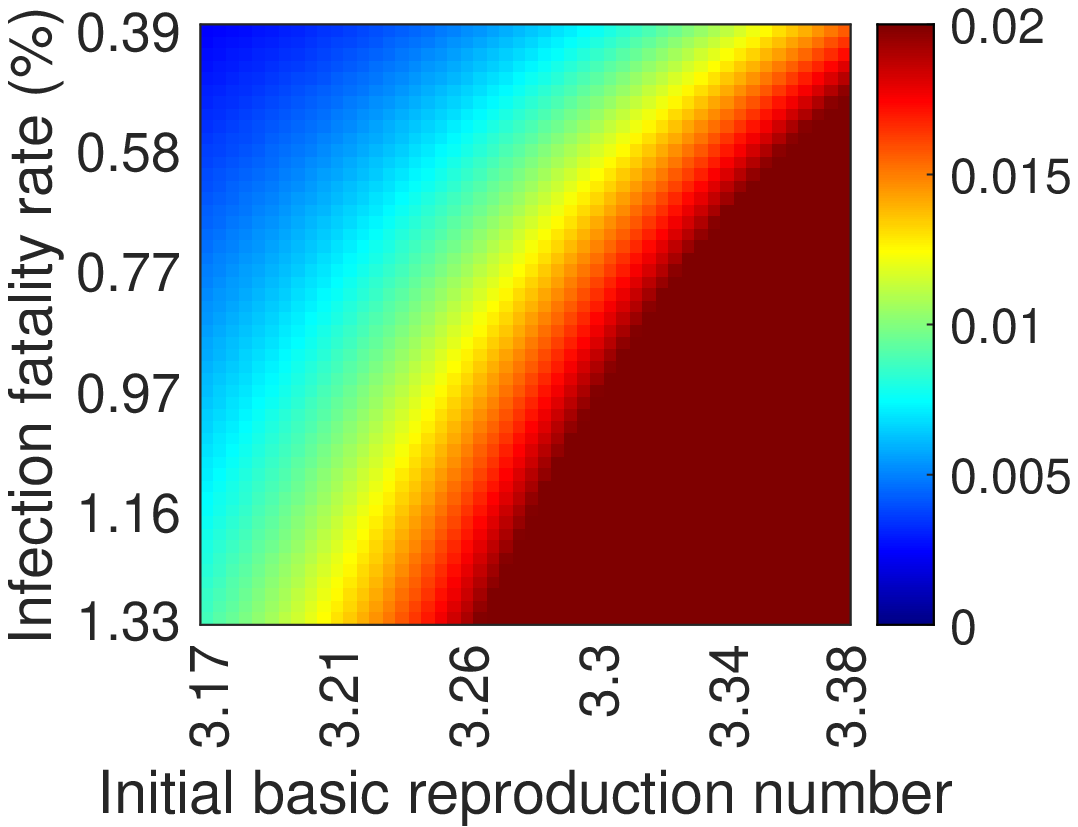} 
\label{Heat_6}
\end{subfigure}
\hspace{0.125\textwidth}
\begin{subfigure}{0.20\textwidth}
\centering
\includegraphics[scale=0.497]{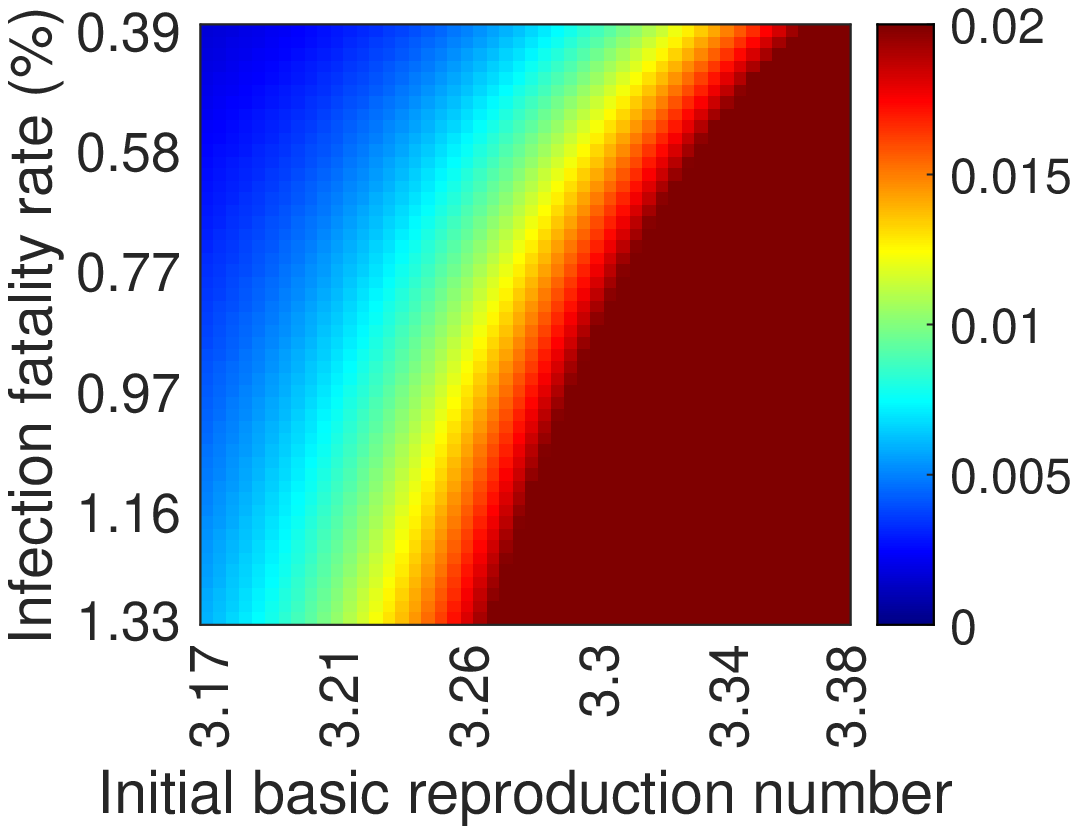}
\label{Heat_7}
\end{subfigure}
\hspace{0.125\textwidth}
\begin{subfigure}{0.20\textwidth}
\centering
\includegraphics[scale=0.497]{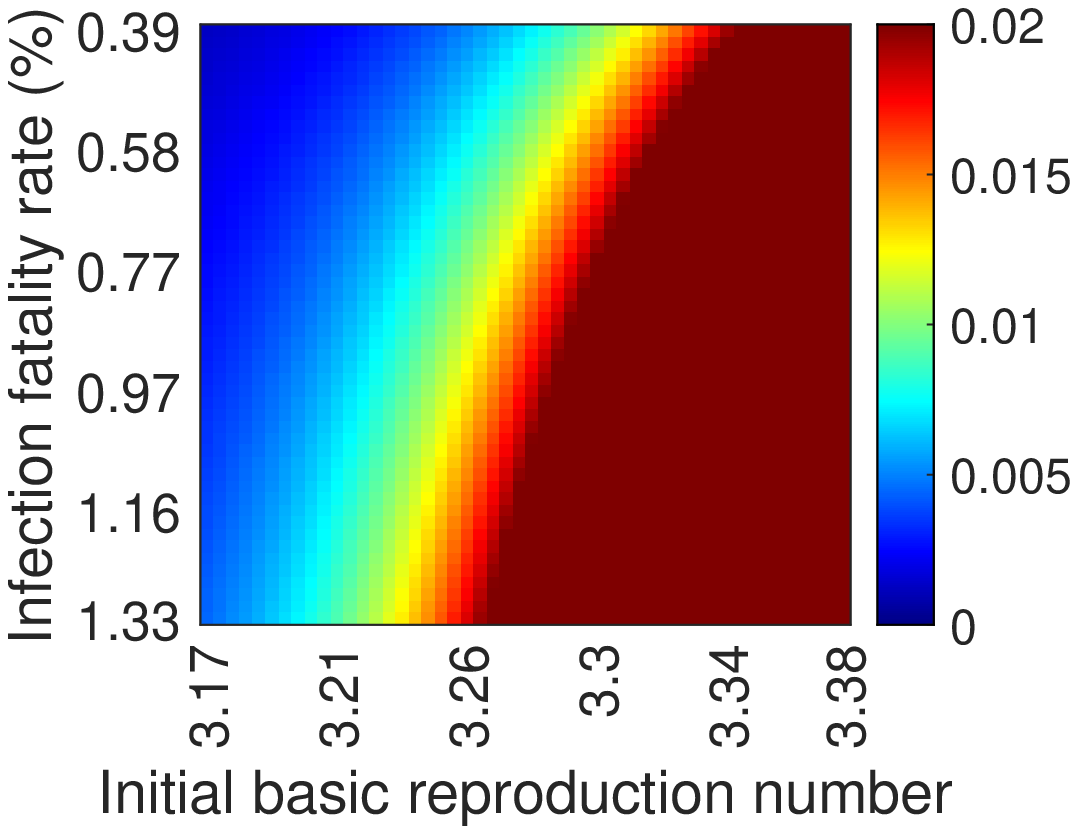}
\label{Heat_8}
\end{subfigure}
}
\vspace{-3mm}
\caption{\textbf{\ka{Effect of uncertainty in the initial reproduction and infection fatality rates.}} \an{Portion of deceased} for $\bar{R}_0 \in [3.17, 3.38]$ and infection fatality rate ranging between $0.39\%$ and $1.33\%$ associated with
decease tolerances of $1\%$ (\ka{\textbf{a}}), $0.1\%$ (\ka{\textbf{b}}) and $0.01\%$ (\ka{\textbf{c}}) \ka{when} no (left) slow (center) and fast (right) testing policies \ka{are implemented}.
\ka{Fast testing always \an{limits the deceases to} less than $1\%$ and hence there is no corresponding case.}
The \ka{table on (\textbf{a})} reports the worst case \an{portion of} deceased  associated with each strategy, obtained
at the maximum considered values of $\bar{R}_0$ and the infection fatality rate.
\ka{The values where acquired by applying in each case the optimal continuous government intervention strategy obtained with $\bar{R}_0 = 3.27$ and infection fatality rate of $0.66\%$.}
 \ka{Dark red and dark blue colours correspond to aggregate deceases that are at least twice and less than half the adopted tolerance levels}. 
 Additional results that demonstrate the impact of parametric uncertainty in forming effective government mitigation strategies are provided in the Appendix.}
\vspace{-3mm}
\label{Heat_All}
\end{figure*}

\subsection{Effect of {parametric} uncertainty}

{The design of optimal control strategies relies on the use of mathematical models. Therefore, a} critical aspect in designing government mitigation strategies is their dependence on {parametric} uncertainty, i.e. the extend of the effect of inaccurately estimating model parameters.
In this section, we consider how the uncertainty in {the value of} the initial basic reproduction number $\bar{R}_0$ and the infection fatality rate  affect the amount of deceases resulting from the strategies presented in the previous section.

In particular, we considered the effect on the \an{aggregate deceases} when the value of $\bar{R}_0$ ranges between $3.17$ and $3.38$ and when the infection fatality rate ranges between $0.39\%$ and $1.33\%$, which correspond to $95\%$ confidence intervals, as reported in \cite{yuan2020monitoring} and \cite{verity2020estimates} respectively.

Figure \ref{Heat_All} demonstrates the effect of {parametric} uncertainty on the \an{portion of} deceased  \an{resulting} from each of the $8$ considered government intervention strategies {presented in Fig. \ref{Policy_All}}, i.e. it depicts the \an{portion of} deceased from implementing the selected strategies when $\bar{R}_0$ and \ka{the infection fatality rate} have been imprecisely estimated.
Figure \ref{Heat_All} demonstrates that the level of adopted decease tolerance is crucial when it comes to the effect of model uncertainty in the decease rate.

In particular, when a $1\%$ decease tolerance level is adopted, the obtained policy enabled a moderate percentage increase in the \an{total} deceases, resulting to $1.75\%$ deceases in a worst case scenario, as demonstrated in Fig. \ref{Heat_All} (\ka{\textbf{a}}).
 When the decease tolerance level was decreased to $0.1\%$, then {parametric} uncertainty  could result in  up to $0.29\%$ deceases, i.e. about $3$ times higher, as demonstrated in Fig. \ref{Heat_All} (\ka{\textbf{b}}) and the table in Fig. \ref{Heat_All} (\ka{\textbf{a}}). 
 Finally, when a requirement for decease tolerance level of $0.01\%$ was imposed, then the considered uncertainty  could result in up to $11.6$ times higher \an{values}, as depicted in Fig. \ref{Heat_All} (\ka{\textbf{c}}) and the table in Fig. \ref{Heat_All} (\ka{\textbf{a}}).
These show that {parametric} uncertainty should also be taken into account in forming government policies, particularly when the adopted policies are stricter, i.e. when a low decease tolerance level is imposed.

\section{Discussion}\label{sec_discussion}

\ka{Following the COVID-19 outbreak, governments across the world have adopted  strict intervention policies to contain the pandemic.
However, the high economic costs resulting from these policies have sparred debates on the necessity of the measures and on how these could be relaxed without risking a new wave of infections.}
Recently, several approaches have been proposed  to control the spread of the COVID-19 pandemic.
In particular,  optimal  intervention strategies that simultaneously minimize the number of fatalities and the economic costs are presented in \cite{djidjou2020optimal} and \cite{rowthorn2020cost}.
A similar problem has been considered in \cite{kohler2020robust}, which in addition investigated model predictive control approaches.
Moreover, \cite{tsay2020modeling} and \cite{bin2020fast} explored the use of on-off policies to mitigate the effects of the pandemic, proposing control strategies that  alternate between  no measure and full measure policies.
{Furthermore, \cite{della2020network} considered regional instead of national interventions to alleviate the effects of the COVID-19 pandemic.}
In addition, \cite{alvarez2020simple} considered the problem of selecting the optimal lockdown portion that jointly minimizes the number of fatalities of the COVID-19 pandemic and the economic costs associated with the proposed policy by considering an adapted SIR model.
A similar problem has been studied in \cite{acemoglu2020optimal}, which in addition considered the effect of imposing different policies to different age groups.

\ka{An important missing aspect that we consider in this study is the fact that governments can only impose a limited amount of intervention policies, for practicality and implementability issues and to reduce the social fatigue resulting from frequently changing policies.
We demonstrate that a small number of distinct policies and policy changes yields a close to optimal government strategy.
In particular, our results suggest that {the additional cost incurred from implementing $4$ policies and $6$ policy changes  is  less than $1\%$ compared to {the optimal} continuously changing strategy.}
Using tools from optimal control theory, we provide an approach, analytically described in the appendix (see Section \ref{sec_sol_meth}), which allows to select \ka{in an optimized fashion}, which policies should be implemented and when should a government switch to a different policy.} 
\ka{Our approach can be easily applied to different types of models and cost functions, which might give emphasis to other aspects of the pandemic.}

\ka{A further contribution of this paper is the study of the uncertainty in the parameters associated with the initial reproduction number and the infection fatality rate. 
Our results suggest that parametric uncertainty has a more significant effect  when
stricter government policies, aiming for lower decease rates, are adopted.
}

\section{Methods}\label{Sec_methods}


{We} consider a SIDARE model, which is a  variation of the SIR model, to describe the evolution of the COVID-19 pandemic, where the population is divided in six categories: (i) Susceptible to be infected, (ii) Infected but undetected, (iii) infected and Detected, (iv) Acutely symptomatic - threatened, (v) Recovered and (vi) Extinct - deceased.
Note that we use the terms threatened and acutely symptomatic, as well as deceased and extinct interchangeably.



 The dynamics of the SIDARE model are given by
\begin{subequations}\label{SIDARE}
\begin{gather}
\dot{s} = -\beta s i, \label{SIDARE_s} \\
\dot{i} = \beta s i - \gamma_i i -\xi_i i -\nu i, \label{SIDARE_i} \\
\dot{d} = \nu i - \gamma_d d - \xi_d d, \label{SIDARE_d} \\
\dot{a} = \xi_i i + \xi_d d - \gamma_a a - \mu a, \label{SIDARE_a} 
\\
\dot{r} = \gamma_i i + \gamma_d d + \gamma_a a, \label{SIDARE_r} 
\\
\dot{e} = \mu a, \label{SIDARE_e} 
\end{gather}
\begin{gather}
s(0) = s_0, i(0) = i_0, d(0) = d_0, a(0) = a_0, r(0) = r_0, e(0) = e_0,
\end{gather}
\end{subequations}
where $s, i, d, a, r, e \in [0,1]$ are the states of the system describing the portions of susceptible, infected - undetected, infected - detected, threatened,  recovered and deceased population respectively.
Moreover, $s_0, i_0, d_0, a_0,  r_0, e_0 \in [0, 1]$ denote the initial values for $s, i, d, a, r, e$ respectively. 
The model parameters are briefly summarized below:
\begin{itemize}
\item $\beta$ describes the infection  rate for susceptible individuals.
\item $\gamma_i, \gamma_d$ and $\gamma_a$ describe the recovery rates for infected undetected, infected detected and threatened individuals.
\item $\nu$ denotes the rate of detection of infected individuals, associated with the adopted level of testing.
\item $\xi_i$ and $\xi_d$ describe the rates at which {infected undetected and infected detected} individuals become acutely symptomatic.
\item $\mu$  describes the mortality rate of the disease, i.e. the rate at which acutely symptomatic individuals decease.
\end{itemize}
Note that all  model parameters are assumed non-negative and constant.  
The SIDARE model is based on the following assumptions:\\[1mm]
$(i)$ Those recovered are no longer susceptible to the disease. \\[1mm]
$(ii)$ The considered population is constant, i.e. no births or deaths not attributed to COVID-19 are taken into account.\\[1mm]
$(iii)$ The considered country (or region) is isolated, i.e. no imported cases are taken into account.\\[1mm]
$(iv)$ Infected individuals that are detected are assumed to be \ka{quarantined}, i.e. they do not contribute to new infections, something justified by existing practices. \\[1mm]
$(v)$ Infected individuals become acutely symptomatic before they decease. \\[1mm]
$(vi)$ Acutely symptomatic individuals  require hospitalization since they are considered threatened for decease.

The assumption of constant population suggests that the states satisfy $s + i + d + a + r + e = 1$ at all times
and hence that one state is redundant since it can be described by the remaining states at all times. 
In the analysis below, we select $r$ to be the redundant state, satisfying $r = 1 - s - i - d - a - e$.

\subsection{Impact of healthcare capacity on mortality rate}
An important aspect that we consider is the impact of the healthcare system capacity on the mortality rate. 
It is evident that when the healthcare capacity is exceeded, then the mortality rate of the population increases.
The latter is modelled in equation \eqref{heathcare_mu}, 
which suggests that the mortality rate depends on the portion of acutely symptomatic population by the relation 
\begin{equation}\label{heathcare_mu}
\overline{\mu}(a)  = \begin{cases}
\mu a, \text{ if } a \leq \overline{h}, \\
\mu \overline{h} + \hat{\mu}(a - \overline{h}), \text{ if }  a > \overline{h},
\end{cases}
\end{equation}
where the function $\overline{\mu} : [0,1] \rightarrow \mathbb{R}_+$  describes the mortality  of the acutely symptomatic population.
The values of $\mu$ and $\hat{\mu}$ satisfy $\mu < \hat{\mu}$ and correspond to the mortality rates when the healthcare system satisfies the demand and when the healthcare capacity is exceeded by much. This means that when the infected population increases, the mortality rate tends to $\hat{\mu}$.
In addition, the value of $\overline{h}$ describes the existing healthcare capacity.
Note that for simplicity we assume a constant value of $\overline{h}$, although its value could  rise in the future due to a possible increase in the healthcare system  capacity.

\subsection{Modelling government interventions on the SIDARE model}

To account for the effect of the government  actions to mitigate the spread of the pandemic, we introduce an intervention input $u$ to the SIDARE model.
Its value affects  the infection rate of the disease, $\beta$, resulting in a slower spread.
 The controlled SIDARE model follows from equation \eqref{SIDARE} when $\beta$ is replaced by $\beta(1-u)$ and in addition includes the healthcare capacity impact on the mortality rate, described by equation \eqref{heathcare_mu}. Its dynamics are given by
\begin{subequations}\label{SIDARE_u}
\begin{gather}
\dot{s} = -\beta s i (1 - u), \label{SIDARE_s_u} \\
\dot{i} = \beta s i (1 - u) - \gamma_i i - \xi_i i - \nu i, \label{SIDARE_i_u} \\
\dot{d} = \nu i - \gamma_d d - \xi_d d, \label{SIDARE_d_u} \\
\dot{a} = \xi_i i + \xi_d d - \gamma_a a - \overline{\mu} (a), \label{SIDARE_a_u} \\
\dot{e} = \overline{\mu} (a), \label{SIDARE_e_u} \\
s(0) = s_0, i(0) = i_0, d(0) = d_0, a(0) = a_0, e(0) = e_0, \\
s + i + d + a + r + e = 1,
\end{gather}
\end{subequations}
where $u \in \mathcal{U} = [0, \bar{u}]$ and \ka{$\bar{u} \leq 1$} is a  positive constant that denotes the maximum value that the intervention policy $u$ is allowed to take.
Since government actions should only aid in curtailing the effects of the pandemic, $u$ is not allowed to take negative values.
The dynamics of the controlled SIDARE model are depicted in Fig.~\ref{SIDARE_figure}.

The value of $u(t)$ corresponds to the government intervention policy at time $t$, with higher values of $u$ corresponding to stricter intervention policies.
 For example, when a government does not take any action, then $u = 0$ and when a government takes the strictest possible measures, e.g. when implementing a full scale lockdown, then $u = \bar{u}$.
%

\subsection{A multi objective optimization problem}
A suitable government strategy should aim to simultaneously minimize the number of fatalities and the costs associated with implementing intervention policies.
The number of the aggregate fatalities during the considered period is given by $e(T)$, where the constant $T > 0$ denotes the considered timeframe.


Moreover, any policy  $u$ comes with a cost associated with the social and economic side effects from its implementation.
For example, a lockdown policy has an economic cost due to the inability of a portion of the population to work and a social cost associated with restricting the population movements and interactions.
In addition, we consider the cost associated with the acutely symptomatic population. The latter describes the costs resulting from people requiring additional care, including possible hospitalization. These motivate the following cost functional,
\begin{equation*}
C(u,a) =  \int_{0}^{T}  \frac{1}{2}u(t)^2 dt + \theta_a \int_{0}^{T}  \frac{1}{2}a(t)^2 dt,
\end{equation*}
where the non-negative parameter $\theta_a $ describes the weight given on the cost associated with the threatened population.
The proposed cost functional sets a penalty analogous to the square of  the intervention effort $u$,   set by the government to mitigate the effects of the disease, and  the square of the aggregate threatened population $a$. 
Note that a quadratic cost is considered in order to enable a close estimate to the non-linear cost effects arising from intense government strategies and from having a large portion of the population being in a threatened state.

However, there is a {trade-off} between minimizing the economic cost of government policies and the number of fatalities.
 The above motivates the following
 optimization problem
 \begin{equation}\label{Problem_to_min}
 \begin{aligned}
& \min_{u \in \mathcal{U}} J(a,e,u) = C(u,a) + \theta_e e(T) \\
& \text{s.t. } \eqref{SIDARE_u},    
\end{aligned}
 \end{equation} 
where $\theta_e$  describes the weight given to the total number of deaths  in comparison with the cost associated with the threatened population and government intervention effort.  
The values of weight coefficients $\theta_e$ and $\theta_a$ are key to form the optimal policy.
For example, if  $\theta_e = \theta_a = 0$, then the focus of the government is to minimize the {socio-economic} cost of the intervention strategy, which trivially results to $u = 0$ for all times. On the other hand, when $\theta_a$ and $\theta_e$ are large, then the focus becomes to minimize the number of fatalities and the number of acutely symptomatic individuals, which results in a value of  $u$ that is close to $\bar{u}$  at all times. Since there is a {trade-off} between these objectives, selecting suitable values for $\theta_e$ and $\theta_a$ is highly important. 
Furthermore, note that the relative ratio between $\theta_a$ and $\theta_e$ enables an extra degree of freedom in the choice of the optimization problem and a richer set of solutions.
The approach to solve \ka{the above optimization problem}, using tools from optimal control theory, is explained in the appendix \an{(see Section \ref{sec_sol_meth})}.



\subsection{Implementing a limited {number} of policies and policy changes}


\an{A government can only implement a limited number  of policies  and  policy changes  over the time span of the pandemic, for practicality and to avoid the  social fatigue resulting from frequent policy changes.
To account for this,} we restrict {both the number} of distinct implemented policies and policy changes in the {previously} considered optimization problem.

{We} denote the set of possible policies by $\mathcal{U}_d$, the set of distinct policies within strategy $u$ by $\mathcal{R}(u) = \{\tilde{u}: \exists t \in [0,T] \text{ s.t. } u(t) = \tilde{u}\}$ and the set of  switching instants by $\mathcal{T} = \{t \in [0, T]: \lim_{\epsilon \rightarrow 0} u(t - \epsilon) \neq \lim_{\epsilon \rightarrow 0} u(t + \epsilon)\}$.
The revised problem is given by
 \begin{equation}\label{Problem_to_min_discrete}
 \begin{aligned}
& \min_{u \in \mathcal{U}_d} J(a,e,u) \\[-1mm]
& \text{s.t. } \eqref{SIDARE_u},     |\mathcal{R}(u)| \leq \hat{n}_1, |\mathcal{T}| \leq \hat{n}_2,
\end{aligned}
 \end{equation}
 where $\hat{n}_1$ and $\hat{n}_2$  denote the maximum number of policies that $u$ is allowed to take from the set $\mathcal{U}_d$ and the maximum allowed number of changes in the intervention strategy over the considered timeframe respectively.
The solution approach \ka{for problem} \eqref{Problem_to_min_discrete} is analytically described within the Appendix.

 \subsection{Simulated parameters}

In this section, we describe and justify the parameters considered in the \ka{simulation results presented in Section \ref{sec_results}}. 
We use the controlled SIDARE model, \ka{described by equation \eqref{SIDARE_u},} for our simulations, with a time horizon of $T = 365$ days.
The selected initial conditions correspond to the very early stage of the disease, where $0.001\%$ of the population has been infected and there are no detected cases, acutely symptomatic cases, fatalities or recoveries yet.
When possible, data associated with the pandemic in Italy where used \ka{for consistency}.

The values of $\gamma_i$ and $\gamma_d$ were selected following \cite{who_report}, which  suggests a median time of disease onset to recovery for mild cases of approximately two weeks.
The value of $\gamma_a$ was selected following \cite{wang2020phase}, which reported a median time between hospitalization and recovery of $12.4$ days.
Furthermore, to select the values for $\xi_i$ and $\xi_d$, 
we used the findings from \cite{verity2020estimates} on hospitalization rate per age group and  data for the Italian population age distribution  \cite{Italy_UN}.
In addition, \an{for the results presented in Fig. \ref{Policy_All}--\ref{Heat_All},} we considered \ka{a healthcare capacity of} $333$ care beds per $100,000$ habitants following \cite{rhodes2012variability}, \an{which corresponded to the full capacity case}.
The maximum allowed value for $u$, given by {$\bar{u}$}, was selected to be $0.8$ to account for the fact  that complete isolation is impossible, since always some critical units will need to remain operational.

The value of $\beta$ was selected following an  initial basic reproduction number  of approximately $3.27$ as estimated in \cite{yuan2020monitoring} \an{and the relation $\overline{R}_0 = \beta s_0/ (\gamma_i + \xi_i + \nu)$ which is analytically shown in the appendix \an{(see Section \ref{sec_SIDARE})}, assuming $\nu = 0$ at $t = 0$ days.}
The value of $\mu$  was selected following a median infection fatality rate of $0.66\%$, as reported  in \cite{verity2020estimates}. 
We let $\hat{\mu}$, which corresponds to the fatality rate when the healthcare system capacity is overloaded, be $5$ times higher than $\mu$, motivated by the findings in \cite{catena2020relocate}.
\ka{In addition, the values for $\theta_a$ associated with no, low and high emphasis on acutely symptomatic population where $0, 5 \times 10^4$ and $10^5$ respectively. 
For each case, a broad range of cost coefficients associated with the deceased population was considered, letting $\theta_e \in [0, 2.5 \times 10^4]$.
The no, slow and fast testing policies corresponded to values of $\nu$ of $0, 0.05$ and $0.10$ respectively.}
Additional explanations on the simulated parameters are provided in the Appendix.

 \section{Conclusion}
We considered the problem of forming government intervention strategies that optimize the {trade-off} between the number of deceases and the social and economic costs. 
We demonstrate the relation between the number of fatalities and cost of optimal government intervention, and how this depends on the adopted testing policy and the healthcare system capacity.
Moreover, we determine that a small number of policies and policy changes suffices for a close to optimal intervention strategy.
In particular, \ka{our results suggest} that implementing $5$ policies and $8$ policy changes  enables a cost difference of less than $1\%$ compared to a continuously changing strategy.
Finally, we considered the impact of uncertainty in the initial reproduction number and infection fatality rate and demonstrated that its effect is more severe when intense government strategies are implemented.

\newpage

\setcounter{section}{0}
\renewcommand{\thesection}{\Alph{section}}

\section*{Appendix}

This appendix describes the methodology used to obtain the main results of the paper.
In particular, we provide analytical results that describe the behaviour of the SIDARE model and explain how we obtain the optimal government intervention strategy for given parameters. Moreover, we present a detailed algorithm that yields an optimized strategy with a limited number of policies and policy changes. Furthermore, we provide additional details regarding the presented numerical simulations  and include extra simulation results with aim to enable additional clarity and intuition on the presented findings.


\subsection*{Notation}

The set of natural numbers is denoted by $\mathbb{N}$.
The projection of a scalar $x$ to a {set} $X$ is denoted by $[x]_X = \argmin_{y \in X} |y - x|$. 
The interior of a set $X$ is denoted by $X^{\mathrm{o}}$.
The cardinality of a discrete set $\Sigma$ is denoted by $|\Sigma|$.
A monotonically increasing (respectively decreasing) function $f:\mathbb{R}\rightarrow \mathbb{R}$ satisfies $f(x) \leq f(y)$ (respectively $f(x) \geq f(y)$) for all $x,y$ such that $x \leq y$.
{The notation $x \succeq 0$ implies that all elements of vector $x$ are  non-negative.}
In addition, we use $\dot{x}$ to denote the time derivative of a signal $x$.

\ak{
\section{Analysis of the SIDARE model}\label{sec_SIDARE}
}

\ak{
\subsection{Analytical results for the SIDARE model}
}

The SIDARE model, given by  \ak{equation} \eqref{SIDARE}, is a bilinear system with six states.
Below, we demonstrate that all states \ak{in equation} \eqref{SIDARE} take non-negative values, provided that the initial conditions are non-negative, which is clearly in line with intuition.

\begin{lemma}\label{lemma_non_negative}
Consider \ak{equation} \eqref{SIDARE} and let all states take non-negative values at $t = 0$. Then, all states take non-negative values for all times.
\end{lemma}

\emph{Proof of Lemma \ref{lemma_non_negative}:}
If all states takes non-negative values at the initial conditions, then from \ak{equation} \eqref{SIDARE_s}, it follows that $s$ can not take negative values by continuity and the fact that $\dot{s} = 0$ when $s = 0$.
It then follows that $i$ is also non-negative from $s$ being non-negative and \ak{equation} \eqref{SIDARE_i}.
Similar arguments hold for states $d$ and $a$, i.e. the non-negativity of $i$ and $(i,d)$ respectively and the dynamics in \ak{equations} \eqref{SIDARE_d}--\eqref{SIDARE_a} suffice so that $d$ and $a$ are non-negative.
Finally, the non-negativity of $(i,d,a)$ implies that $\dot{r} \geq 0$ and $\dot{e}  \geq 0$ from \ak{equations} \eqref{SIDARE_r} and \eqref{SIDARE_e} respectively and hence that $r$ and $e$ are also non-negative.
\hfill $\blacksquare$

It should also be noted that from \ak{equations} \eqref{SIDARE_s}, \eqref{SIDARE_r}, and \eqref{SIDARE_e} and Lemma \ref{lemma_non_negative}, it follows that the state $s$ is monotonically decreasing and the states $r$ and $e$ are monotonically increasing.

The following lemma characterises the equilibria of \ak{equation} \eqref{SIDARE}.
Note that we let $x = (s, i, d, a, r, e)$ for convenience. In addition, we denote an equilibrium of \eqref{SIDARE} by $x^* = (s^*, i^*, d^*, a^*, r^*, e^*)$.

\begin{lemma}\label{equilibrium_lemma}
The set of equilibria of  \eqref{SIDARE} is given by $S = \{x^* : i^*= d^*= a^* = 0 \}$.
\end{lemma}

\emph{Proof of Lemma \ref{equilibrium_lemma}:}
Any equilibrium of \ak{equation} \eqref{SIDARE} should satisfy $s^* i^* = 0$ from \ak{equation} \eqref{SIDARE_s}. This suggests from \ak{equation} \eqref{SIDARE_i} that $i^* = 0$, which in turn suggests from  \ak{equation} \eqref{SIDARE_d} that $d^* = 0$. An equilibrium with $i^* = d^* = 0$  also satisfies $a^* = 0$  from \ak{equation} \eqref{SIDARE_a}.
Moreover, consider any point in the state space of \ak{equation} \eqref{SIDARE} with $i = d = a = 0$. This point is necessarily an equilibrium since all \ak{equations} \eqref{SIDARE_s}--\eqref{SIDARE_e} are zero.
\hfill $\blacksquare$

The following proposition enables to define the basic reproduction number $\bar{R}_0$ in terms of model parameters. 
Note that the definitions of  stable and unstable equilibria follow from \cite[Definition 4.1]{khalil1996nonlinear}.
{Within Proposition \ref{lemma_R0}, we make use of the set $X_{\succeq 0} = \{x: x \succeq 0\}$.}

\begin{proposition}\label{lemma_R0}
Consider an equilibrium point of \eqref{SIDARE}, given by $x^{\ast} = (s^{\ast}, 0, 0, 0, r^{\ast}, e^{\ast})$ and let $x(0)$ take non-negative values.
Then,
\begin{itemize}
\item[(i)]  if  $s^{\ast} < (\gamma_i + \xi_i + \nu)/ \beta$, then there exists a neighbourhood $Z$ of $x^{\ast}$ such that solutions initiated in {$Z \cap X_{\succeq 0}$} asymptotically converge to an equilibrium point within $Z$.
\item[(ii)]  if  $s^{\ast} > (\gamma_i + \xi_i + \nu)/ \beta$, then $x^\ast$ is an unstable equilibrium point.
\end{itemize}
\end{proposition}

\emph{Proof of Proposition \ref{lemma_R0}:}
We first prove part (i), using Lyapunov based arguments, and then part (ii).

\emph{Part (i):}
{First note that since all states take non-negative values at $t = 0$, then they take non-negative values at all times from Lemma \ref{lemma_non_negative}.}
Furthermore,  from \ak{equations} \eqref{SIDARE_s}, \eqref{SIDARE_r}, and \eqref{SIDARE_e} it follows that the state $s$ is monotonically decreasing and the states $r$ and $e$ are monotonically increasing.

For the system \ak{described by equation} \eqref{SIDARE}, we consider an equilibrium point $x^*$ and the Lyapunov candidate function 
\begin{equation*}
V(x) = V_1(s,r,e) + V_2(i,d,a),
\end{equation*} 
where
\begin{equation*}
V_1(s,r,e) = \frac{1}{2} (s - s^*)^2 +  \frac{1}{2} (r - r^*)^2 + \frac{1}{2} (e - e^*)^2,
\end{equation*}
and
\begin{equation*}
V_2(i,d,a) =  i + \theta_1 d + \theta_2 a, 
\end{equation*}
where $\theta_1$ and $\theta_2$ are positive constants.

The dynamics of \ak{equation} \eqref{SIDARE} suggest that 
\begin{equation*}
\dot{V}_1(x) =  (s - s^*)(-\beta s i) + (r - r^*)(\gamma_i i + \gamma_d d + \gamma_a a) + (e - e^*)\mu a \leq 0
\end{equation*}
where the inequality holds since $s(t) \geq s^*$, $r(t) \leq r^*$ and $e(t) \leq e^*$ for all times since $s$ is monotonically decreasing and $r$ and $e$ are monotonically increasing and all states {$i,d,a$} and parameters $\beta, \gamma_i, \gamma_d, \gamma_d$ and $\mu$ take non-negative values.

In addition, it follows that along trajectories of \ak{equation} \eqref{SIDARE} it holds that
\begin{equation*}
\dot{V}_2(x) = \beta s i - \gamma_i i -\xi_i i -\nu i + \theta_1(\nu i - \gamma_d d - \xi_d d) + \theta_2(\xi_i i + \xi_d d - \gamma_a a - \mu a).
\end{equation*}
Note that for solutions initiated sufficiently close to $x^*$ it holds that $s(t) < (\gamma_i + \xi_i + \nu)/ \beta$ at all times due to the monotonicity of $s$.
Therefore, letting $\rho = (\gamma_i + \xi_i + \nu) - \beta s(0)$ it follows that $\dot{i} \leq -\rho i$ at all times and consequently
\begin{equation*}
\dot{V}_2(x) \leq -\rho i + \theta_1(\nu i - \gamma_d d - \xi_d d) + \theta_2(\xi_i i + \xi_d d - \gamma_a a - \mu a).
\end{equation*}
Hence, there exist positive $\theta_1, \theta_2$ satisfying $\theta_1 \nu + \theta_2 \xi_i < \rho$ and $\theta_2 \xi < \theta_1(\gamma_d + \xi_d)$ such that
\begin{equation}\label{dot_V}
\dot{V} = \dot{V}_1 + \dot{V}_2 \leq -\phi_1 i - \phi_2 d - \phi_3 a \leq 0,
\end{equation}
where $\phi_1, \phi_2$ and $\phi_3$ are positive constants.
Hence, there exists a compact connected set $\Xi$ which includes $x^*$, given by $\Xi = \{x: V(x) \leq r\}$, for sufficiently small $r > 0$, such that solutions initiated within $\Xi$ remain in $\Xi$ for all times. 
Since $r$ can be selected to be arbitrarily small, it follows that $x^*$ is a stable equilibrium of  \eqref{SIDARE}.

LaSalle's invariance principle \cite{khalil1996nonlinear}[Theorem 4.4] can now be applied on the compact and positively invariant set $\Xi$.
This guarantees that solutions to \ak{equation} \eqref{SIDARE} initiated in $\Xi$ converge to the largest invariant set within
 $\Xi \cap \{x: \dot{V}(x) = 0\}$.
If $\dot{V} = 0$ holds within $\Xi$ it follows from \ak{equation} \eqref{dot_V} that $(i,d,a) = (0,0,0)$.
The latter implies convergence of the states $(s,r,e)$ from \ak{equations} \eqref{SIDARE_s} and \eqref{SIDARE_r}--\eqref{SIDARE_e} to the set of equilibria within $\Xi$.

In addition, note that if $r$ in the definition of $\Xi$ is selected to be sufficiently small, then $\Xi$ contains only stable equilibria, satisfying $s^* <(\gamma_i + \xi_i + \nu)/ \beta$.
This allows to deduce that solutions initiated within $\Xi$ asymptotically converge to an equilibrium point by using arguments analogous to those in \cite[Prop. 4.7, Thm 4.20]{haddad2011nonlinear}.

\emph{Part (ii):} Consider any equilibrium $x^\ast$ such that $s^\ast > (\gamma_i + \xi_i + \nu)/ \beta$. 
Then, the Jacobian matrix associated with the linearisation of \ak{equation} \eqref{SIDARE} at $x^\ast$ is given by
\begin{equation*}
J_{x^{\ast}} = 
\begin{bmatrix}
0 & -\beta s^\ast & 0 & 0 & 0 & 0 \\
0 & \beta s^\ast - \gamma_i - \xi_i - \nu & 0 & 0 & 0 & 0 \\
0 & \nu & -\gamma_d - \xi_d & 0 & 0 & 0 \\
0 & \xi_i & \xi_d & -\gamma_a - \mu & 0 & 0 \\
0 & \gamma_i & \gamma_d & \gamma_a & 0 & 0 \\
0 & 0 & 0 & \mu & 0 & 0
\end{bmatrix}.
\end{equation*}
It can easily be shown that $J_{x^{\ast}}$ has three zero eigenvalues and three more given by $-\gamma_d - \xi_d$, $-\gamma_a - \mu$ and $\beta s^\ast - \gamma_i - \xi_i - \nu$.
When $s^\ast > (\gamma_i + \xi_i + \nu)/ \beta$, the last eigenvalue is positive which suggests that the linearisation at $x^\ast$ is unstable.
\hfill $\blacksquare$

The value of $\bar{R}_0$ corresponds to the basic reproduction rate.
Stability of the dynamics \ak{described by equation} \eqref{SIDARE} imply a basic reproduction rate of less than or equal to $1$.
Hence, for the considered case it holds that $\bar{R}_0 \bar{s} = 1$.
Therefore, from Proposition \ref{lemma_R0} it follows that the value of $\bar{R}_0$ is given by $\bar{R}_0 = 1/\bar{s} = \beta/(\gamma_i + \xi_i + \nu)$.

It should be noted that all all analytical results presented above, can be trivially extended when the term $\mu a$ in \ak{equations} \eqref{SIDARE_a} and \eqref{SIDARE_e} is replaced with  $\overline{\mu}(a)$.

\ak{
\section{Optimal control methodology}\label{sec_sol_meth}
}

In this section we present the solution approach that yields the optimal continuous strategy and optimized strategies with a limited number of policies and policy changes.

To improve the presentation of the results, we note that the states $a$ and $e$ are uniquely defined by the initial conditions and $u$, and hence a functional $\mathcal{J}$ can be defined such that $\mathcal{J}(u) = J(a,e,u)$ for given initial conditions. Therefore, the optimization problem Problem (\ref{Problem_to_min}) can be equivalently written as
 \begin{equation}\label{Problem_to_min_equiv}
 \begin{aligned}
& \min_{u \in \mathcal{U}} \mathcal{J}(u) \\
& \text{s.t. } \eqref{SIDARE_u}.    
\end{aligned}
 \end{equation}

\subsection{Solution approach}\label{sec_solution}

We call  $\hat{u}$ a solution to Problem (\ref{Problem_to_min}) if $\mathcal{J}(\hat{u}) = \min_{u \in \mathcal{U}} \mathcal{J}(u) = \min_{u \in \mathcal{U}} J(a,e,u)$.
We call the solution $\hat{u}$ and the corresponding state $\hat{x}$ an optimal pair.
The existence of a solution to Problem (\ref{Problem_to_min})  is established by the following lemma.

\begin{lemma}\label{lemma_uniqueness}
There exists a solution to Problem (\ref{Problem_to_min}).
\end{lemma}

\emph{Proof of Lemma \ref{lemma_uniqueness}:}
To prove Lemma \ref{lemma_uniqueness}, we apply  \cite[Thm. 2.1, p. 63]{fleming2012deterministic}, which provides sufficient conditions for the existence of an optimal solution to a general optimal control problem. 
In particular, the conditions in \cite[Thm. 2.1, p. 63]{fleming2012deterministic} are satisfied for Problem (\ref{Problem_to_min}) since:
\begin{itemize}
\item[(i)] \ak{Equation} \eqref{SIDARE_u} is continuously differentiable, 
\item[(ii)] there exist a feasible  solution  to \ak{equation} \eqref{SIDARE_u},
\item[(iii)] $\mathcal{U}$ is a compact set,
\item[(iv)]  $F(t,x) = \{ [-\beta si(1 - u) \;,\; \beta s i (1 - u) - \gamma_i i - \xi_i i - \nu i \;,\; \nu i - \gamma_d d - \xi_d d \;,\; \xi_i i + \xi_d d - \gamma_a a - \overline{\mu} (a) \;,\; \overline{\mu} (a)]^T : u \in \mathcal{U}\}$ is convex for all $(t, s, i, d, a, e) \in [0, T] \times [0,1]^5$.
\end{itemize}

Hence, there exists an optimal solution $\hat{u}$ to Problem (\ref{Problem_to_min}).
 \hfill $\blacksquare$

To obtain the optimal intervention strategy $\hat{u}$, it will be convenient to form the Hamiltonian for Problem (\ref{Problem_to_min}), as below
\begin{equation}
\label{Hamiltonian}
H(x,u,\lambda) =  \frac{1}{2} u^2 + \frac{1}{2} \theta_a a^2 + \lambda^T (f_0(x) +  f_1(x)u)
\end{equation}
where  $\lambda \in \mathbb{R}^5$ is called the co-state of the system and $f_0(x)$ and $f_1(x)$ follow from \ak{equation} \eqref{SIDARE_u} and are given\footnote{Note that the state $r$ within $x$ does not appear in either $f_0(x)$ or $f_1(x)$. We  defined $f_0$ and $f_1$ as functions of $x$ to avoid introducing extra notation. }  by
\begin{equation*}
f_0(x) = \begin{bmatrix}
-\beta s i \\
\beta s i  - \gamma_i i - \xi_i i - \nu i \\
\nu i - \gamma_d d - \xi_d d \\
 \xi_i i + \xi_d d - \gamma_a a - \overline{\mu} (a) \\
 \overline{\mu} (a)
\end{bmatrix}, \;
 f_1(x) = \begin{bmatrix}
\beta s i \\
- \beta s i \\
 0 \\
 0 \\
 0
\end{bmatrix}.
\end{equation*}

Below, we provide necessary optimality conditions for Problem (\ref{Problem_to_min}), which are a result of Pontryagin's minimum principle.

\begin{proposition}\label{Proposition_Pontryagin}
Let $(\hat{x}, \hat{u})$ be a locally optimal pair to Problem (\ref{Problem_to_min}).
Then, there exists a co-state function $\hat{\lambda} :[0,T] \rightarrow \mathbb{R}^5$ such that the following conditions hold for almost all $t \in [0,T]$:
\begin{subequations}\label{Pontryagin_equations}
\begin{equation}
\dot{\hat{\lambda}}^T = - (\hat{\lambda}^T ( \nabla f_0(\hat{x}) + \nabla f_1(\hat{x}) \hat{u} )),
\end{equation}
\begin{equation}
\hat{\lambda}(T) = [0 \; 0 \; 0 \; 0 \; \theta_e]^T,
\end{equation}
\begin{equation}\label{Pontryagin_equations_u}
 \hat{u} =  [- \hat{\lambda}^T f_1(\hat{x})]_\mathcal{U}
\end{equation}
\end{subequations}
\end{proposition}

\emph{Proof of Proposition \ref{Proposition_Pontryagin}:}
The proof follows from applying Pontryagin's minimum principle \cite{pontryagin2018mathematical} in Problem (\ref{Problem_to_min}).
The principle states that for a trajectory  $(\hat{x}, \hat{u})$ that solves Problem (\ref{Problem_to_min}), there exists a function $\hat{\lambda}: [0, T] \rightarrow \mathbb{R}^5$ that that:
\begin{subequations}\label{Pontryagin_equations_proof}
\begin{equation}
\dot{\hat{x}}^T = \frac{\partial H}{\partial \lambda}(\hat{x}(t), \hat{u}(t), \hat{\lambda}(t)), \;  \hat{x}(0) = x_0,
\end{equation}
\begin{equation}
\dot{\hat{\lambda}}^T = - \frac{\partial H}{\partial x}(\hat{x}(t), \hat{u}(t), \hat{\lambda}(t)), 
\end{equation} 
\begin{equation}
\hat{\lambda}(T) = \frac{\partial S}{\partial x} (\hat{x}(T))^T,
\end{equation}
\begin{equation}\label{Pontryagin_equations_proof_u}
H(\hat{x}(t), \hat{u}(t), \hat{\lambda}(t)) \leq H (\hat{x}(t), u(t), \hat{\lambda}(t)),
\quad \forall u \in \mathcal{U},
\end{equation}
where the Hamiltonian $H$ is given by \ak{equation} \eqref{Hamiltonian} and $S$ is associated with the final cost, and is given by 
$S(x) = \theta_e e$. 
Hence, \ak{equation} \eqref{Pontryagin_equations} follows directly\footnote{
Note that there is a technical issue in the definition of $\frac{\partial H}{\partial x}$ when $a = \bar{h}$ since $\frac{d\bar{\mu}(a)}{da}$ is not well defined at this point. 
This issue can be resolved by considering the subdifferntial of $\bar{\mu}(a)$ at $a = \bar{h}$ \cite{rockafellar2015convex}, and define the solutions for $x$ and $\lambda$ in the sense of Filippov \cite{filippov1988differential}.  
We refrain from properly defining these concepts to avoid the introduction of extensive technical notions and to keep the focus of the paper on the practical aspects of the results.
} 
from \ak{equations} \eqref{Hamiltonian} and \eqref{Pontryagin_equations_proof}.
 Note that \ak{equation} \eqref{Pontryagin_equations_u} follows from \ak{equation} \eqref{Pontryagin_equations_proof_u}, which results to  $\frac{\partial H}{\partial u}(\hat{x}(t), \hat{u}(t), \hat{\lambda}(t)) = 0$ when $\hat{u} \in \mathcal{U}^{\mathrm{o}}$, i.e. when $\hat{u}$ lies in the interior of $\mathcal{U}$.
\end{subequations}
 \hfill $\blacksquare$

We make use of \ak{equation} \eqref{Pontryagin_equations} and the \ak{controlled} SIDARE model dynamics, \ak{described by equation} \eqref{SIDARE_u} to obtain the optimal solution to Problem (\ref{Problem_to_min}). 
To tackle the challenges associated with dealing with initial and final value constraints in numerical simulations, we used an adapted forward-backward sweep method (e.g. \cite[Ch. 21]{lenhart2007optimal}).

\subsection{Implementing a limited amount of policies and policy changes}

An implementable government {strategy} requires that $u$ takes a small number of distinct values, i.e. there exists a finite  set of distinct possible intervention policies $\mathcal{U}_d = \{\tilde{u}_1, \tilde{u}_2, \dots, \tilde{u}_{\bar{m}} \} \subset \mathcal{U}$ such that $u(t) \in \mathcal{U}_d$ for all $t \geq 0$.
{For notational convenience, it is assumed that $\tilde{u}_i < \tilde{u}_j$ for $i < j$.} 
Moreover, for each strategy $u$, we consider the set of policy change instants $\mathcal{T} = \{t_1, t_2, \dots t_n\}$ satisfying $0 < t_1 < t_2 < \dots < t_n < T$, 
where $T > 0$ denotes the considered timeframe,
 such that $u(t) = u_j \in \mathcal{U}_d, t \in [t_{j-1}, t_j), {1 \leq j \leq n+1}$, where $t_0 = 0$ and $t_{n+1} = T$.
In addition, it is assumed that $u_{j} \neq u_{j+1}, 1 \leq j \leq n-1$.
Note that the cardinality of the set $\mathcal{T}$ is important as it describes the amount of changes between policies.
For any intervention strategy $u$, we let $\mathcal{R}(u) = \{\tilde{u}: \exists t \in [0,T] \text{ s.t. } u(t) = \tilde{u}\} $ satisfying $\mathcal{R}(u) \subseteq \mathcal{U}_d$ denote the set of {policies within $u$}. 
In addition, $|\mathcal{R}(u)|$ describes the number of distinct policies within {strategy} $u$. 
 
 
Drawing a practical government intervention strategy motivates the introduction of two constraints: (i) on the number of distinct policies, so that a small set of rules is implemented by the population,
(ii) on the number of policy changes, 
since frequent policy changes may result in social fatigue, decreasing the responsiveness of the population to the policy instructions.
To account for these, we have defined the  optimization problem \eqref{Problem_to_min_discrete}, presented also below for convenience:
 \begin{equation*}
 \begin{aligned}
& \min_{u \in \mathcal{U}_d} J(a,e,u) \\
& \text{s.t. } \eqref{SIDARE_u},    |\mathcal{R}(u)| \leq \hat{n}_1, |\mathcal{T}| \leq \hat{n}_2,
\end{aligned}
 \end{equation*}
where $\hat{n}_1$ and $\hat{n}_2$ respectively denote the maximum number of policies that $u$ is allowed to take from the set $\mathcal{U}_d$ and the maximum allowed number of changes in the intervention policy over the considered timeframe.
In analogy to (\ref{Problem_to_min_equiv}), Problem (\ref{Problem_to_min_discrete}) can be equivalently written as
 \begin{equation}
 \begin{aligned}
& \min_{u \in \mathcal{U}_d} \mathcal{J}(u) \\
& \text{s.t. } \eqref{SIDARE_u},   |\mathcal{R}(u)| \leq \hat{n}_1, |\mathcal{T}| \leq \hat{n}_2.
\end{aligned}
 \end{equation}

To solve the above problem, we make use of the solution to Problem (\ref{Problem_to_min}), which is a continuous relaxation to Problem (\ref{Problem_to_min_discrete}), and hence provides a lower bound to the cost of Problem (\ref{Problem_to_min_discrete}).
Our approach to obtain a solution to Problem (\ref{Problem_to_min_discrete}) is described below: \\[2mm]
\emph{(i)} Obtain $\hat{u}$ that solves Problem (\ref{Problem_to_min}).\\[2mm]
\emph{(ii)} Using $\hat{u}$, obtain an optimized solution to Problem (\ref{Problem_to_min_discrete}), denoted by $\hat{u}_d $, resulting from a local minimum cost search algorithm.

Part (i) of the above approach is explained in Section \ref{sec_solution} above (see also Proposition \ref{Proposition_Pontryagin}).
Part (ii) is described in detail below.

\emph{Part (ii):}
The methodology in this section is split in two steps.
First, we obtain a strategy that approximates $\hat{u}$ and simultaneously  satisfies the constraints of Problem (\ref{Problem_to_min_discrete}).
Using this strategy for initialization, we implement Algorithm \ref{algorithm_1} which produces the optimized policy $\hat{u}_d$.
\ak{Note the in the presented process, it is assumed that $\hat{n}_1 > 1$ and $\hat{n}_2 \geq 1$. The case where $\hat{n}_1 = 1$ which necessarily results in $|\mathcal{T}| = \emptyset$ can be trivially solved by evaluating the costs of implementing each value of $u \in \mathcal{U}_d$.}

\textbf{Initialization:}
The first step  aims to obtain an intervention strategy $\bar{u}_d$ that approximates the continuous strategy $\hat{u}$ and simultaneously satisfies the constraints of the problem.
To obtain such strategy, we let $u_{max} = \max_{t \in [0,T]} \hat{u}$ and $u_{min} = \min_{t \in [0,T]} \hat{u}$ and define $\hat{u}_i = [u_{min} + (i-1)*(u_{max} - u_{min})/(\hat{n}_1 - 1)]_{\mathcal{U}_d}$ and $U_s = \{\hat{u}_i: i \in \{1, \dots, \hat{n}_1\}\}$, where $\hat{n}_1$ is given in Problem (\ref{Problem_to_min_discrete}).
Consequently, we let $v_{d}(t) = \argmin_{y \in U_s} |y - \hat{u}(t)|$  for all $t \in [0, T]$, i.e. we project the values of the continuous policy $\hat{u}$ onto the discrete set $U_s$.
The latter leads to $\mathcal{R}(v_{d}) = U_s$.
Moreover, we define the set $\bar{\mathcal{T}}^i = \{t \in [0, T] : \lim_{\epsilon \rightarrow 0} v_d(t + \epsilon) \neq \lim_{\epsilon \rightarrow 0} v_d(t - \epsilon)\}$, i.e. the set of all time instants when a change occurs in the intervention policy $v_d$.

If $|\bar{\mathcal{T}}^i| \leq \hat{n}_2$, then we  select $\bar{u}_d = v_d$.
Otherwise,  we define 
{$|\bar{\mathcal{T}}^i|$} strategies $v_{d,j}$, such that strategy $v_{d,j}$ satisfies $v_{d,j}(t) = v_d(t), t \in [0, T]\setminus [t_{j}, t_{j+1})$ and  $v_{d,j}(t) = v_d(t_{j-1}), t \in  [t_{j}, t_{j+1})$, i.e. $v_{d,j}$ follows by omitting the $j$th switch in $v_d$.
We then calculate the cost  for each of the above policies, given by $\mathcal{J}(v_{d,j})$,
 and select the $q = |\bar{\mathcal{T}}^i| - \hat{n}_2$ {strategies} with the least cost.
 Finally, we {neglect} the switches associated with these $q$ {strategies} from $v_d$, i.e. the $q$ switches which individually cause the least increase in cost once removed, to construct $\bar{u}_d$.
The latter enables the construction of $\bar{u}_d$ which satisfies all the constraints of Problem (\ref{Problem_to_min_discrete}).

\textbf{Algorithm \ref{algorithm_1}:}
After obtaining $\bar{u}_d$, as explained above, then an optimized solution $\hat{u}_d$ is obtained from Algorithm \ref{algorithm_1}, as described below.
To initialise Algorithm \ref{algorithm_1}, we first consider the set of switching instants  for strategy $\bar{u}_d$, which we denote by  $\bar{\mathcal{T}}$, and select a sufficiently small $\delta > 0$ which corresponds to the tolerance of Algorithm \ref{algorithm_2}, implemented within Algorithm \ref{algorithm_1}.
In addition, the vector $p$, with $p \in \mathbb{N}^{|\mathcal{R}(\hat{u}_d)|}$, represents the indices of implemented policies from $\mathcal{U}_d$ in strategy $\bar{u}_d$ in decreasing order, i.e. $p_i = k$ means that the $i$th largest implemented policy is $\tilde{u}_k$, reminding that $\mathcal{U}_d = \{\tilde{u}_1, \tilde{u}_2, \dots, \tilde{u}_{\bar{m}}\}$.

Algorithm \ref{algorithm_1} aims to obtain an optimized intervention strategy, denoted by $\hat{u}_d$, 
 by making incremental changes in the test policy $\hat{u}_t$, as seen in \ak{equations} \eqref{alg_1}--\eqref{alg_2}, e.g. 
when $m=1$ and $n=1$ then Algorithm \ref{algorithm_1} replaces the largest policy within $\hat{u}_t$ with the immediately more relaxed policy within $\mathcal{U}_d$, in analogy when $n=2$ then the attempted policy is the  one which is the immediately more strict. 
Then, Algorithm \ref{algorithm_2}, described in detail below, makes use of the test policy $\hat{u}_t$ and set of switching times $\mathcal{T}$ and provides a new test policy with an optimized set of switching times, denoted by $\tilde{\mathcal{T}}$, and a corresponding cost, denoted by $C$.
If $C$ is less than the lowest cost obtained by that stage, denoted by $C_m$, then $\hat{u}_d$ is updated as in \ak{equation} \eqref{alg_41}.
In addition, the supplementary variables $C_m, p$ and $\bar{\mathcal{T}}$ are updated to facilitate the progression of the algorithm, as described in \ak{equations} \eqref{alg_42}--\eqref{alg_44}.

The convergence variable $\theta$ takes the value of $0$, which {imposes that} a new set of  iterations need to be performed for the convergence of the algorithm. 
If there is no improvement in $C_m$ for all $m \in \{1, \dots, |\mathcal{R}(\hat{u}_d)|\}$ and $n \in \{1,2\}$, then Algorithm \ref{algorithm_1} converges.
Convergence of Algorithm \ref{algorithm_1} suggests that any incremental change in any of the policies within $\hat{u}_t$ will result in a higher cost, and hence that a local minimum has been reached.
For improved clarity, a flowchart of Algorithm \ref{algorithm_1} is presented in Fig. \ref{Flowchart_1}.

\textbf{Algorithm \ref{algorithm_2}:}
Algorithm \ref{algorithm_2} aims to obtain an optimized set of switching times, denoted by $\tilde{\mathcal{T}}$ and a corresponding  intervention strategy $\hat{u}_t$ by making use of the {discrete test} policy $\hat{u}_t$, the set of switching times ${\mathcal{T}}$ and the tolerance level $\delta$ provided from Algorithm \ref{algorithm_1} for its initialization.
In particular, it alters the switching times by $\delta$ in either direction, as demonstrated in \ak{equation} \eqref{alg_21} and updates $\bar{u}_t$ as described in \ak{equations} \eqref{alg_22}--\eqref{alg_23}.
Then, a new cost is calculated, based on the new strategy $\bar{u}_t$ by implementing \ak{equation} \eqref{alg_24}.
If the new cost, denoted by $\tilde{C}$ is lower than the cost obtained by policy $\hat{u}_t$, then the cost, the policy and the set of switching times are updated, as demonstrated in \ak{equations} \eqref{alg_251}--\eqref{alg_253}. Otherwise, the policy $\bar{u}_t$ reverts to its previous value.

The variable $\phi$ serves as a convergence variable. In particular, when any iteration occurs that allows a decrease in cost then its value is set to $0$, which results in a new set of iterations for $\tilde{\mathcal{T}}$.
Convergence of the algorithm suggests that a local minimum of the cost $C$ is reached, where changing any switching time in $\tilde{\mathcal{T}}$ does not result in a lower cost.
For improved intuition,  Fig. \ref{Flowchart_2} depicts a flowchart of Algorithm \ref{algorithm_2}.

Algorithm \ref{algorithm_1} creates a monotonically decreasing sequence of values for $C$, which is lower bounded by $C^* = \mathcal{J}(\hat{u})$, i.e. the cost associated with the optimal continuous policy obtained in part (i).
Therefore, the sequence  of updates in $C$  and  hence Algorithm \ref{algorithm_1} converge, as directly follows by the monotone convergence theorem (e.g. \cite[Theorem 2.4.2]{abbott2001understanding}).
A similar argument follows for the convergence of Algorithm \ref{algorithm_2}.

%
%

\begin{figure}[H]
\centering
\includegraphics[scale=0.6]{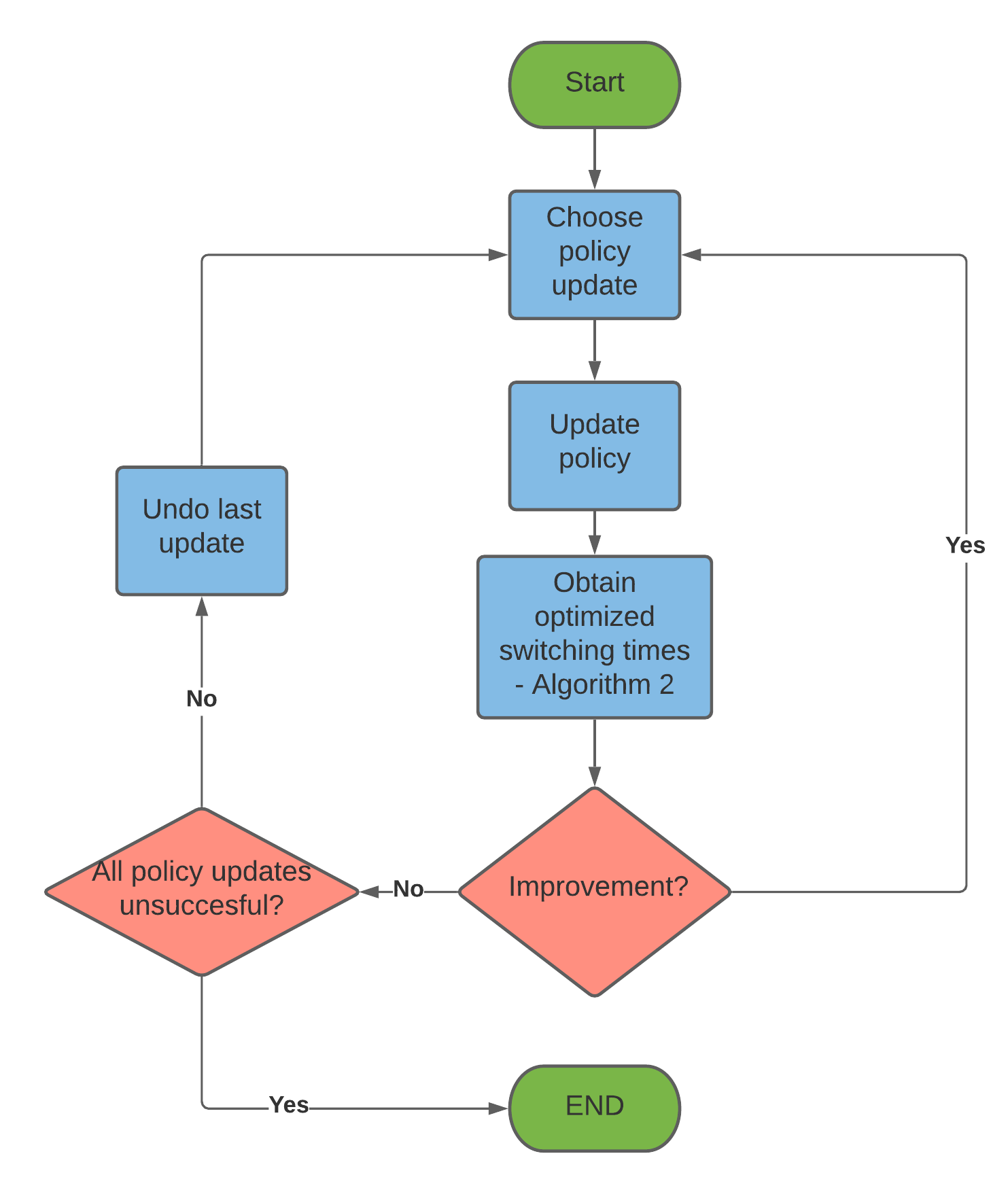}
\caption{Flowchart for Algorithm \ref{algorithm_1}.}
\label{Flowchart_1}
\end{figure}

{\centering
\begin{minipage}{.55\linewidth}
\begin{algorithm}[H]\label{algorithm_1}
\textbf{Inputs:} $\bar{\mathcal{T}}, \bar{u}_d, \delta, \mathcal{U}_d, p$.

\textbf{Output:} $\hat{u}_d $.

 \textbf{Initialization:}
$C = \mathcal{J}(\bar{u}_d),C_m = C, \mathcal{T} = \bar{\mathcal{T}}, \hat{u}_d = \hat{u}_t = \bar{u}_d, \theta = 0.$ \\
\While{$\theta = 0$}
{
$\theta = 1$, \\
\textbf{for} $m = 1:|\mathcal{R}(\hat{u}_d)|$ \\
\text{ \quad }  
\textbf{for} $n = 1:2$ 
\begin{align*}
&  \text{ \quad }\hat{\mathcal{T}} = \{t: \hat{u}_d(t) = \tilde{u}_{p_m}\}, 
\tag{A.1.1}  \hspace{20mm}  \label{alg_1}\\
& \text{ \quad }\hat{u}_t(t) = \tilde{u}_{[p_m + (-1)^n]_{[1, \bar{m}]}}, t \in \hat{T}, \tag{A.1.2}  \label{alg_2} \\
&\text{ \quad } [\tilde{\mathcal{T}}, \hat{u}_t, C] = \text{Algorithm 2}({\mathcal{T}}, \hat{u}_t, \delta) \tag{A.1.3}  \label{alg_3} \\
& \text{ \quad } \textbf{if }
 C < C_m,  \\
& \text{ \quad }\hspace{4mm}  \hat{u}_d = \hat{u}_t, \label{alg_41} \tag{A.1.4.1} \\
&\text{ \quad } \hspace{4mm}  C_m = C,  \label{alg_42} \tag{A.1.4.2} \\
& \text{ \quad } \hspace{4mm}  p_m = [p_m + (-1)^n]_{[1, \bar{m}]}, 
\label{alg_43} \tag{A.1.4.3}\\
& \text{ \quad } \hspace{4mm}  {\mathcal{T}} = \tilde{\mathcal{T}}, 
\label{alg_44} \tag{A.1.4.4}\\
& \text{ \quad }\hspace{4mm}  \theta = 0, \label{alg_45} \tag{A.1.4.5} \\
&  \text{ \quad }\textbf{else }  \nonumber \\
& \text{ \quad }\hspace{4mm} \hat{u}_t = \hat{u}_d, \label{alg_46} \tag{A.1.4.6} \\
& \text{ \quad }\textbf{end}
\end{align*}
    \text{ \quad }  \textbf{end}  \\
    \textbf{end} 
}
\caption{Scheme to obtain $\hat{u}_d $.}
\end{algorithm}
\end{minipage}
\par
}

\begin{figure}[H]
\centering
\includegraphics[scale=0.65]{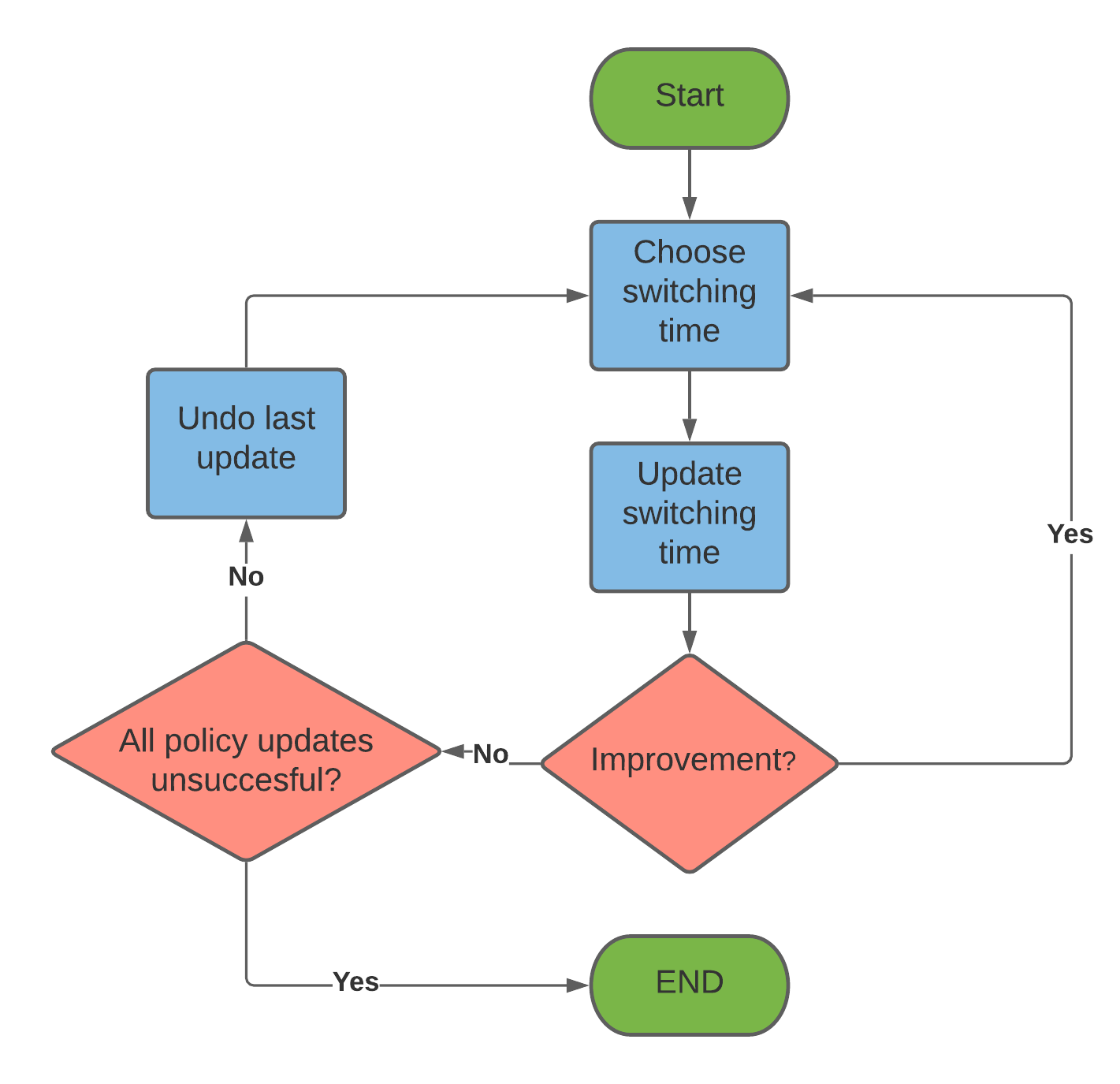}
\caption{Flowchart for Algorithm \ref{algorithm_2}.}
\label{Flowchart_2}
\end{figure}

{\centering
\begin{minipage}{.55\linewidth}
\begin{algorithm}[H]\label{algorithm_2}
\textbf{Inputs:} ${\mathcal{T}}, \hat{u}_t, \delta$.

\textbf{Output:} $\tilde{\mathcal{T}}_i, {1\leq i \leq |\mathcal{T}|},\hat{u}_t, C$.

 \textbf{Initialization:}
$C = \mathcal{J}(\hat{u}_t),  \bar{u}_t = \hat{u}_t, \tilde{\mathcal{T}} = {\mathcal{T}},  \phi = 0, {\tilde{\mathcal{T}}_0 = 0, \tilde{\mathcal{T}}_{|\mathcal{T}| + 1} = T}.$ \\

\While{$\phi = 0$}
{
$\phi = 1,$ 
\\ 
\textbf{for} $k = 1:|{\mathcal{T}}|$ \\
\text{ \quad } \textbf{for} $j = 1:2$
 \begin{align*}
   \hat{\delta} =  \hspace{2mm} & \delta \times (-1)^j,
      \tag{A.2.1}
   \label{alg_21}
   \\
      \tilde{\mathcal{S}} = &
      \begin{cases}
      [[\tilde{\mathcal{T}}_k + \hat{\delta}]_{[\tilde{\mathcal{T}}_{k-1}, \tilde{\mathcal{T}}_k]}, \tilde{\mathcal{T}}_k], \text{ if } j = 1, \\
      [\tilde{\mathcal{T}}_k, [\tilde{\mathcal{T}}_k + \hat{\delta}]_{[\tilde{\mathcal{T}}_k, \tilde{\mathcal{T}}_{k+1}]}], \text{ if } j = 2,
      \end{cases}
      \tag{A.2.2}
   \label{alg_22}
   \\
\bar{u}_t(t)& = {\hat{u}_t(\tilde{\mathcal{T}}_{k + 1 - j}),} \forall t \in \tilde{\mathcal{S}},
   \tag{A.2.3} 
   \label{alg_23}
   \\
    \tilde{C}  =  \hspace{1mm} & \mathcal{J}(\bar{u}_t),
      \tag{A.2.4}
   \label{alg_24} \\
    \textbf{if } \tilde{C} & < C \\
   & C = \tilde{C}, \label{alg_251} \tag{A.2.5.1} \\
   & \hat{u}_t = \bar{u}_t, \label{alg_252} \tag{A.2.5.2}\\
   & \phi = 0, \label{alg_253} \tag{A.2.5.3} \\
        & \tilde{\mathcal{T}}_{k} =  \hspace{1mm}  [\tilde{\mathcal{T}}_{k} + \hat{\delta}]_{{[\tilde{\mathcal{T}}_{k-1},\tilde{\mathcal{T}}_{k+1}]}},
    \tag{A.2.5.4}
   \label{alg_254}
   \\
   \textbf{else}& \\
   &\bar{u}_t = \hat{u}_t, \label{alg_255} \tag{A.2.5.5}\\
     \textbf{end} &
    \end{align*}
    \text{ \quad } \textbf{end}  \\
    \textbf{end} 
}
\caption{Scheme to obtain $\tilde{\mathcal{T}}, \hat{u}_t, C$.}
\end{algorithm}
\end{minipage}
\par
}

\section{Supplementary Results}

In this section,  we present additional results that supplement the presented findings.
In particular, we present optimized strategies and the corresponding decease rates when $4$, $7$ and $10$ distinct policies and $6$, $12$ and $18$ policy changes are allowed.
Furthermore, we complement our results associated with the effect of {parametric} uncertainty, depicted in Fig. \ref{Heat_All}, with additional results demonstrating the optimal strategies and corresponding \ak{aggregate deceases} for a range of values for  the infection fatality rate and initial basic reproduction number $\bar{R}_0$. 
These results aim to provide additional intuition on the effect of parametric uncertainty in forming efficient government strategies.

\subsection{Simulated parameters}

The controlled SIDARE model, \ak{described by equation} \eqref{SIDARE_u} has been used in our simulations, with a time horizon of $T = 365$ days.
The selected initial conditions correspond to the very early stage of the disease, where $0.001\%$ of the population has been infected and there are no detected cases, acutely symptomatic cases, fatalities or recoveries yet.

The values of $\gamma_i$ and $\gamma_d$ were selected following \cite{who_report}, which  suggests a median time of disease onset to recovery for mild cases of approximately two weeks.
The value of $\gamma_a$ was selected following \cite{wang2020phase}, which reported a median time between hospitalization and recovery of $12.4$ days.

The value of $\beta$ was selected to be $0.251$, corresponding to an  initial basic reproduction number $\overline{R}_0 = \beta s_0/ (\gamma_i + \xi_i + \nu)$ of approximately $3.27$ following \cite{yuan2020monitoring}, and assuming $\nu = 0$ at $t = 0$ days.
The value of $\bar{u}$, corresponding to the maximum allowed value for the input $u$ was selected to be $0.8$.
Furthermore, to select the values for $\xi_i$ and $\xi_d$ for the case of Italy, 
we used the findings from \cite{verity2020estimates} on hospitalization rate per age group and  data for the Italian population age distribution  \cite{Italy_UN}.
In addition, we considered $333$ care beds per $100,000$ habitants following \cite{rhodes2012variability}.

The value of $\mu$  was selected to be $0.0085$, being associated with an infection fatality rate of $0.66\%$ as reported  in \cite{verity2020estimates}. The latter is in agreement with various studies that report a mortality rate close or below $1\%$ \cite{mallapaty2020deadly}, \cite{hallal2020remarkable}, \cite{salje2020estimating}.
We let $\hat{\mu}$, which corresponds to the fatality rate when the healthcare system capacity is overloaded, be $5$ times higher than $\mu$, motivated from \cite{catena2020relocate} which compared the fatality rates in two regions in Italy, Lombardy where more than $80\%$ of the healthcare capacity was held by COVID-19 patients and Veneto where up to $40\%$ was held, and deduced an approximate five fold increase in the mortality rate. 

To facilitate the reproducibility of the results, the parameter values used in the simulations, and the justification for their selection, are reported in Table \ref{parameters_table}.

 \begin{table}[h!]
\centering
\caption{Model parameter values}
\begin{tabular}{|l|c|c|}
\hline
Parameter & Value & Justification  \\
\hline
 $(s_0, i_0, d_0, a_0, e_0)  $ & $(1 - 10^{-5}, 10^{-5}, 0, 0, 0)$  & Early stage of the pandemic \\  \hline
$\gamma_i, \gamma_d$ & $1/14 = 0.071$ & \cite{who_report} \\ \hline
$\gamma_a$ & $1/12.4$ &  \cite{wang2020phase} \\ \hline
$\beta$ & $0.251$ & $\overline{R}_0$ from \cite{yuan2020monitoring}\\ \hline
$\xi_i, \xi_d$ & 0.0053 
 & \cite{verity2020estimates},\cite{Italy_UN} \\ \hline
$\mu$ & 0.0085 & \cite{verity2020estimates} \\ \hline
 $\overline{h}$ & $0.0481$ & \cite{rhodes2012variability}\\ \hline
 $\hat{\mu}$ & 5$\times \mu$ & \cite{catena2020relocate}\\ \hline
\end{tabular}
 \label{parameters_table}
 \vspace{-2mm}
\end{table}

\subsection{Deceased population versus cost of government intervention} 

In \an{Fig. \ref{Cost_Dth_All}}, we associate the cost of (optimal) government intervention and the portion of deceased population when: (i) different hospital capacity rates, (ii) different testing rates and (iii) different cost \an{emphasis levels associated with} the acutely symptomatic \an{population}, are considered.
{Note that in each case we considered a broad range of cost weight values  associated to the total number of deaths.}
{The considered cases for hospital capacity rates, testing rates, and cost weights attributed to the acutely symptomatic cases and the total number of deaths are presented in Table \ref{case_table}.

 \begin{table}[h!]
\centering
\caption[Simulation parameters]{{Parameters associated with the cases presented in Fig.  \ref{Cost_Dth_All}.}}
\begin{tabular}{|l|c|c|}
\hline
Parameter & Symbol & Cases  \\
\hline
 Hospital capacity rate & $\overline{h}$ & $\{222, 333, 444\} \times 10^{-5}$ \\  \hline
Testing rate & $\nu$ & $\{0, 0.05, 0.10\}$ \\ \hline
Cost weight for acutely symptomatic cases & $\theta_a$ & $\{0, 5\times10^4, 10^5\}$ \\ \hline
Cost weight for aggregate deaths & $\theta_e$ & $[0, 2.5\times10^4]$ \\ \hline
\end{tabular}
 \label{case_table}
\end{table}
}

\ak{The parameters from Table \ref{case_table} associated with each of the $8$ intervention strategies presented in Fig. \ref{Policy_All} are provided in Table \ref{table_parameters} below. In all cases, the considered hospital capacity rate was $333$ care beds per $100,000$ habitants which corresponded to the full capacity case.}

\ak{ 
 \begin{table}[h!]
\centering
\caption[Strategy parameters]{{Parameters associated with the optimal continuous strategies presented in  Fig. \ref{Policy_All}.}}
\begin{tabular}{|c|c|c|c|}
\hline
Tolerance ($\%$) & Testing rate   &  Cost weight for acutely  & Cost weight for  \\
& &    symptomatic cases & aggregate deaths \\ \hline
1 & 0 &  0 &  1600 \\ \hline
1 & 0.05 &  0 & 400 \\ \hline
0.1 & 0 &  $10^5$ &  600 \\ \hline
0.1 & 0.05 &  $10^5$ & 1000  \\ \hline
0.1 & 0.10 &  $5 \times 10^4$ & 1000  \\ \hline
0.01 & 0 &  0 &  $2.5 \times 10^4$ \\ \hline
0.01 & 0.05 &  0 & $1.8 \times 10^4$  \\ \hline
0.01 & 0.10 &  0 & $10^4$  \\ \hline
\end{tabular}
 \label{table_parameters}
\end{table}
}

\subsection{Implementing a limited {number} of policies and policy changes}

In Fig. \ref{Policy_decease_1}--\ref{Policy_decease_8}, we present the intervention strategies and corresponding \ak{aggregate deceases} associated with decease tolerances of $1\%$, $0.1\%$ and $0.01\%$. 
{In particular, we present the optimal continuous strategy and} optimized intervention strategies with $4, 7$ and $10$ distinct policy levels and $6, 12$ and $18$ policy changes respectively. 
These results are associated with Fig. \ref{Policy_All}, which however includes only the case where $4$ distinct policies and $6$  policy changes are allowed.
The results presented in Fig. \ref{Policy_decease_1}--\ref{Policy_decease_8} demonstrate that increasing the number of distinct policies {enables} a decease response that is closer to the optimal {continuous strategy}, which is in agreement with intuition.

\subsection{Effect of {parametric} uncertainty}

In Fig. \ref{robust_R_1}--\ref{robust_R_16}, we present {additional} results associated with {the effect of} uncertainty in knowledge of the initial reproduction number $\bar{R}_0$ and the infection {fatality} rate.
{These results aim to provide additional intuition and clarity on the presented results.}

In particular, we first consider the {case where} the value of $\bar{R}_0$ ranges between $3.17$ and $3.38$, which corresponds to {a} $95\%$ confidence interval, as reported in \cite{yuan2020monitoring} when the infection mortality rate is $0.66\%$, as reported in \cite{verity2020estimates}.
{For the considered case, we obtain the optimal strategies which correspond to each value of $\bar{R}_0$ for the parameters presented in Table \ref{table_parameters}, associated with the $8$ continuous strategies 
\ak{presented in Fig. \ref{Policy_All}.
The ranges of the obtained optimal intervention strategies for each of these $8$ considered cases are}
 depicted in Fig. \ref{robust_R_1}--\ref{robust_R_8} (left).
From these figures, it follows that the uncertainty in $\bar{R}_0$ results in small variations in the optimal strategies when no and slow testing policies are adopted. By contrast, when fast testing policies are implemented, the impact of the uncertainty in $\bar{R}_0$ on the optimal strategies is more substantial.}

{Moreover,  Fig. \ref{robust_R_1}--\ref{robust_R_8} (right)} demonstrate the ranges of decease rates obtained when the same $8$ selected strategies  are applied for  $\bar{R}_0 \in [3.17, 3.38]$ {when the infection mortality rate is $0.66\%$}, i.e. they depict the decease rates from implementing the selected {strategies when the value of} $\bar{R}_0$ has been incorrectly estimated.
The effect of implementing the $8$ selected strategies when in addition the infection mortality rate has been incorrectly estimated is demonstrated  in Fig. \ref{robust_R_9}--\ref{robust_R_16}, which depict the decease rates for the considered range of values for $\bar{R}_0$ when the infection mortality rate is $0.39\%$ and $1.33\%$ respectively.
These values represent the lower and upper bounds of a $95\%$ infection mortality rate confidence interval, as reported in \cite{verity2020estimates}.
{Figures \ref{robust_R_1}--\ref{robust_R_8} (right) and Fig. \ref{robust_R_9}--\ref{robust_R_16} demonstrate that the level of the adopted decease tolerance is crucial when it comes to the effect of model uncertainty in the decease rates. The latter is in agreement with the presented results, which state that parametric uncertainty has a more substantial effect when stricter government strategies, {associated with lower decease tolerances,} are adopted.}

\begin{figure}[H]
\begin{subfigure}{0.5\columnwidth}
\centering
\includegraphics[scale=0.55]{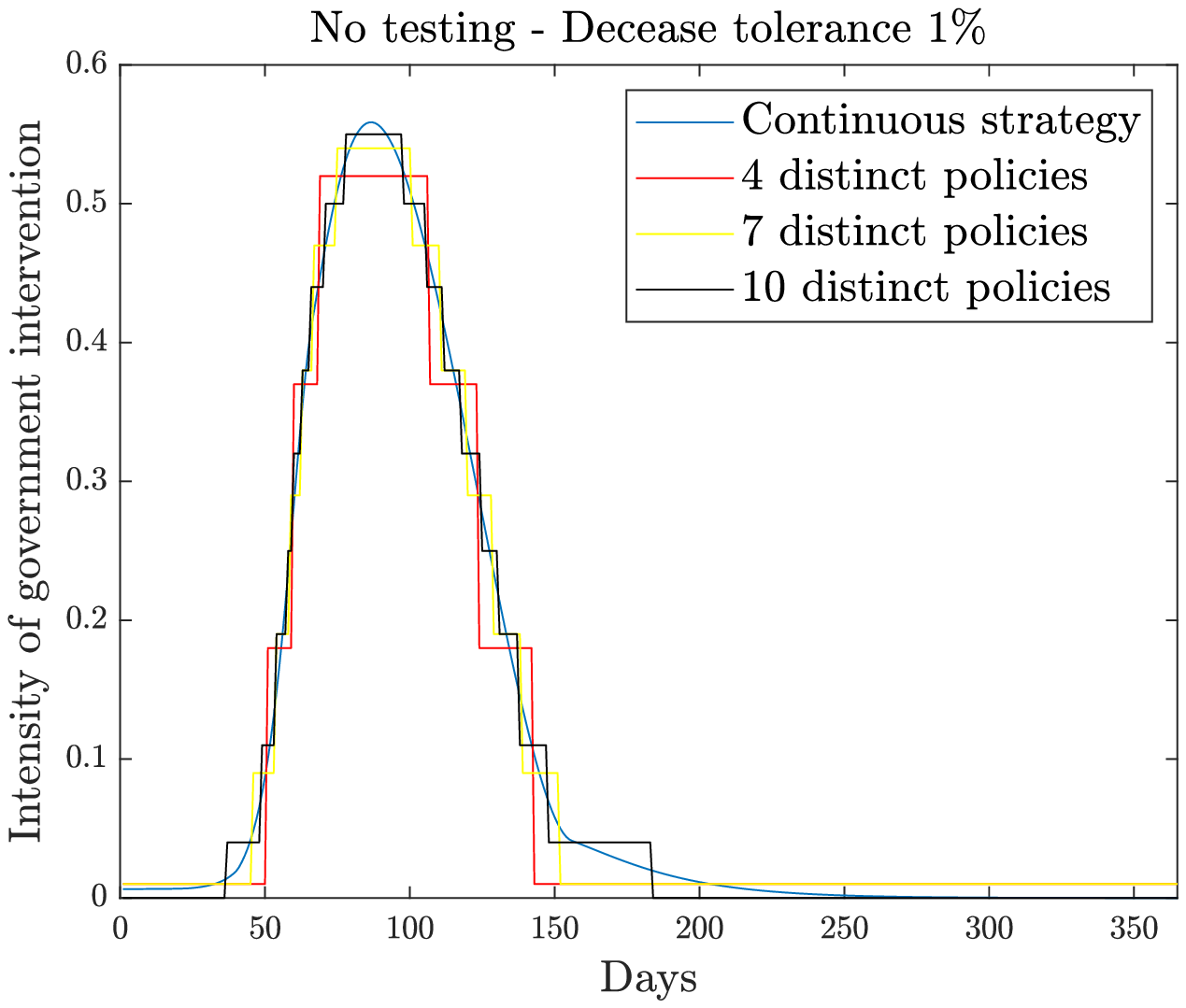}
\end{subfigure} 
\begin{subfigure}{0.5\columnwidth}
\centering
\includegraphics[scale=0.55]{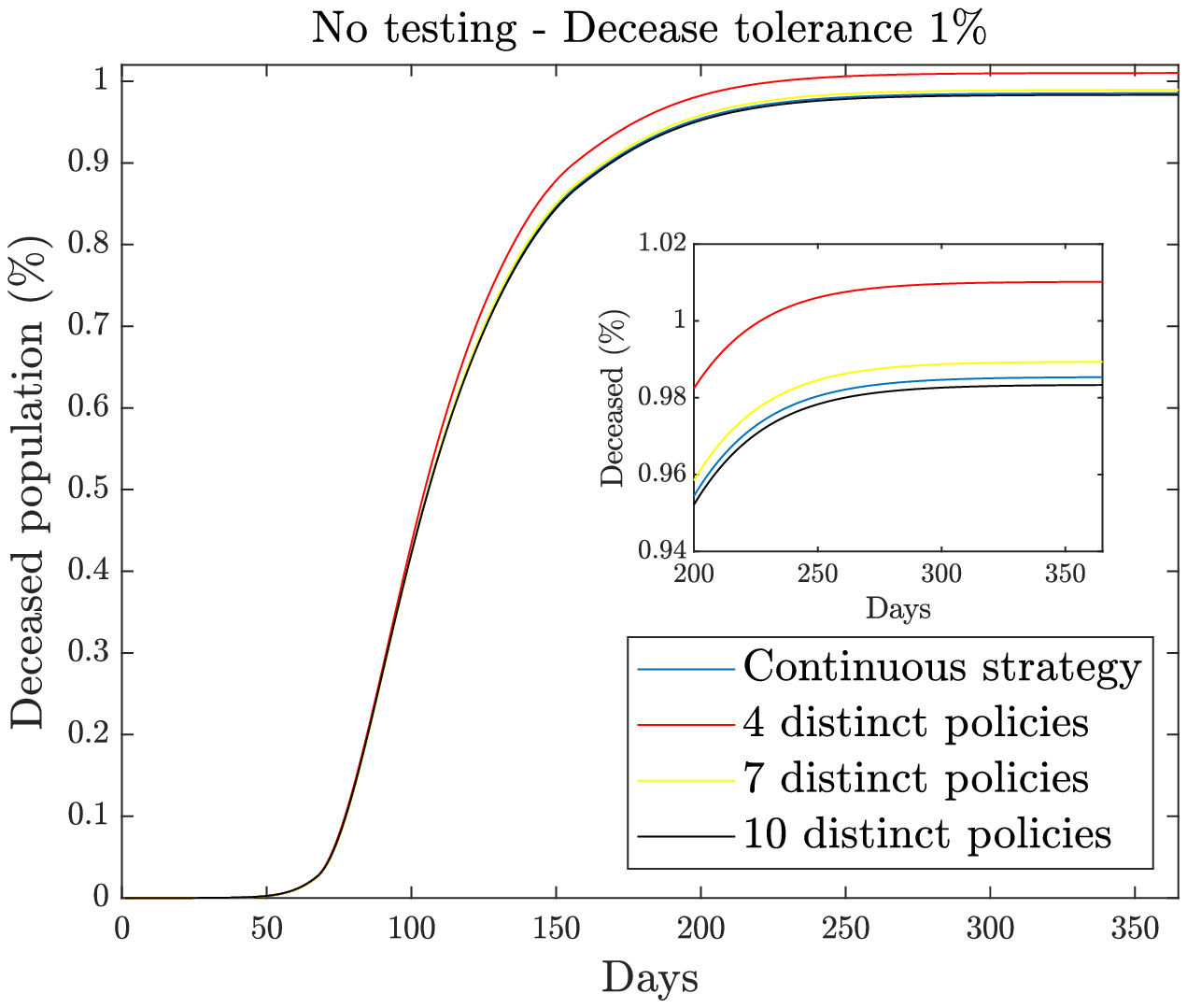}
\end{subfigure}
\vspace{-2mm} \caption[Intervention strategy and decease rate, $1\%$ decease tolerance, no testing]{\textbf{Intervention strategy and decease rate, $1\%$ decease tolerance, no testing.} Intensity of government intervention {(left)} and portion of deceased population {(right)}
when no testing is performed and a decease tolerance of $1\%$ is adopted
  when (i) a continuously changing strategy is selected and (ii) discrete implementations of the selected strategy are considered, allowing $4, 7$ and $10$ policy levels and $6, 12$ and $18$ policy changes respectively. }\vspace{-2mm}
\label{Policy_decease_1}
\end{figure}

\begin{figure}[H]
\begin{subfigure}{0.5\columnwidth}
\centering
\includegraphics[scale=0.55]{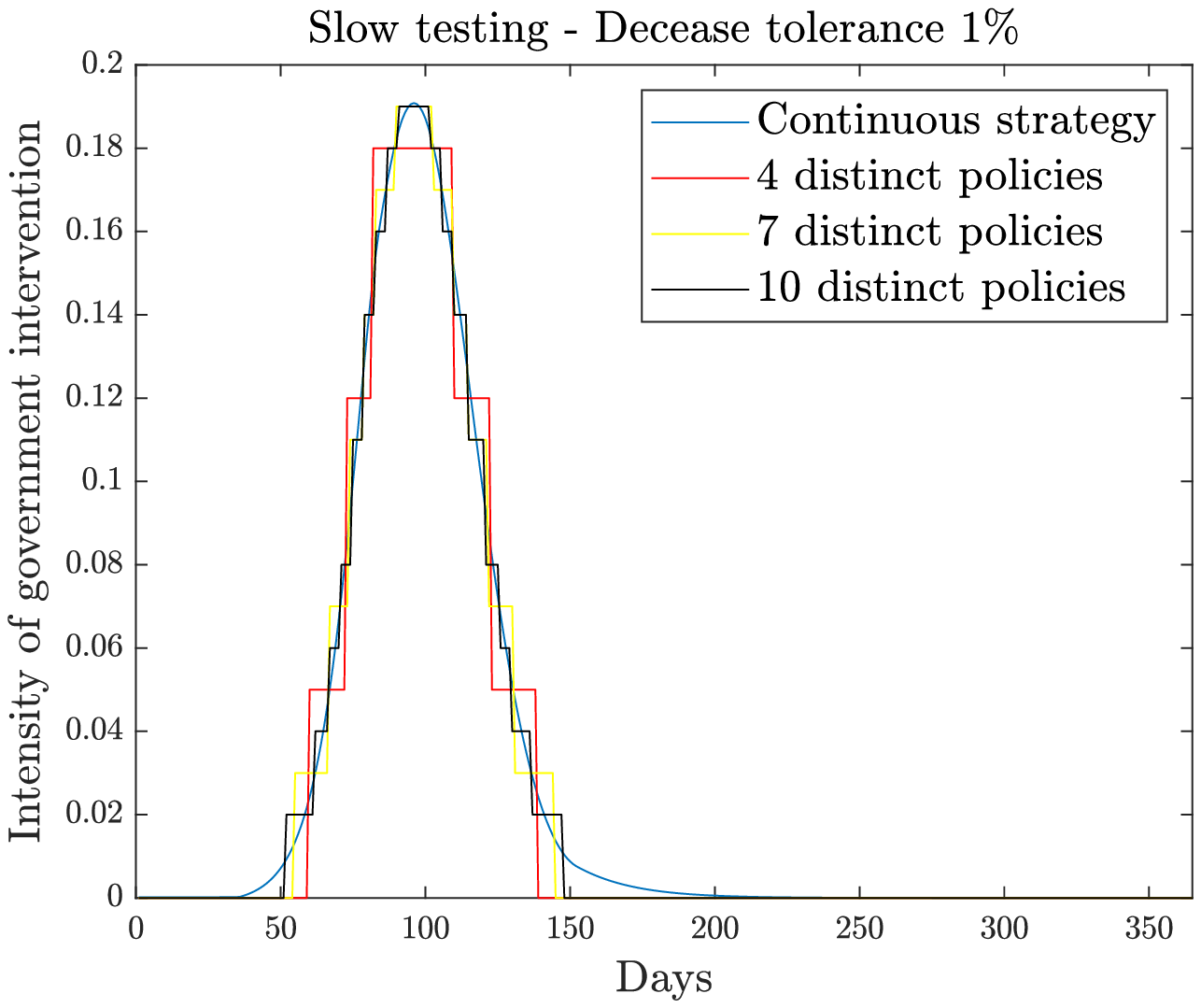}
\end{subfigure} 
\begin{subfigure}{0.5\columnwidth}
\centering
\includegraphics[scale=0.55]{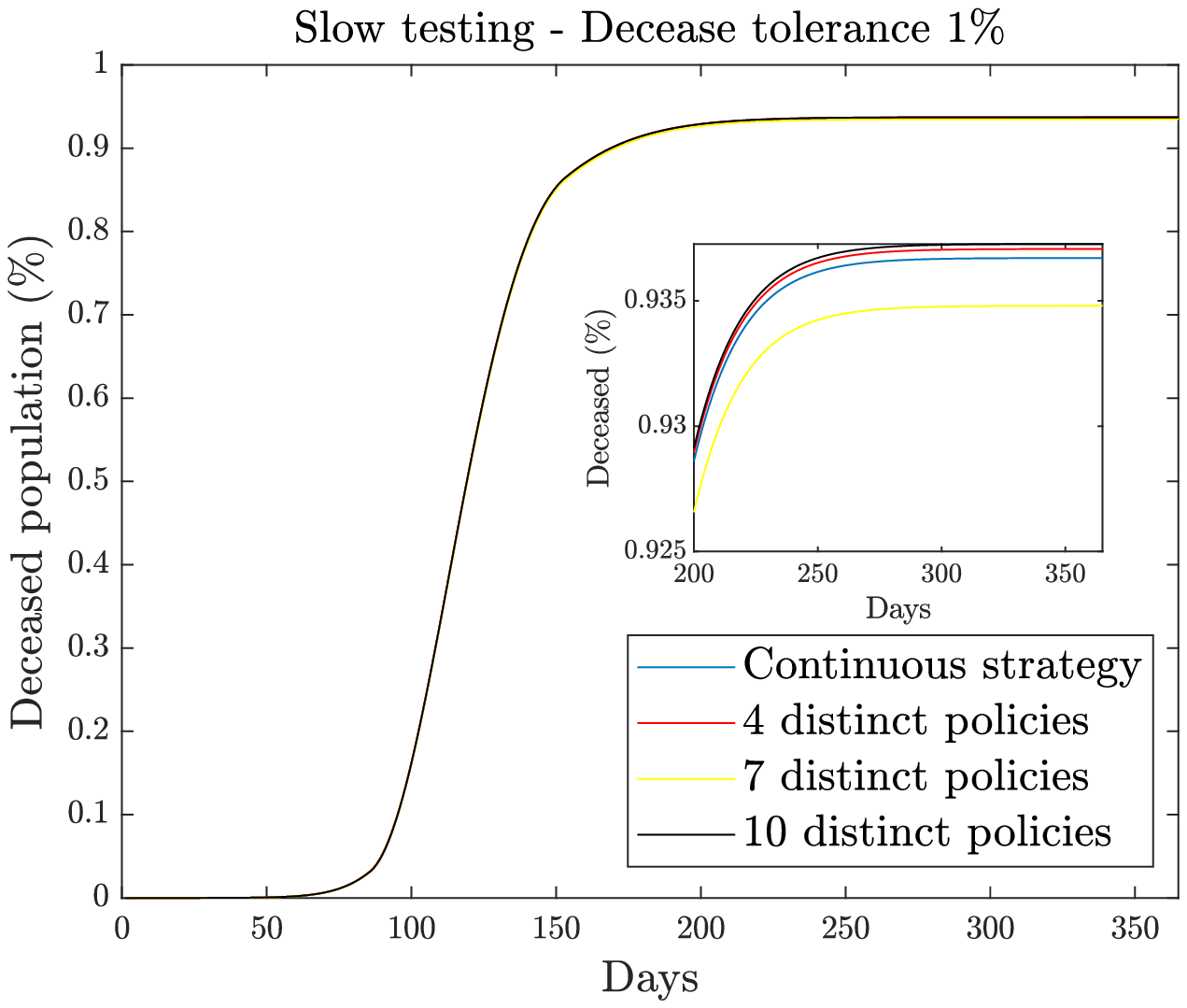}
\end{subfigure}
\vspace{-2mm} \caption[Intervention strategy and decease rate, $1\%$ decease tolerance, slow testing]{\textbf{Intervention strategy and decease rate, $1\%$ decease tolerance, slow testing.} Intensity of government intervention {(left)} and portion of deceased population {(right)}
when a slow testing policy and a decease tolerance of $1\%$ are adopted
  when (i) a continuously changing strategy is select01ed and (ii) discrete implementations of the selected strategy are considered, allowing $4, 7$ and $10$ policy levels and $6, 12$ and $18$ policy changes respectively. }\vspace{-2mm}
\end{figure}

\begin{figure}[H]
\begin{subfigure}{0.5\columnwidth}
\centering
\includegraphics[scale=0.55]{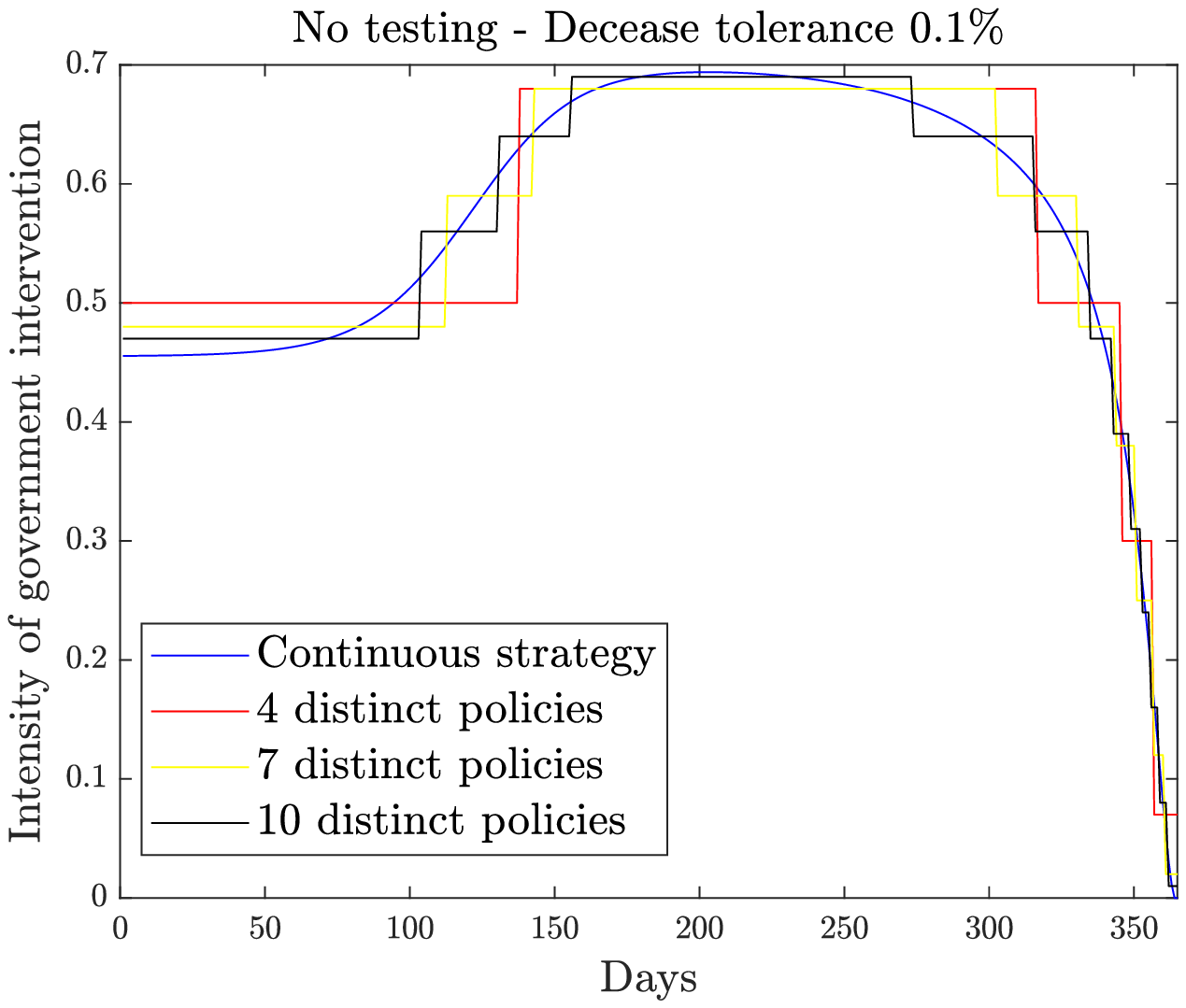}
\end{subfigure} 
\begin{subfigure}{0.5\columnwidth}
\centering
\includegraphics[scale=0.55]{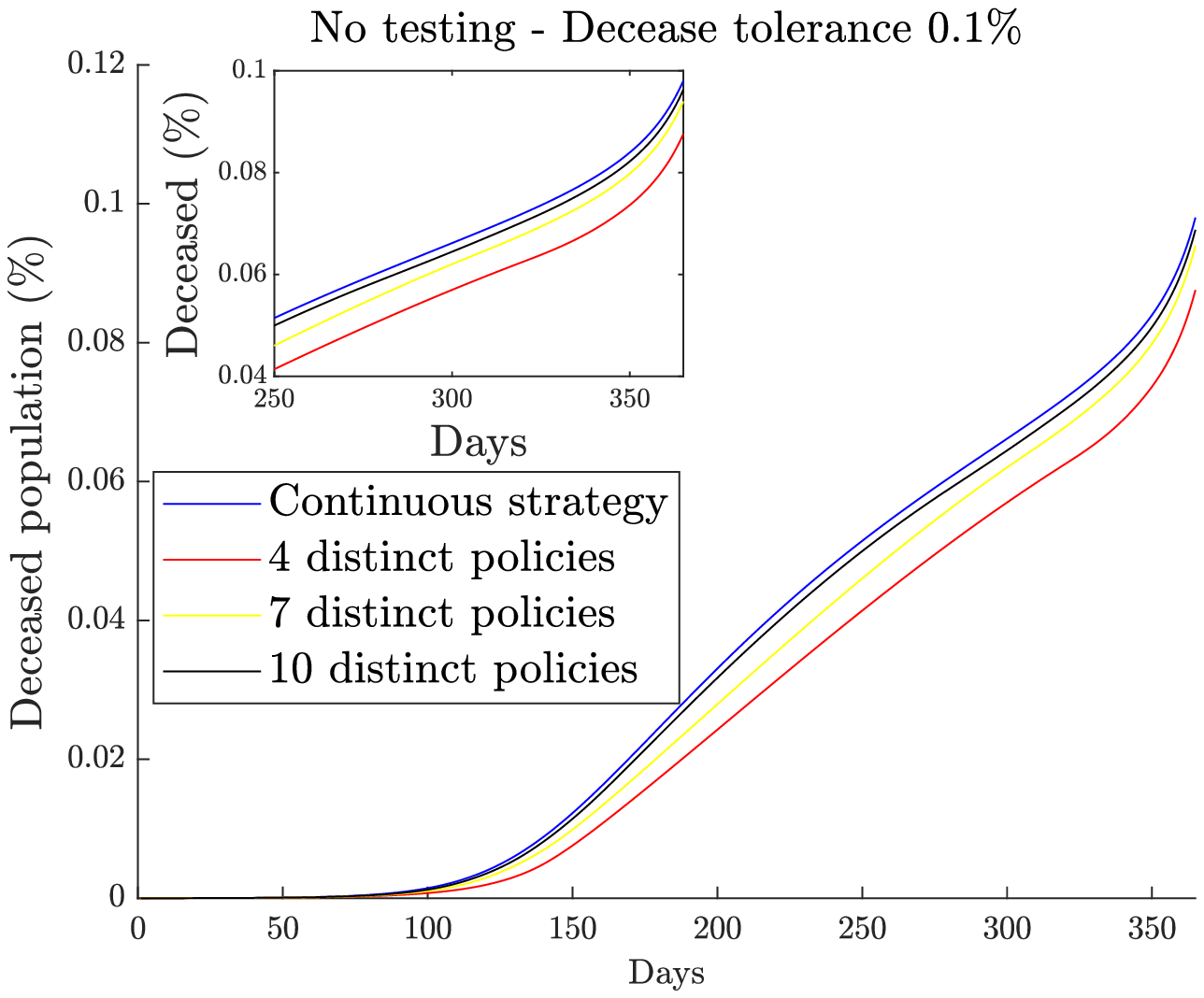}
\end{subfigure}
\vspace{-2mm} \caption[Intervention strategy and decease rate, $0.1\%$ decease tolerance, no testing]{\textbf{Intervention strategy and decease rate, $0.1\%$ decease tolerance, no testing.} Intensity of government intervention {(left)} and portion of deceased population {(right)}
when no testing is performed and a decease tolerance of $0.1\%$ is adopted
  when (i) a continuously changing strategy is selected and (ii) discrete implementations of the selected strategy are considered, allowing $4, 7$ and $10$ policy levels and $6, 12$ and $18$ policy changes respectively. }\vspace{-2mm}
\end{figure}

\begin{figure}[H]
\begin{subfigure}{0.5\columnwidth}
\centering
\includegraphics[scale=0.55]{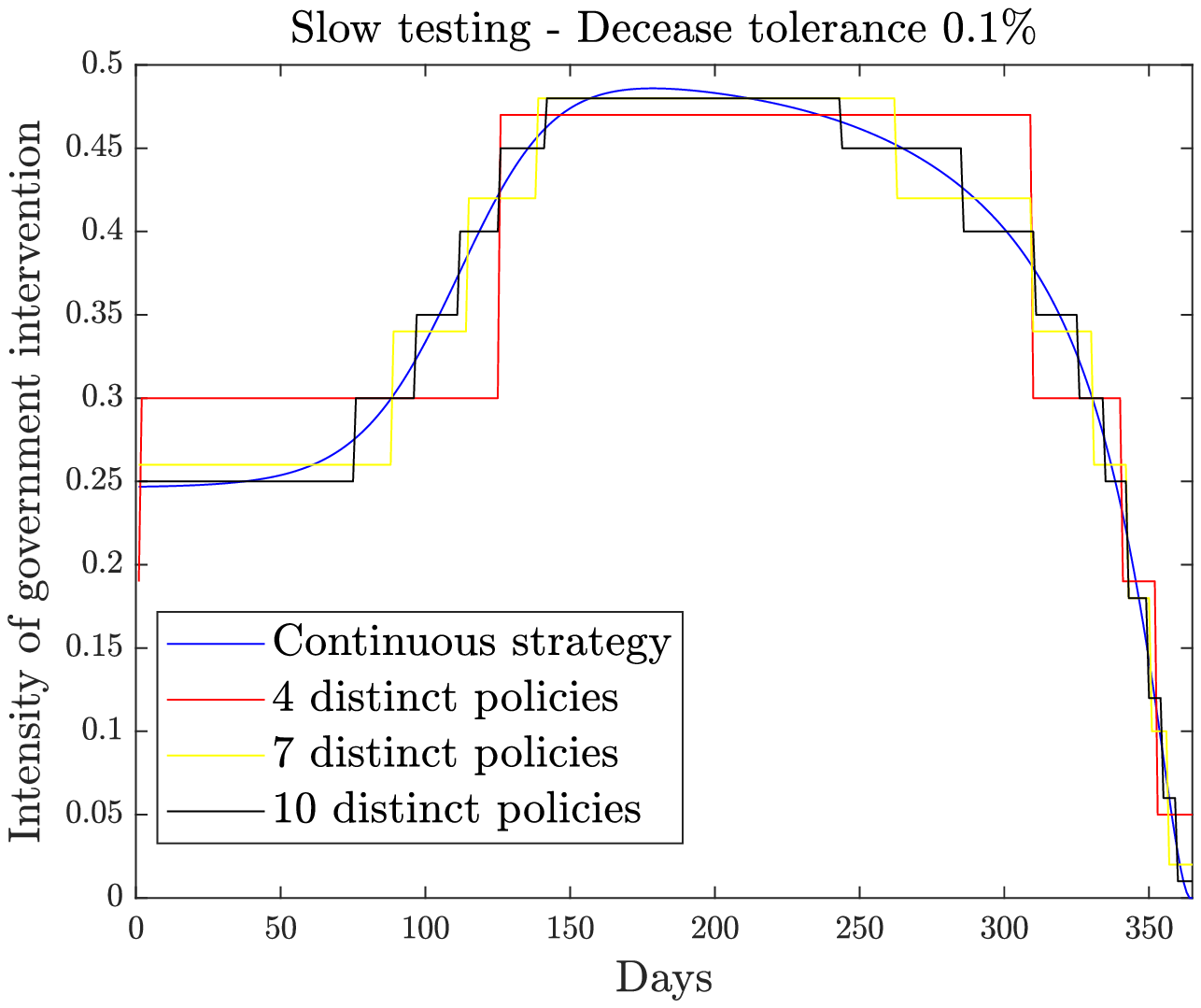}
\end{subfigure} 
\begin{subfigure}{0.5\columnwidth}
\centering
\includegraphics[scale=0.55]{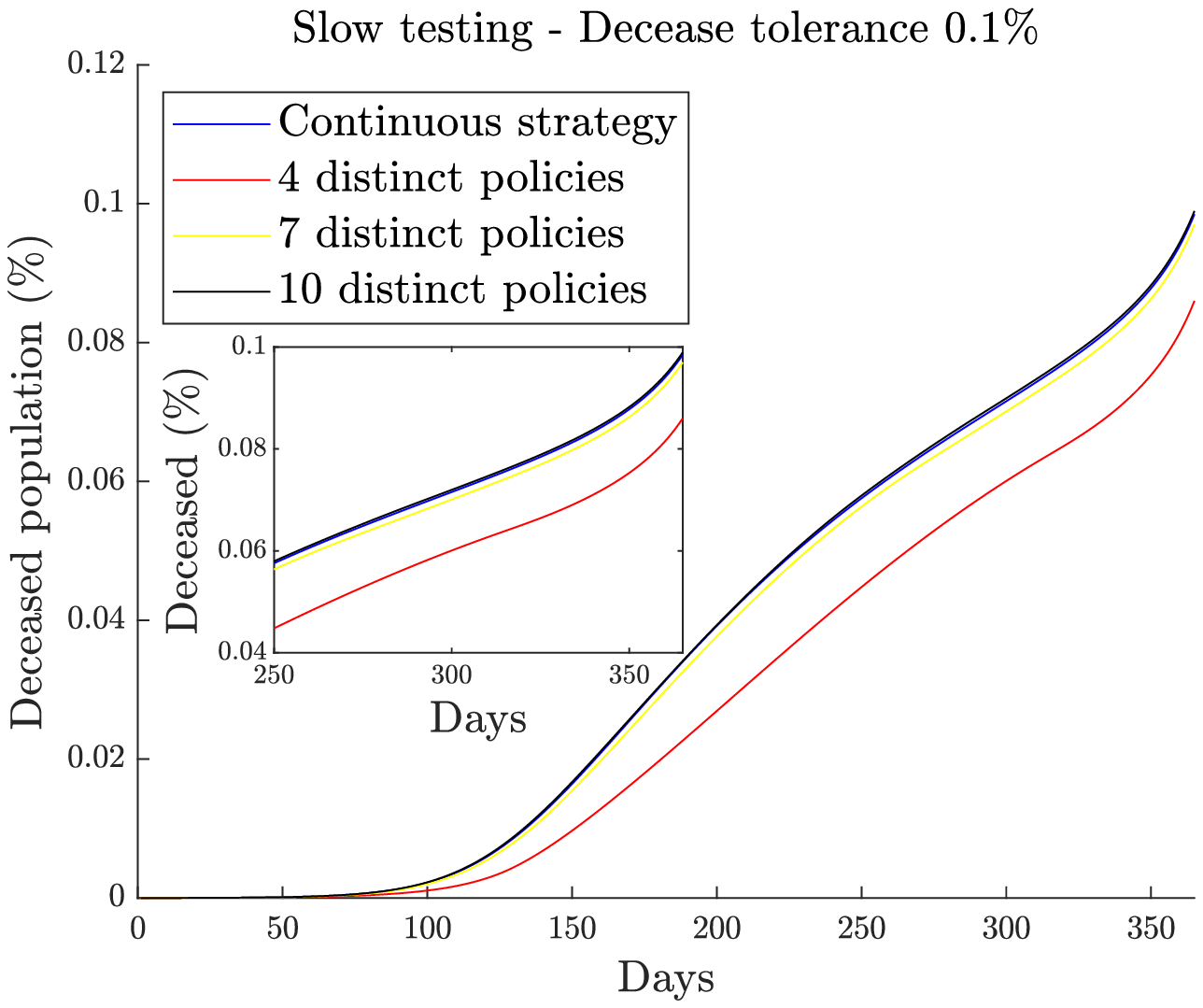}
\end{subfigure}
\vspace{-2mm} \caption[Intervention strategy and decease rate, $0.1\%$ decease tolerance, slow testing]{\textbf{Intervention strategy and decease rate, $0.1\%$ decease tolerance, slow testing.} Intensity of government intervention {(left)} and portion of deceased population {(right)}
when  a slow testing policy and a decease tolerance of $0.1\%$ are adopted
  when (i) a continuously changing strategy is selected and (ii) discrete implementations of the selected strategy are considered, allowing $4, 7$ and $10$ policy levels and $6, 12$ and $18$ policy changes respectively. }\vspace{-2mm}
\end{figure}

\begin{figure}[H]
\begin{subfigure}{0.5\columnwidth}
\centering
\includegraphics[scale=0.55]{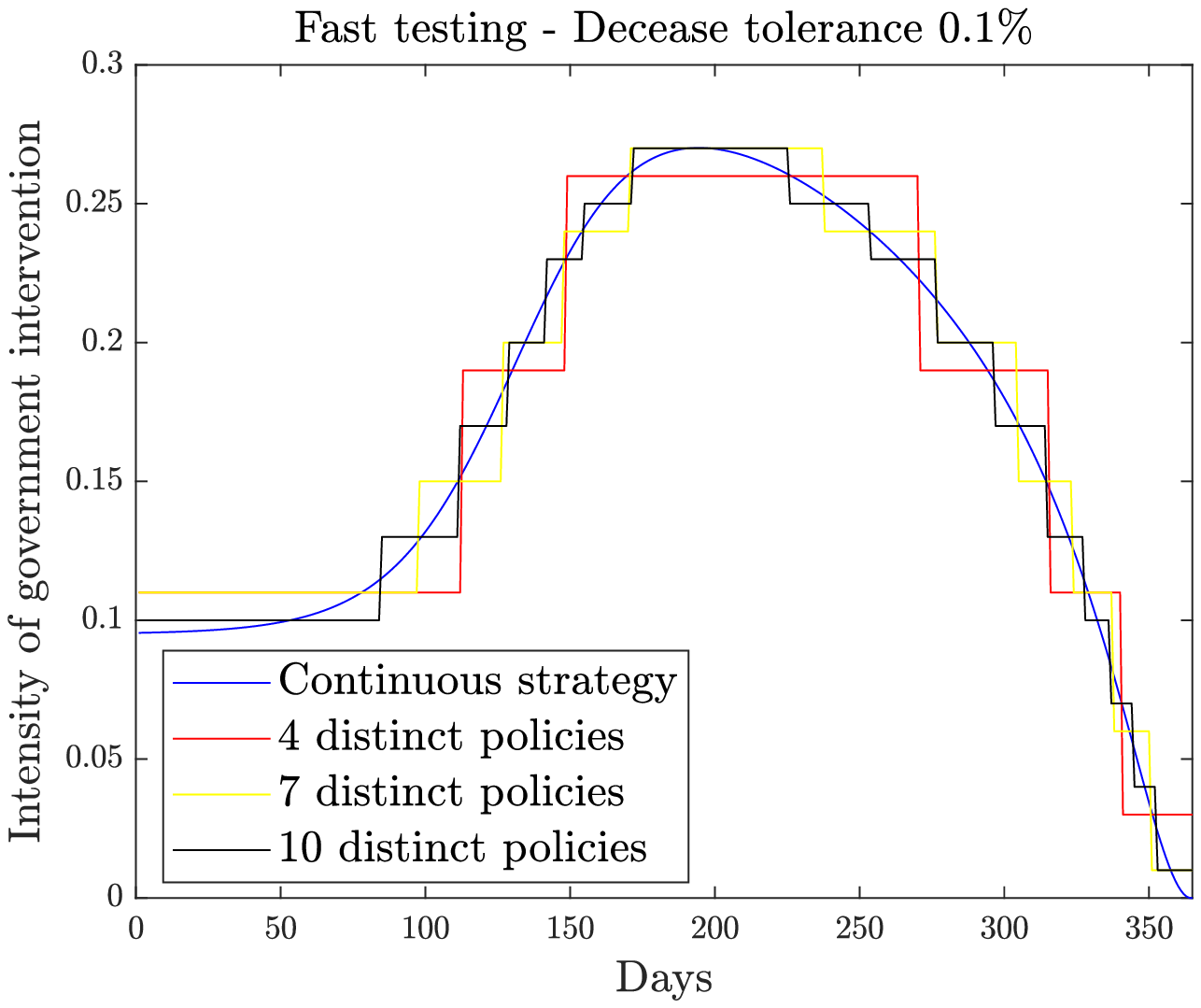}
\end{subfigure} 
\begin{subfigure}{0.5\columnwidth}
\centering
\includegraphics[scale=0.55]{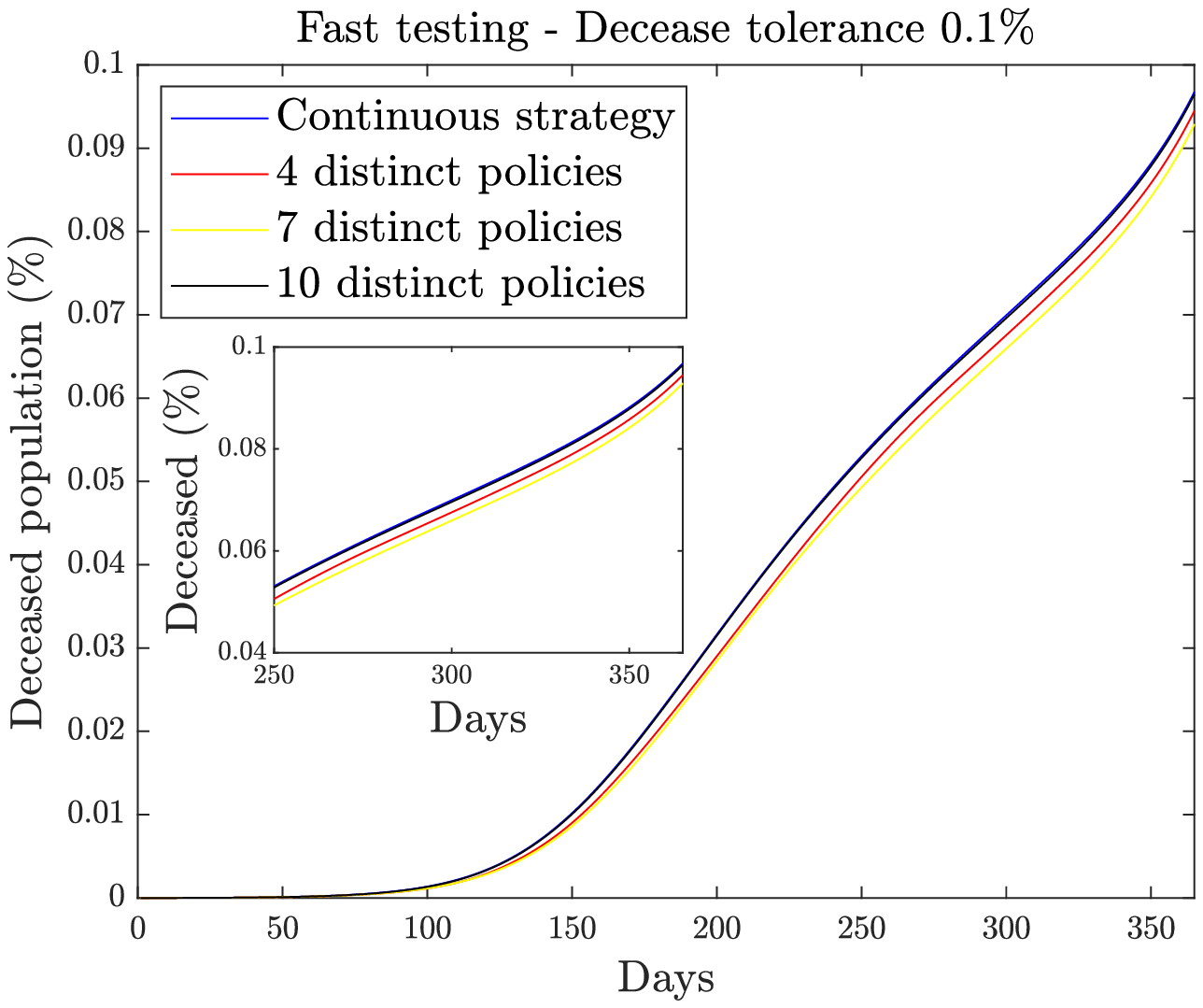}
\end{subfigure}
\vspace{-2mm} \caption[Intervention strategy and decease rate, $0.1\%$ decease tolerance, fast testing]{\textbf{Intervention strategy and decease rate, $0.1\%$ decease tolerance, fast testing.} Intensity of government intervention {(left)} and portion of deceased population {(right)}
when a fast testing policy and a decease tolerance of $0.1\%$ are adopted
  when (i) a continuously changing strategy is selected and (ii) discrete implementations of the selected strategy are considered, allowing $4, 7$ and $10$ policy levels and $6, 12$ and $18$ policy changes respectively. }\vspace{-2mm}
\end{figure}

\begin{figure}[H]
\begin{subfigure}{0.5\columnwidth}
\centering
\includegraphics[scale=0.54]{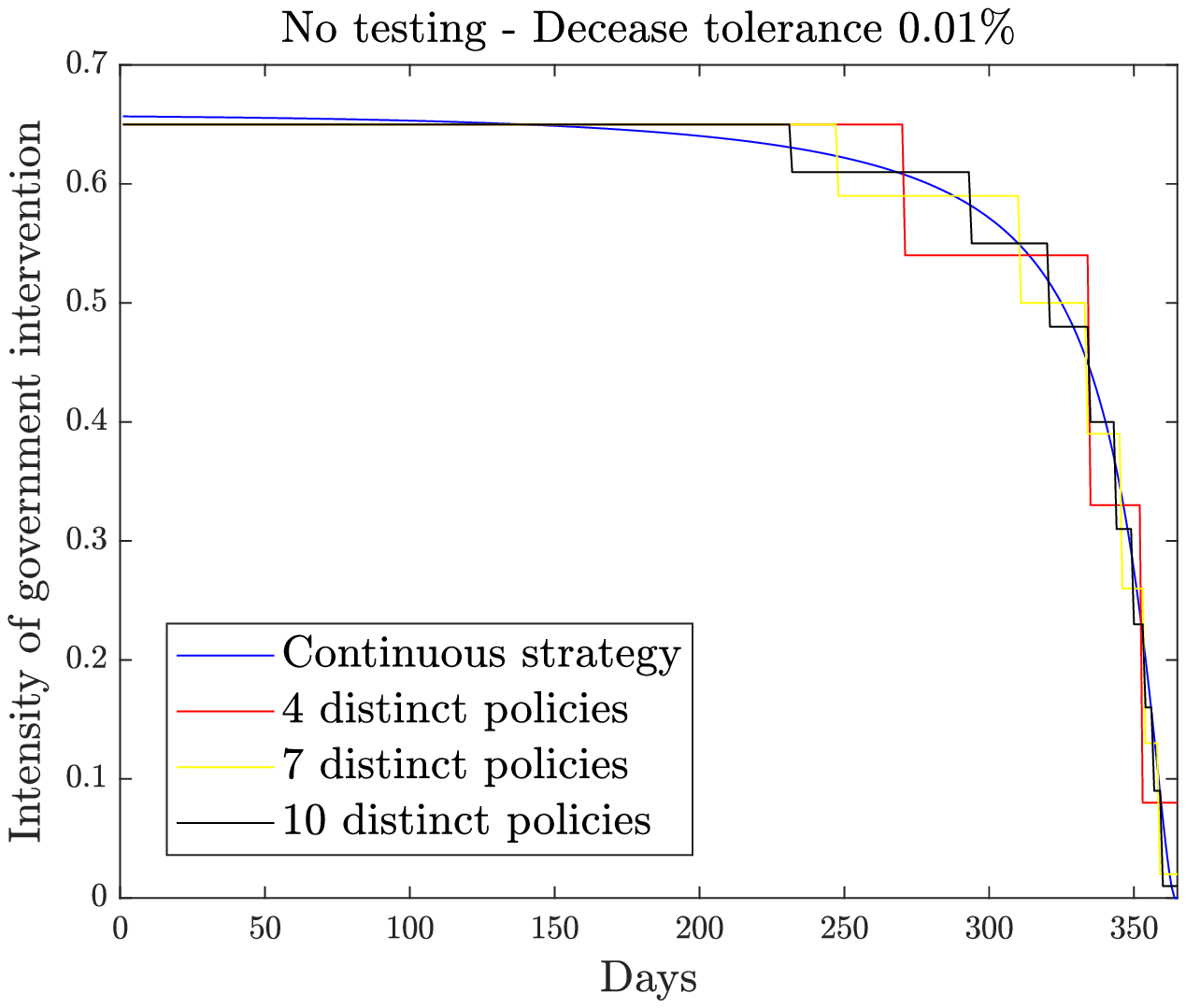}
\end{subfigure} 
\begin{subfigure}{0.5\columnwidth}
\centering
\includegraphics[scale=0.54]{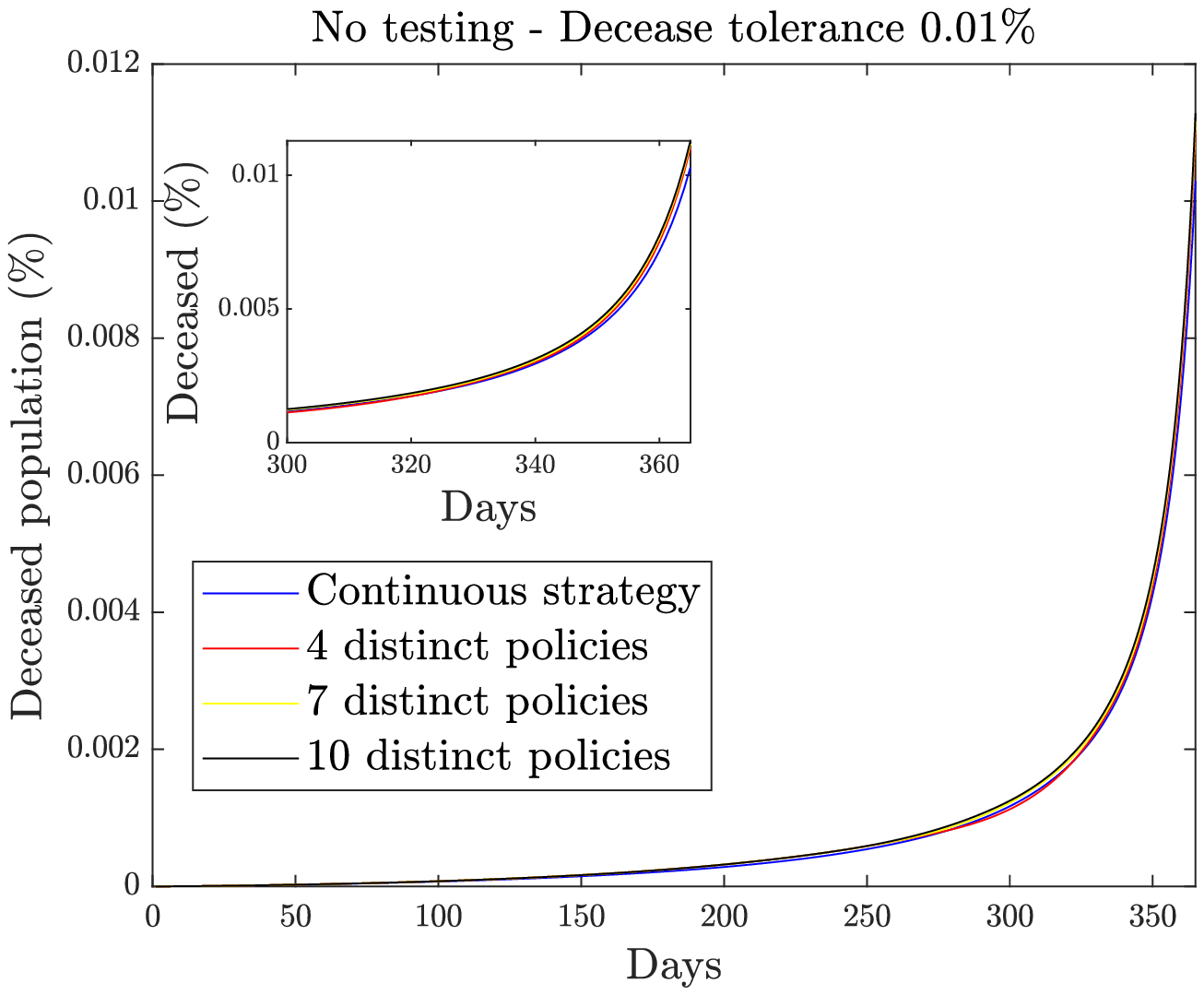}
\end{subfigure}
\vspace{-2mm} \caption[Intervention strategy and decease rate, $0.01\%$ decease tolerance, no testing]{\textbf{Intervention strategy and decease rate, $0.01\%$ decease tolerance, no testing.} Intensity of government intervention {(left)} and portion of deceased population {(right)}
when no testing is performed and a decease tolerance of $0.01\%$ is adopted
  when (i) a continuously changing strategy is selected and (ii) discrete implementations of the selected strategy are considered, allowing $4, 7$ and $10$ policy levels and $6, 12$ and $18$ policy changes respectively. }\vspace{-2mm}
\end{figure}

\begin{figure}[H]
\begin{subfigure}{0.5\columnwidth}
\centering
\includegraphics[scale=0.55]{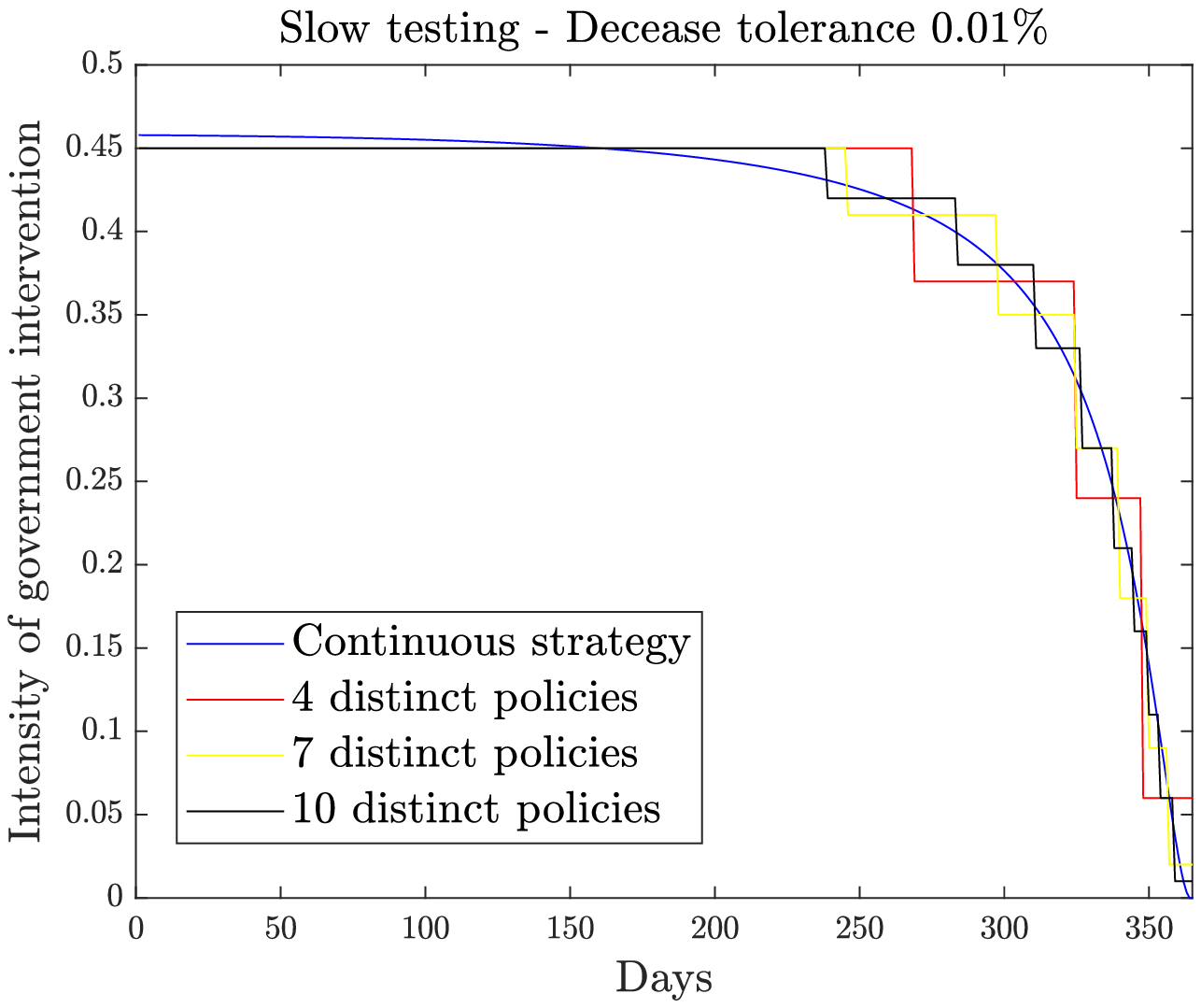}
\end{subfigure} 
\begin{subfigure}{0.5\columnwidth}
\centering
\includegraphics[scale=0.55]{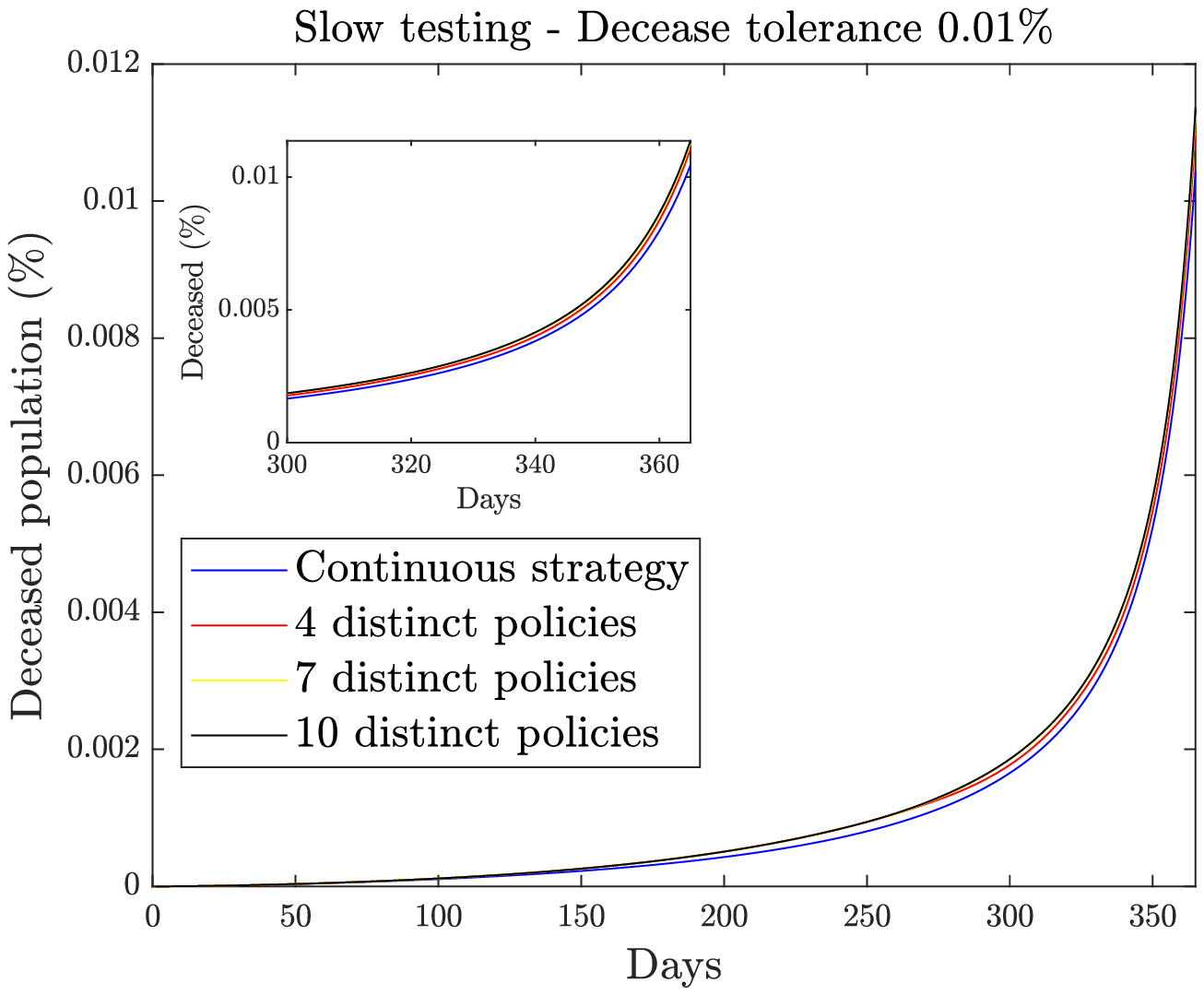}
\end{subfigure}
\vspace{-2mm} \caption[Intervention strategy and decease rate, $0.01\%$ decease tolerance, slow testing]{\textbf{Intervention strategy and decease rate, $0.01\%$ decease tolerance, slow testing.} Intensity of government intervention {(left)} and portion of deceased population {(right)}
when a slow testing policy and a decease tolerance of $0.01\%$ are adopted
  when (i) a continuously changing strategy is selected and (ii) discrete implementations of the selected strategy are considered, allowing $4, 7$ and $10$ policy levels and $6, 12$ and $18$ policy changes respectively. }\vspace{-2mm}
\end{figure}

\begin{figure}[H]
\begin{subfigure}{0.5\columnwidth}
\centering
\includegraphics[scale=0.55]{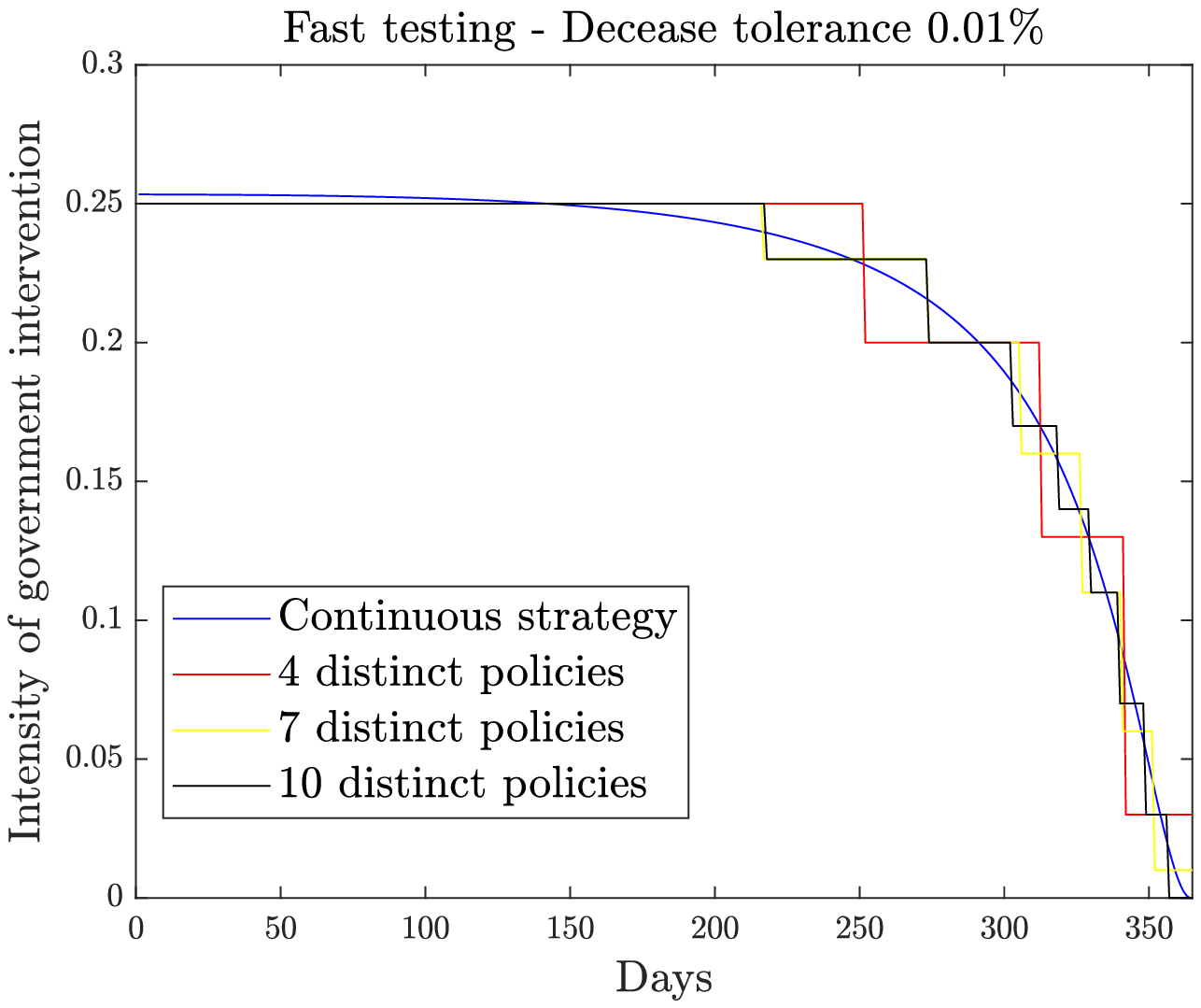}
\end{subfigure} 
\begin{subfigure}{0.5\columnwidth}
\centering
\includegraphics[scale=0.55]{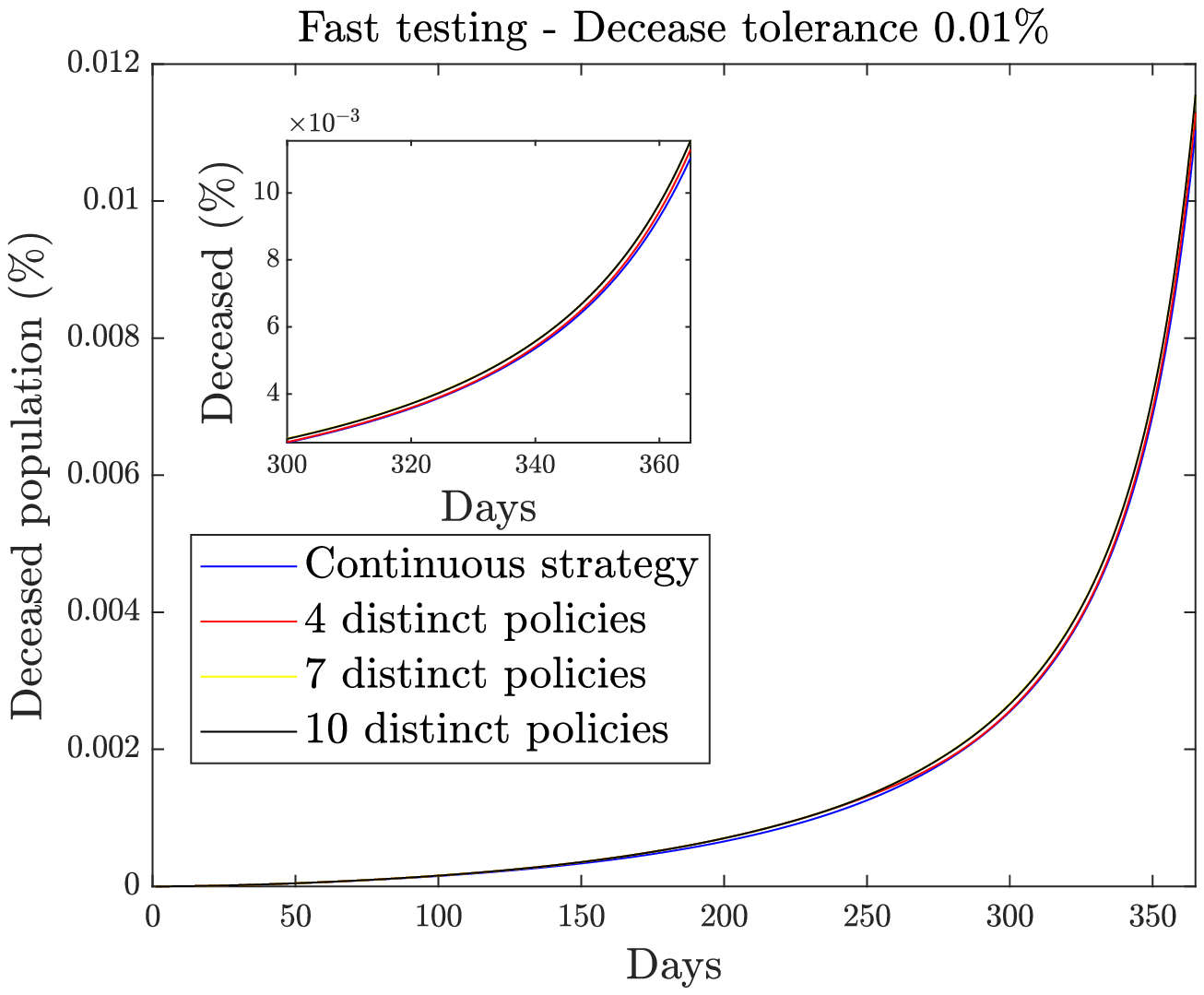}
\end{subfigure}
\vspace{-4mm} \caption[Intervention strategy and decease rate, $0.01\%$ decease tolerance, fast testing]{\textbf{Intervention strategy and decease rate, $0.01\%$ decease tolerance, fast testing.} Intensity of government intervention {(left)} and portion of deceased population {(right)}
when a fast testing policy and a decease tolerance of $0.01\%$ are adopted
  when (i) a continuously changing strategy is selected and (ii) discrete implementations of the selected strategy are considered, allowing $4, 7$ and $10$ policy levels and $6, 12$ and $18$ policy changes respectively. }
  \vspace{-2mm}
\label{Policy_decease_8}
\end{figure}

\begin{figure}[H]
\begin{subfigure}{0.5\columnwidth}
\centering
\includegraphics[scale=0.55]{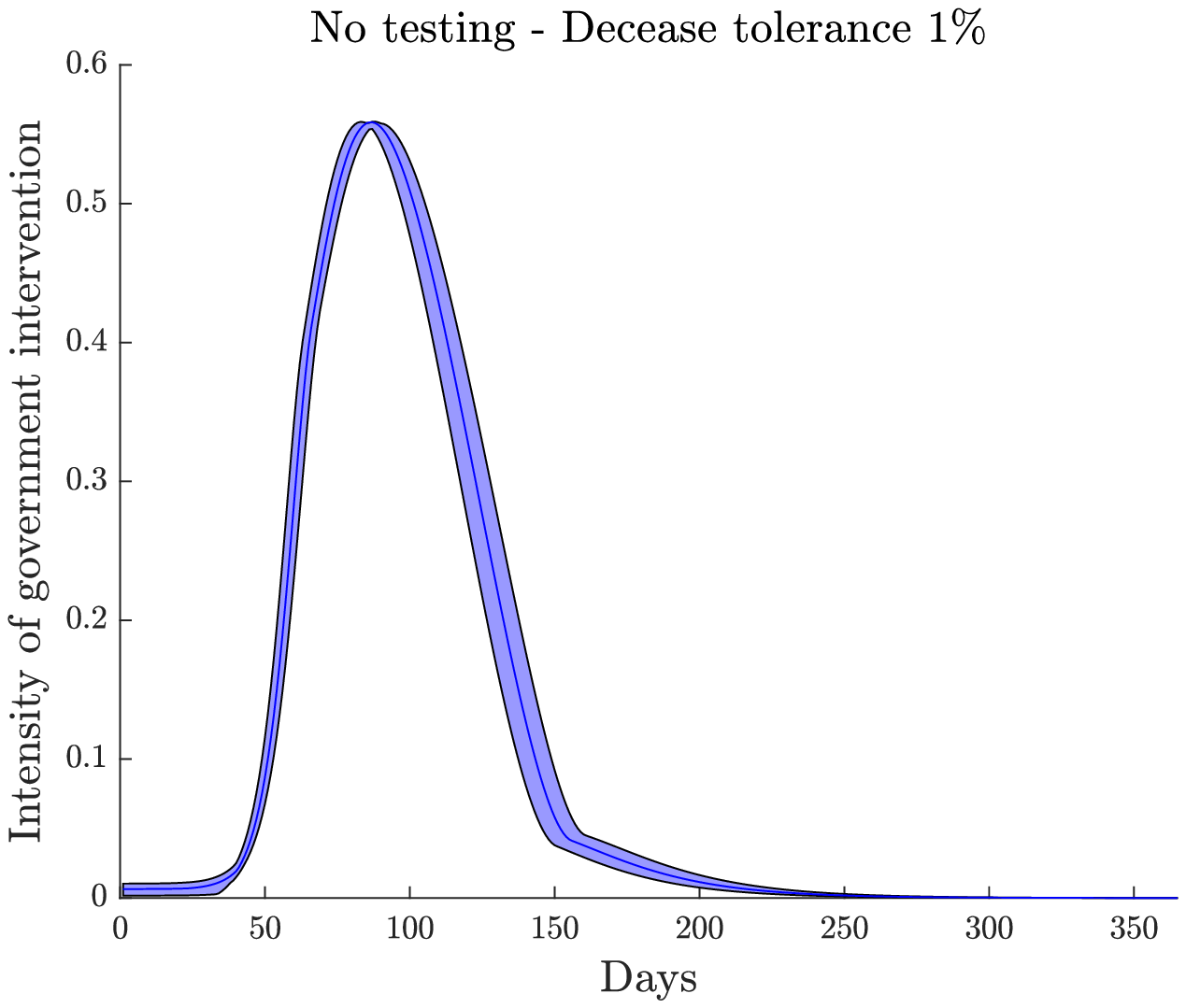}
\end{subfigure} 
\begin{subfigure}{0.5\columnwidth}
\centering
\includegraphics[scale=0.55]{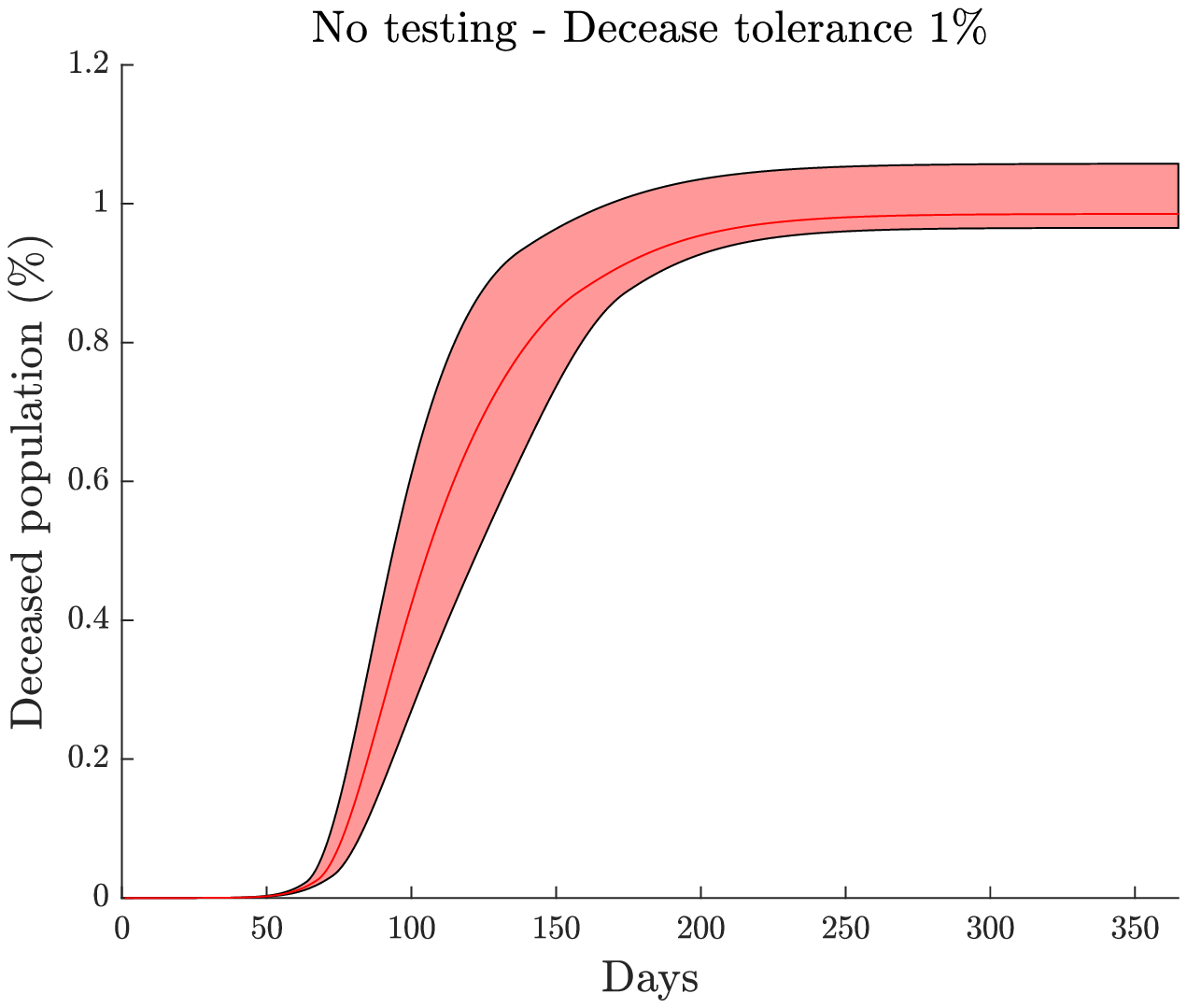}
\end{subfigure}
\vspace{-4mm} \caption[Effect of uncertainty in $\bar{R}_0$ on optimal strategy and decease rate, $1\%$ decease tolerance, no testing]{\textbf{Effect of uncertainty in $\bar{R}_0$ on optimal strategy and decease rate, $1\%$ decease tolerance, no testing.} Both subfigures consider the range $\bar{R}_0 \in [3.17, 3.38]$, a no testing policy and a decease tolerance of $1\%$. 
(left)~Ranges of optimal strategies, (right)  Ranges of \ak{aggregate deceases}  when the optimal strategy obtained  based on $\bar{R}_0 = 3.27$ and infection mortality rate of $0.66\%$ is implemented.
The darker line within the presented ranges corresponds to $\bar{R}_0 = 3.27$.}
\vspace{-2mm}
\label{robust_R_1}
\end{figure}

\begin{figure}[H]
\begin{subfigure}{0.5\columnwidth}
\centering
\includegraphics[scale=0.55]{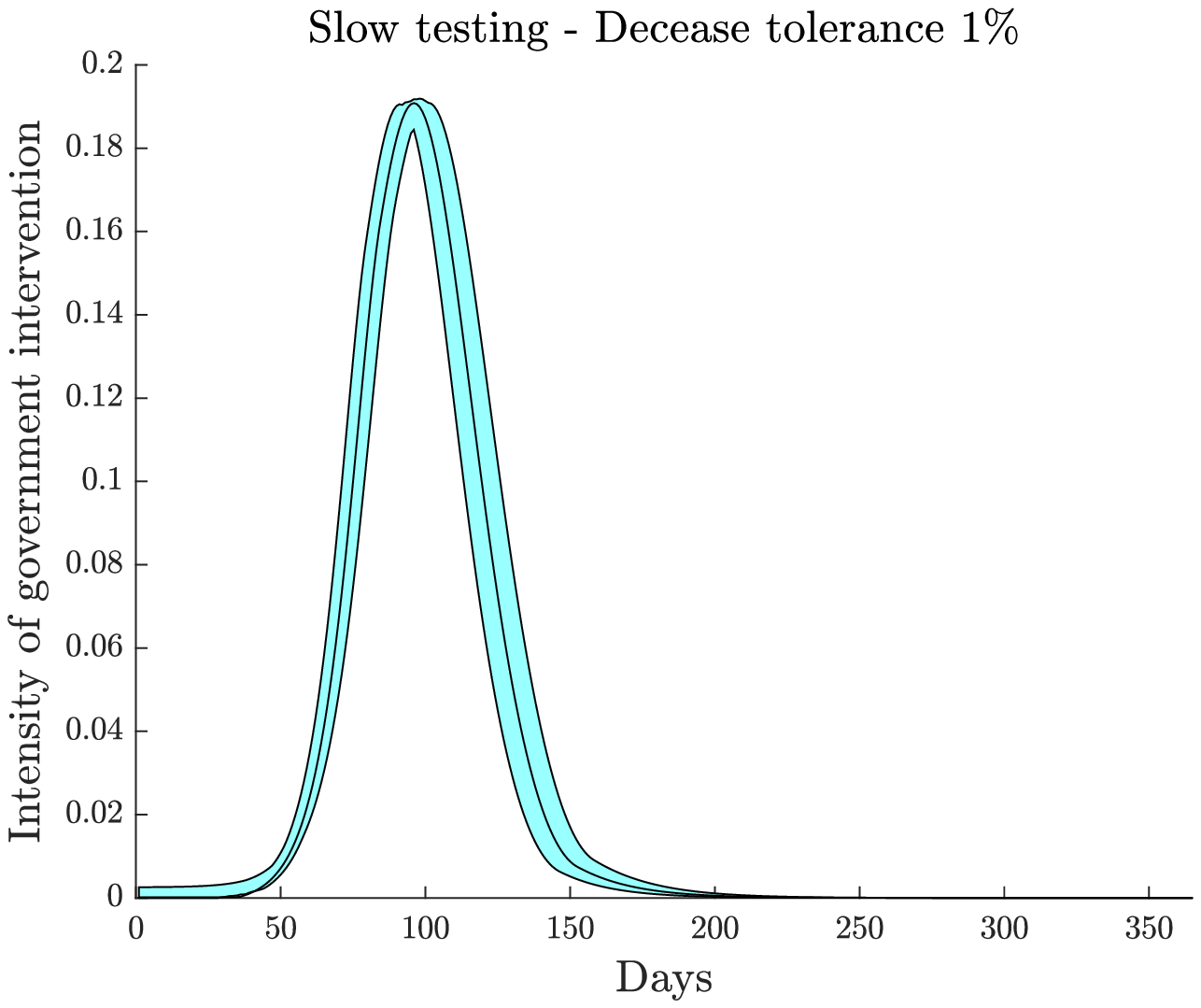}
\end{subfigure} 
\begin{subfigure}{0.5\columnwidth}
\centering
\includegraphics[scale=0.55]{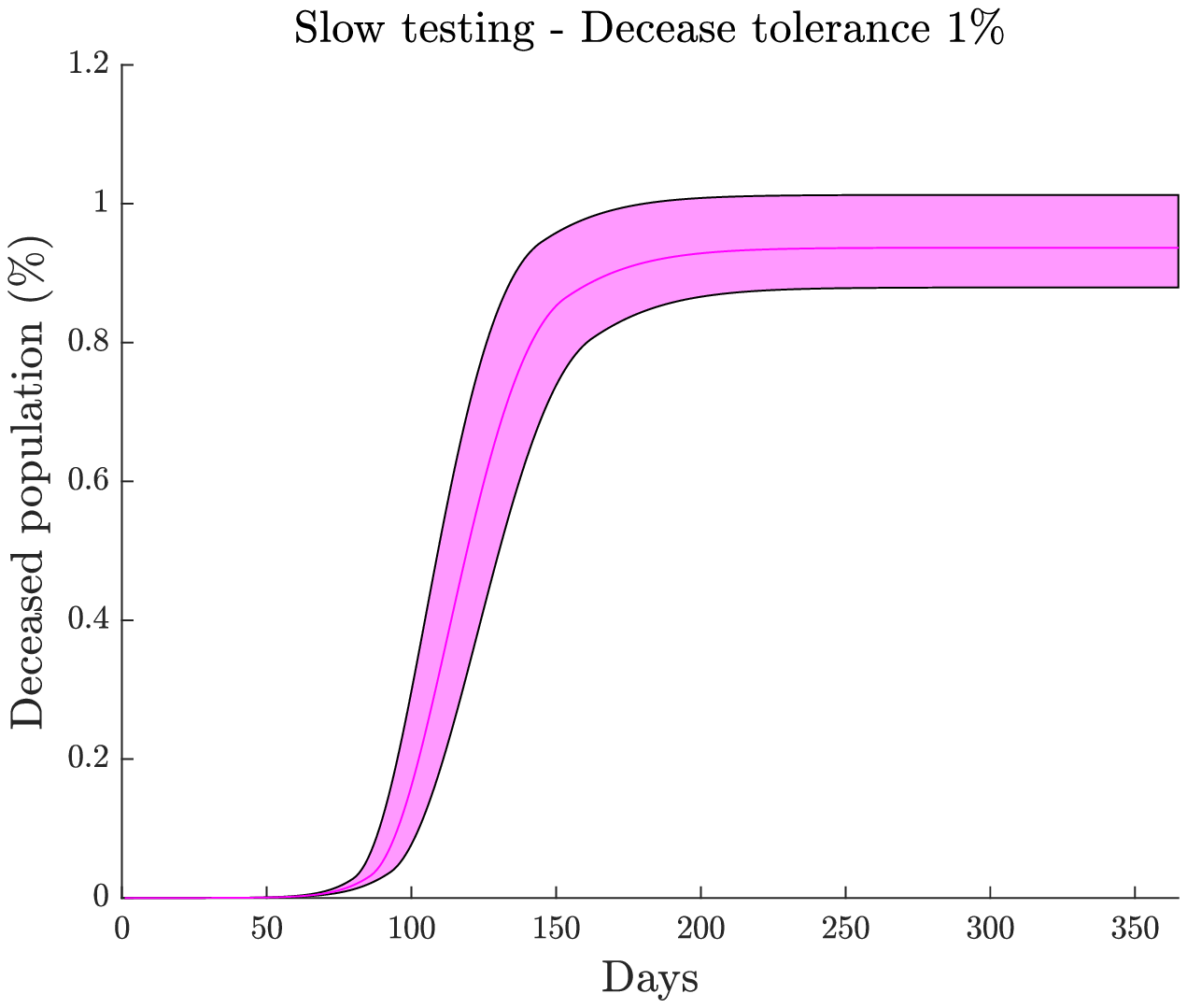}
\end{subfigure}
\vspace{-4mm} \caption[Effect of uncertainty in $\bar{R}_0$ on optimal strategy and decease rate, $1\%$ decease tolerance, slow testing]{\textbf{Effect of uncertainty in $\bar{R}_0$ on optimal strategy and decease rate, $1\%$ decease tolerance, slow testing.} 
Both subfigures consider the range $\bar{R}_0 \in [3.17, 3.38]$,   {a slow testing policy} and a decease tolerance of $1\%$.
(left)~Ranges of optimal strategies, (right)  Ranges of \ak{aggregate deceases}  when the optimal strategy obtained  based on $\bar{R}_0 = 3.27$ and infection mortality rate of $0.66\%$ is implemented.
The darker line within the presented ranges corresponds to $\bar{R}_0 = 3.27$.}
\vspace{-2mm}
\end{figure}

\begin{figure}[H]
\begin{subfigure}{0.5\columnwidth}
\centering
\includegraphics[scale=0.55]{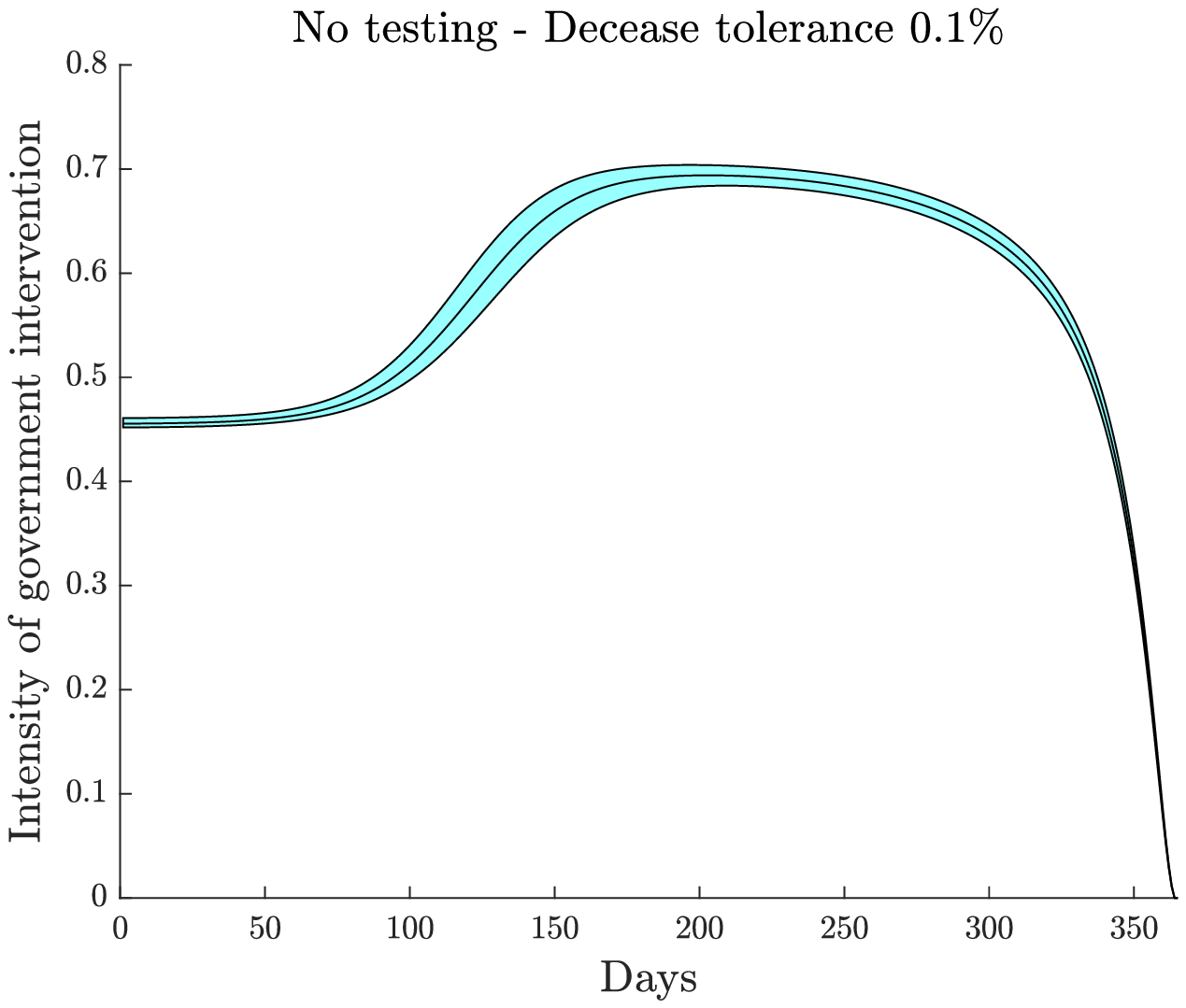}
\end{subfigure} 
\begin{subfigure}{0.5\columnwidth}
\centering
\includegraphics[scale=0.55]{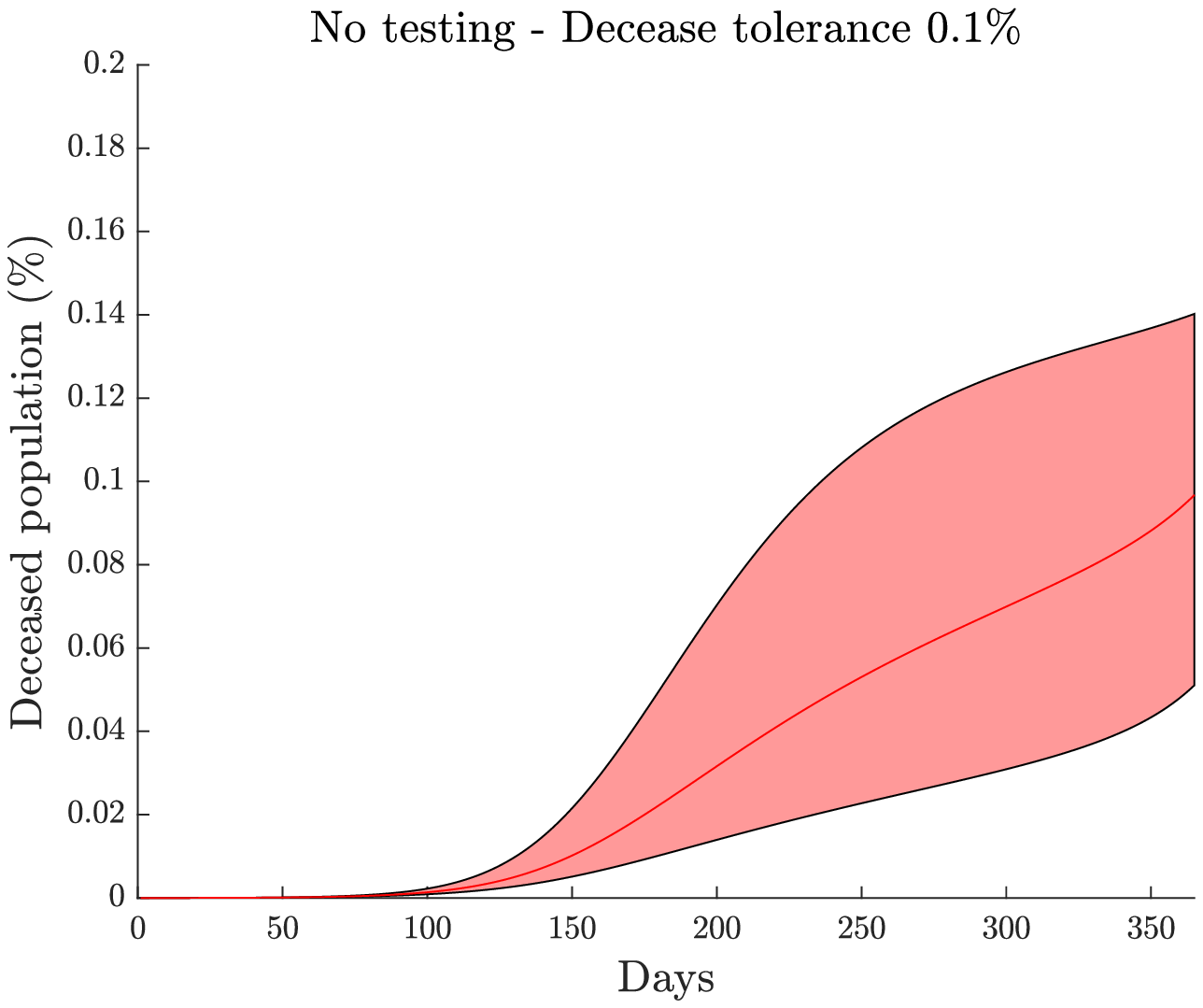}
\end{subfigure}
\vspace{-4mm} \caption[Effect of uncertainty in $\bar{R}_0$ on optimal strategy and decease rate, $0.1\%$ decease tolerance, no testing]{\textbf{Effect of uncertainty in $\bar{R}_0$ on optimal strategy and decease rate, $0.1\%$ decease tolerance, no testing.} 
Both subfigures consider the range $\bar{R}_0 \in [3.17, 3.38]$, a no testing policy and a decease tolerance of $0.1\%$. (left)~Ranges of optimal strategies, (right)  Ranges of \ak{aggregate deceases}  when the optimal strategy obtained  based on $\bar{R}_0 = 3.27$ and infection mortality rate of $0.66\%$ is implemented.
The darker line within the presented ranges corresponds to $\bar{R}_0 = 3.27$.}
\vspace{-2mm}
\end{figure}

\begin{figure}[H]
\begin{subfigure}{0.5\columnwidth}
\centering
\includegraphics[scale=0.55]{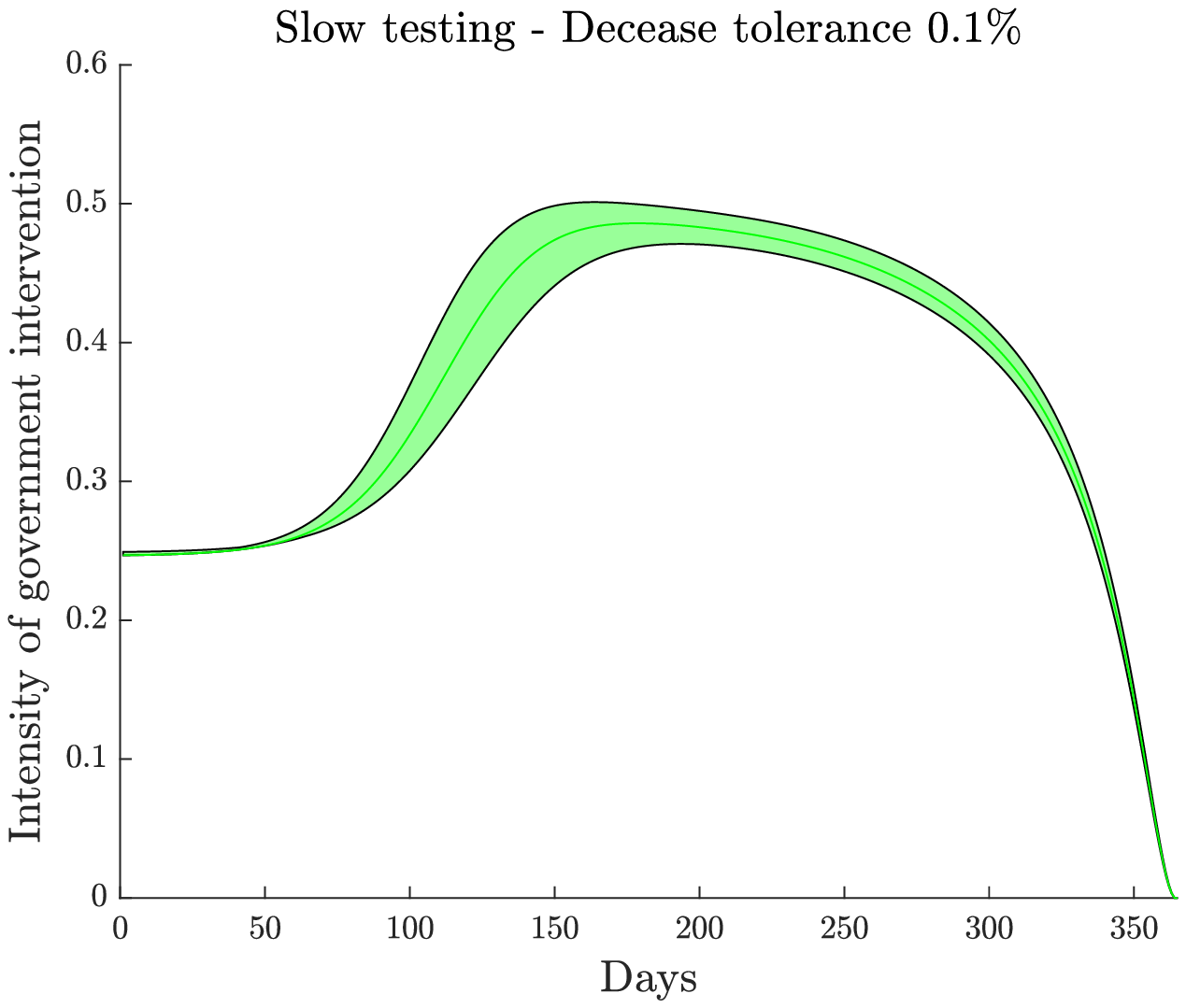}
\end{subfigure} 
\begin{subfigure}{0.5\columnwidth}
\centering
\includegraphics[scale=0.55]{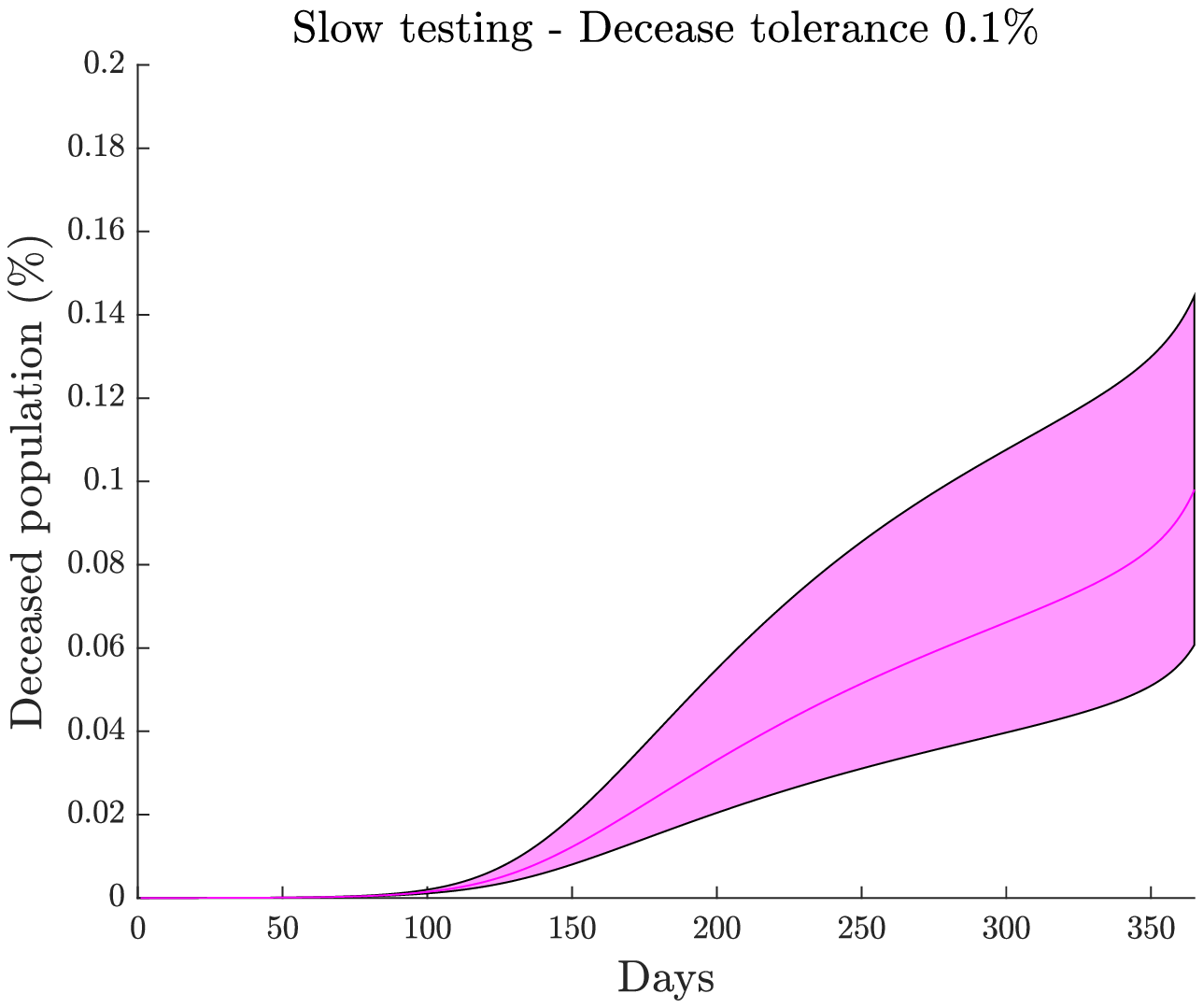}
\end{subfigure}
\vspace{-4mm} \caption[Effect of uncertainty in $\bar{R}_0$ on optimal strategy and decease rate, $0.1\%$ decease tolerance, slow testing]{\textbf{Effect of uncertainty in $\bar{R}_0$ on optimal strategy and decease rate, $0.1\%$ decease tolerance, slow testing.} Both subfigures consider the range $\bar{R}_0 \in [3.17, 3.38]$, {a slow testing policy} and a decease tolerance of $0.1\%$.
(left)~Ranges of optimal strategies, (right)  Ranges of \ak{aggregate deceases}  when the optimal strategy obtained  based on $\bar{R}_0 = 3.27$ and infection mortality rate of $0.66\%$ is implemented.
The darker line within the presented ranges corresponds to $\bar{R}_0 = 3.27$.}
\vspace{-2mm}
\end{figure}

\begin{figure}[H]
\begin{subfigure}{0.5\columnwidth}
\centering
\includegraphics[scale=0.55]{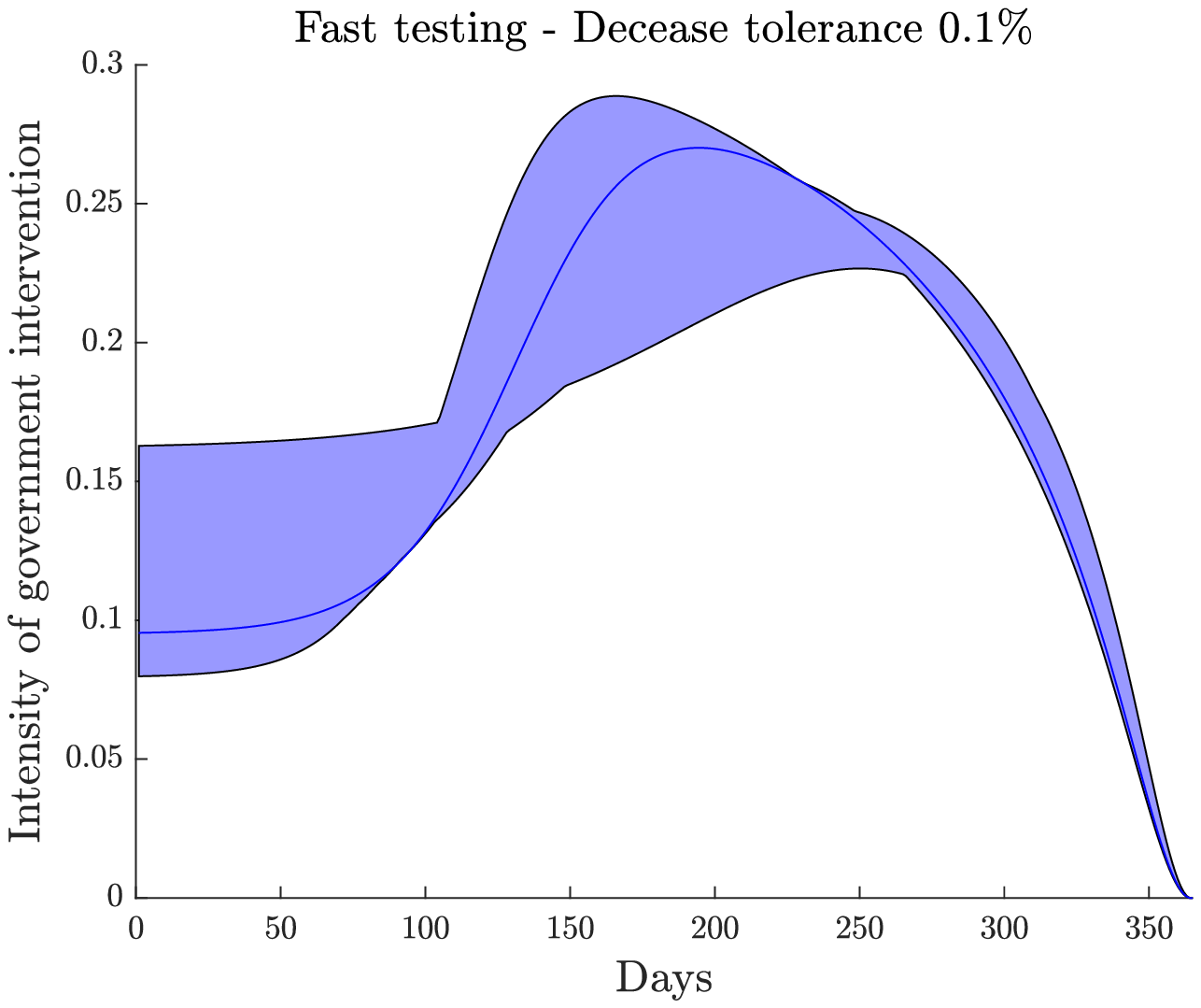}
\end{subfigure} 
\begin{subfigure}{0.5\columnwidth}
\centering
\includegraphics[scale=0.55]{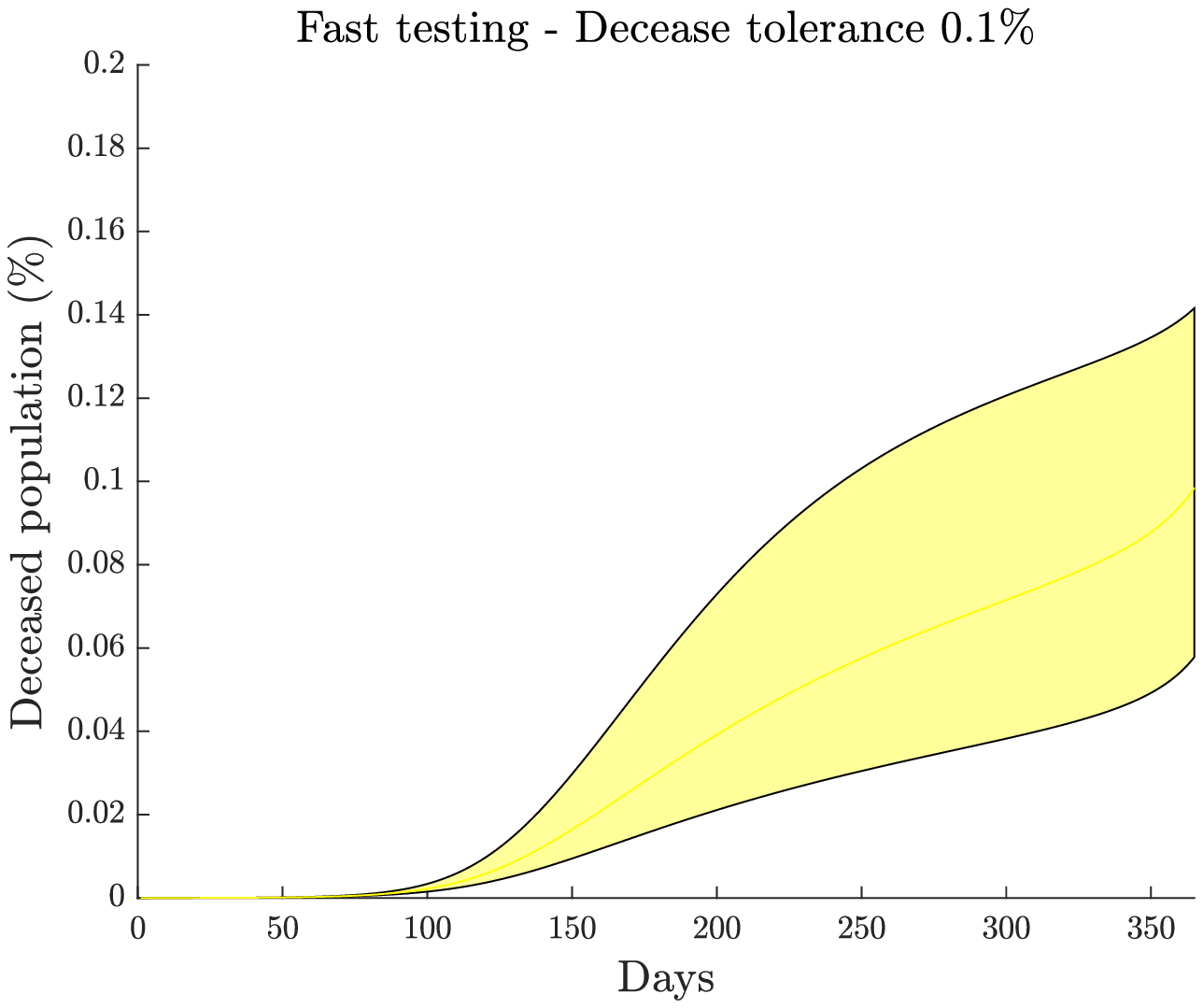}
\end{subfigure}
\vspace{-4mm} \caption[Effect of uncertainty in $\bar{R}_0$ on optimal strategy and decease rate, $0.1\%$ decease tolerance, fast testing]{\textbf{Effect of uncertainty in $\bar{R}_0$ on optimal strategy and decease rate, $0.1\%$ decease tolerance, fast testing.} Both subfigures consider the range $\bar{R}_0 \in [3.17, 3.38]$, {a fast testing policy} and a decease tolerance of $0.1\%$.
(left)~Ranges of optimal strategies, (right)  Ranges of \ak{aggregate deceases}  when the optimal strategy obtained  based on $\bar{R}_0 = 3.27$ and infection mortality rate of $0.66\%$ is implemented.
The darker line within the presented ranges corresponds to $\bar{R}_0 = 3.27$.}
\vspace{-2mm}
\end{figure}

\begin{figure}[H]
\begin{subfigure}{0.5\columnwidth}
\centering
\includegraphics[scale=0.55]{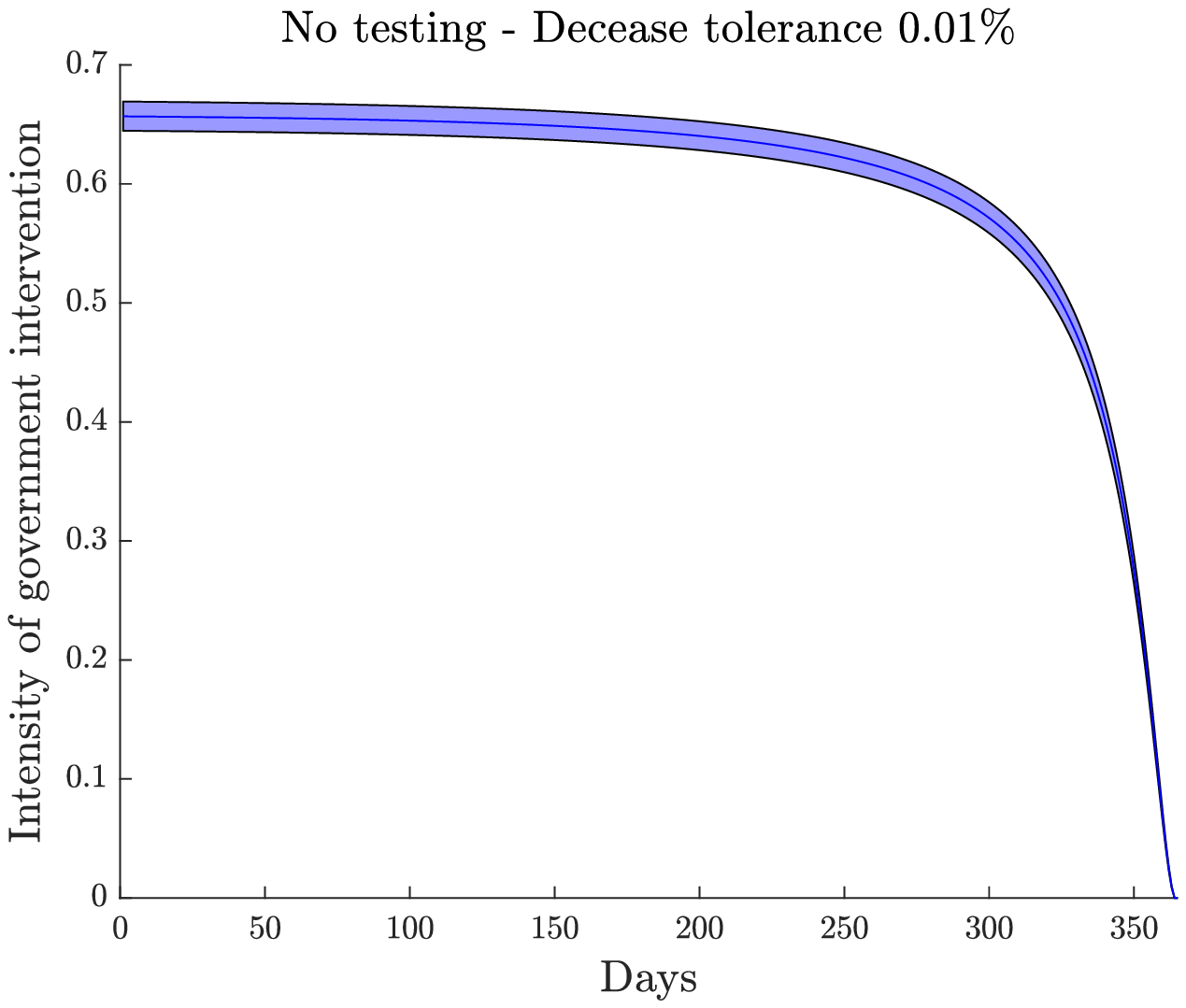}
\end{subfigure} 
\begin{subfigure}{0.5\columnwidth}
\centering
\includegraphics[scale=0.55]{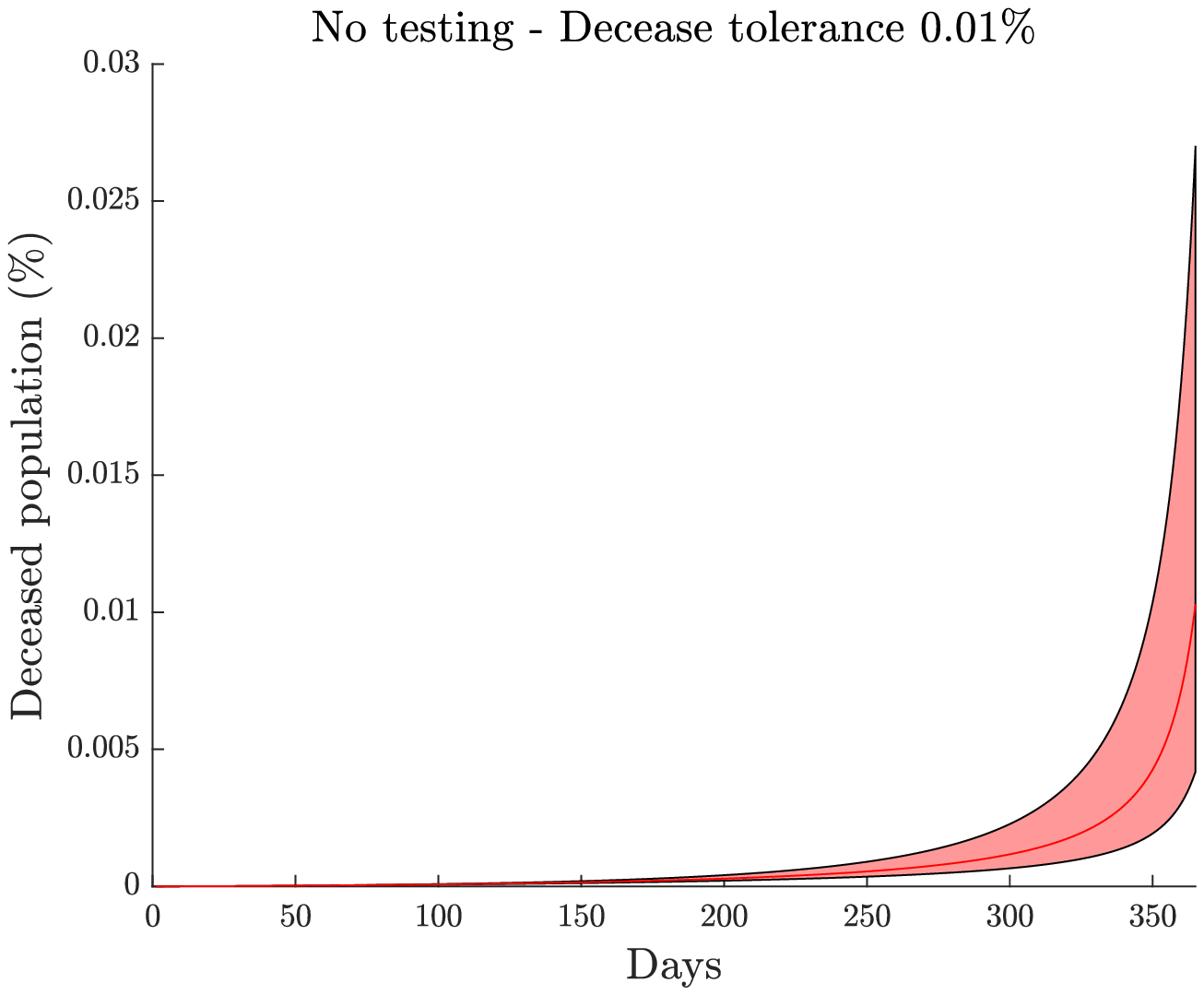}
\end{subfigure}
\vspace{-4mm} \caption[Effect of uncertainty in $\bar{R}_0$ on optimal strategy and decease rate, $0.01\%$ decease tolerance, no testing]{\textbf{Effect of uncertainty in $\bar{R}_0$ on optimal strategy and decease rate, $0.01\%$ decease tolerance, no testing.} Both subfigures consider the range $\bar{R}_0 \in [3.17, 3.38]$, a no testing policy and a decease tolerance of $0.01\%$.
(left)~Ranges of optimal strategies, (right)  Ranges of \ak{aggregate deceases}  when the optimal strategy obtained  based on $\bar{R}_0 = 3.27$ and infection mortality rate of $0.66\%$ is implemented.
The darker line within the presented ranges corresponds to $\bar{R}_0 = 3.27$.}
\vspace{-2mm}
\end{figure}

\begin{figure}[H]
\begin{subfigure}{0.5\columnwidth}
\centering
\includegraphics[scale=0.55]{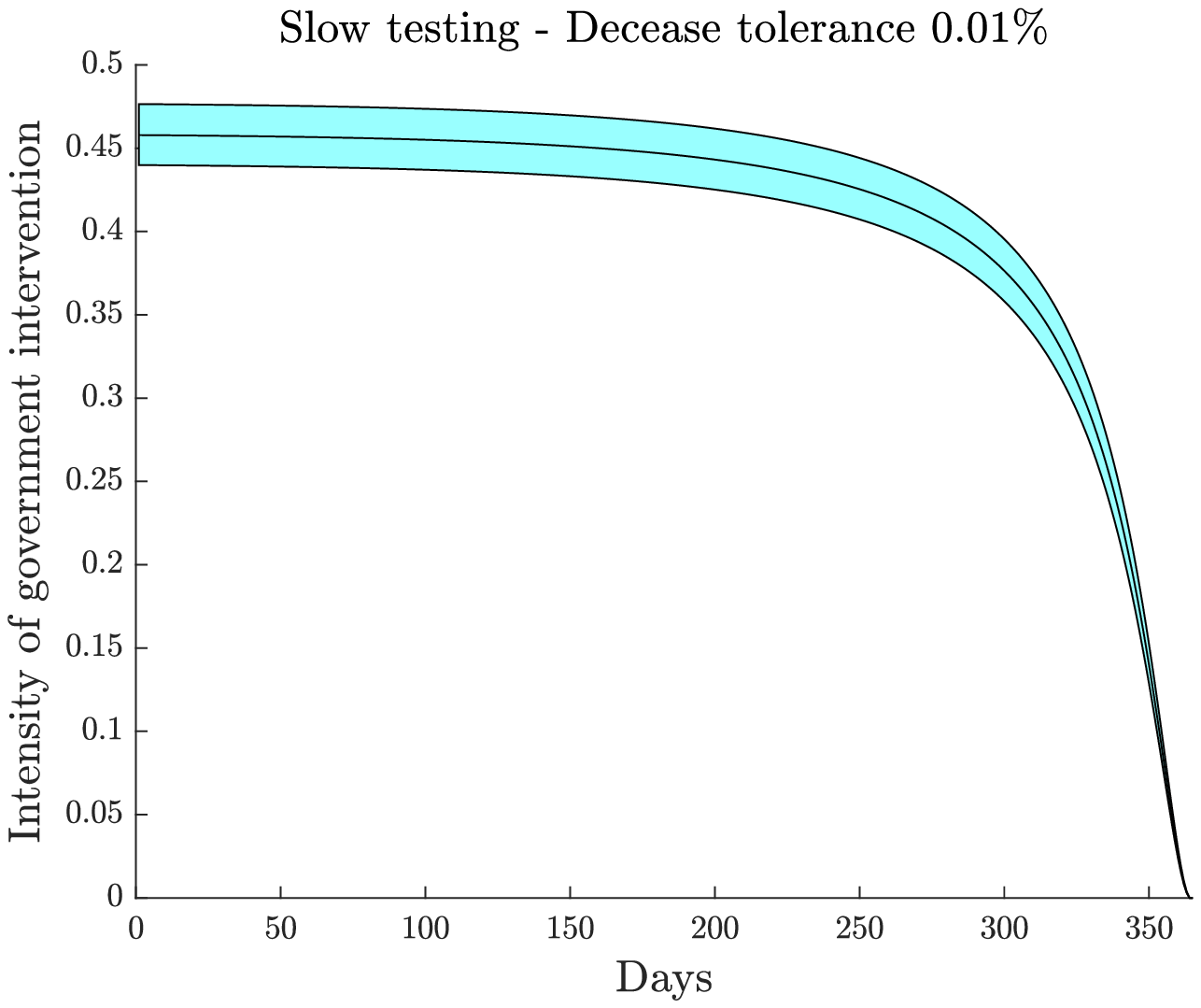}
\end{subfigure} 
\begin{subfigure}{0.5\columnwidth}
\centering
\includegraphics[scale=0.55]{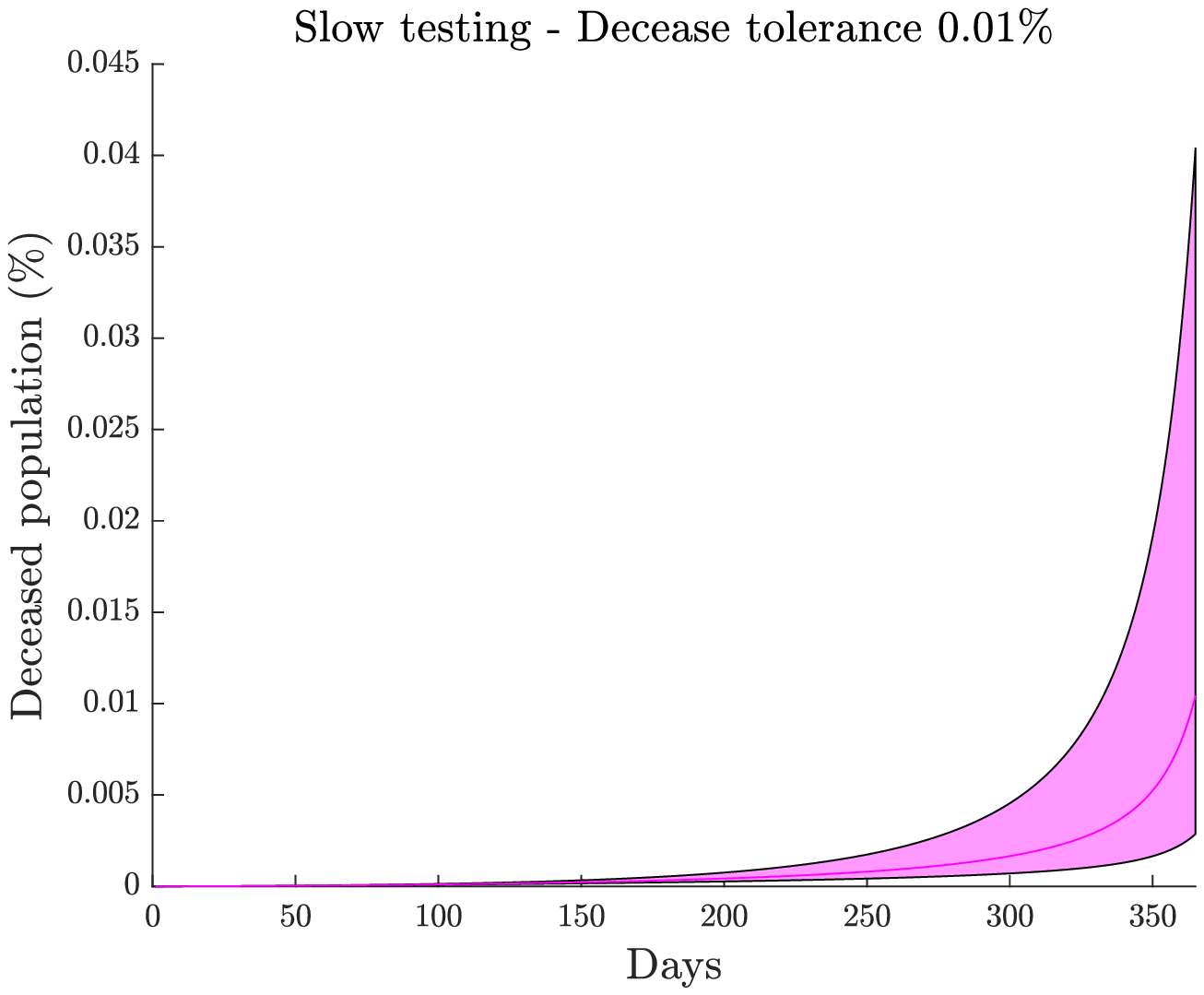}
\end{subfigure}
\vspace{-4mm} \caption[Effect of uncertainty in $\bar{R}_0$ on optimal strategy and decease rate, $0.01\%$ decease tolerance, slow testing]{\textbf{Effect of uncertainty in $\bar{R}_0$ on optimal strategy and decease rate, $0.01\%$ decease tolerance, slow testing.} Both subfigures consider the range $\bar{R}_0 \in [3.17, 3.38]$, {a slow testing policy} and a decease tolerance of $0.01\%$.
(left)~Ranges of optimal strategies, (right)  Ranges of \ak{aggregate deceases}  when the optimal strategy obtained  based on $\bar{R}_0 = 3.27$ and infection mortality rate of $0.66\%$ is implemented.
The darker line corresponds to $\bar{R}_0 = 3.27$.}
\vspace{-2mm}
\end{figure}

\begin{figure}[H]
\begin{subfigure}{0.5\columnwidth}
\centering
\includegraphics[scale=0.55]{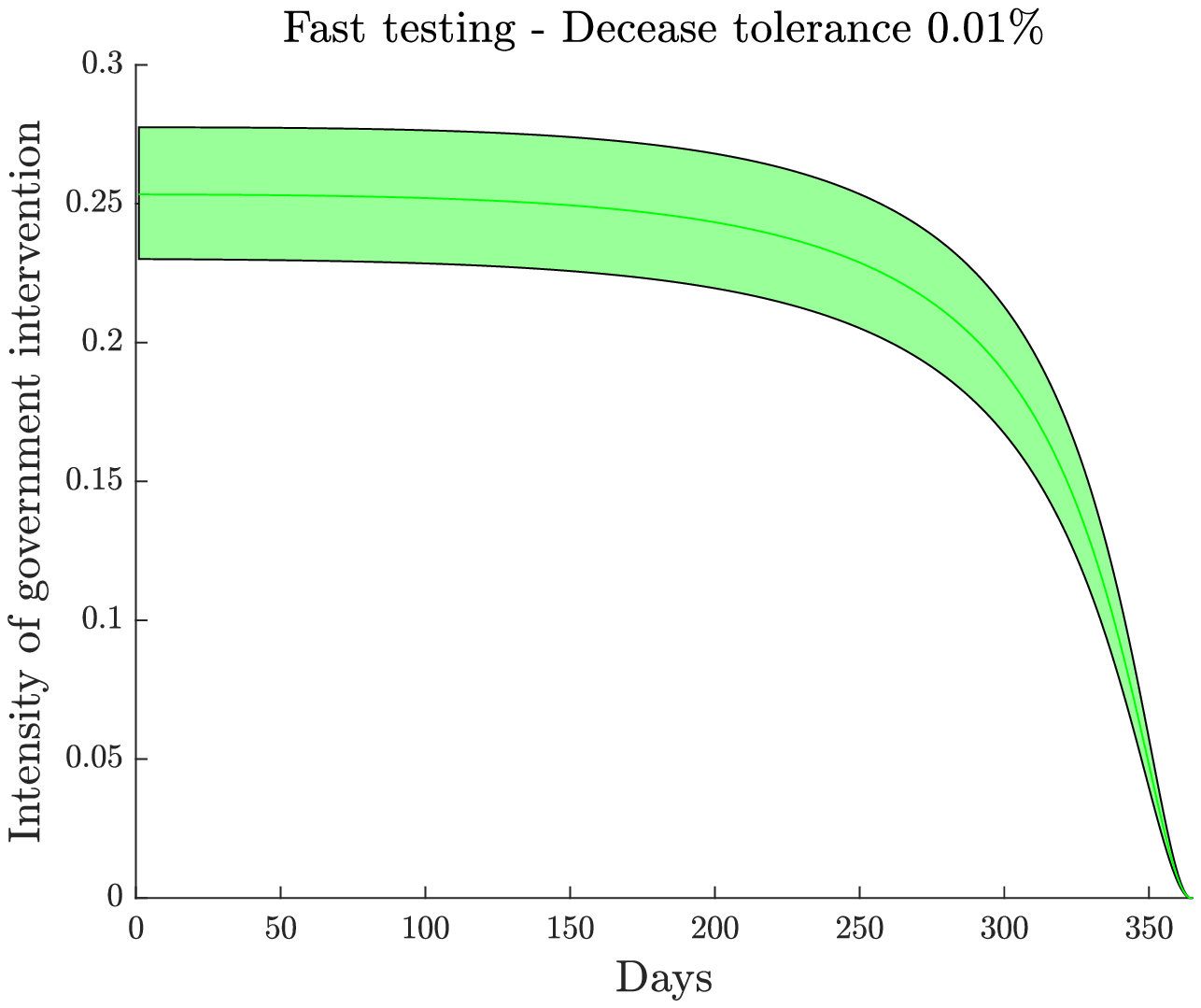}
\end{subfigure} 
\begin{subfigure}{0.5\columnwidth}
\centering
\includegraphics[scale=0.55]{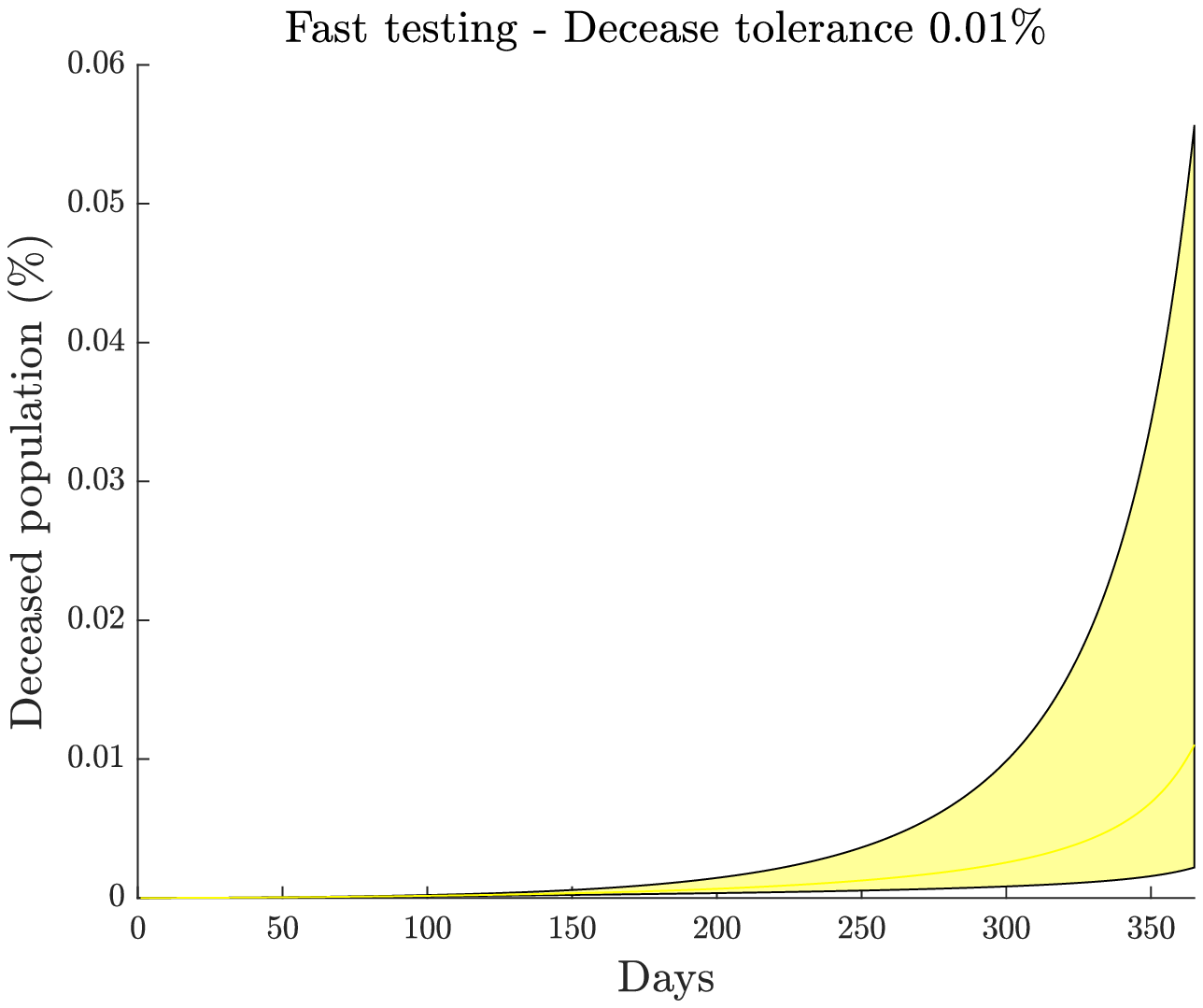}
\end{subfigure}
\vspace{-4mm} \caption[Effect of uncertainty in $\bar{R}_0$ on optimal strategy and decease rate, $0.01\%$ decease tolerance, fast testing]{\textbf{Effect of uncertainty in $\bar{R}_0$ on optimal strategy and decease rate, $0.01\%$ decease tolerance, fast testing.} Both subfigures consider the range $\bar{R}_0 \in [3.17, 3.38]$, {a fast testing policy} and a decease tolerance of $0.01\%$.
(left)~Ranges of optimal strategies, (right)  Ranges of \ak{aggregate deceases}  when the optimal strategy obtained  based on $\bar{R}_0 = 3.27$ and infection mortality rate of $0.66\%$ is implemented.
The darker line corresponds to $\bar{R}_0 = 3.27$.}
\vspace{-2mm}
\label{robust_R_8}
\end{figure}

\begin{figure}[H]
\begin{subfigure}{0.5\columnwidth}
\centering
\includegraphics[scale=0.57]{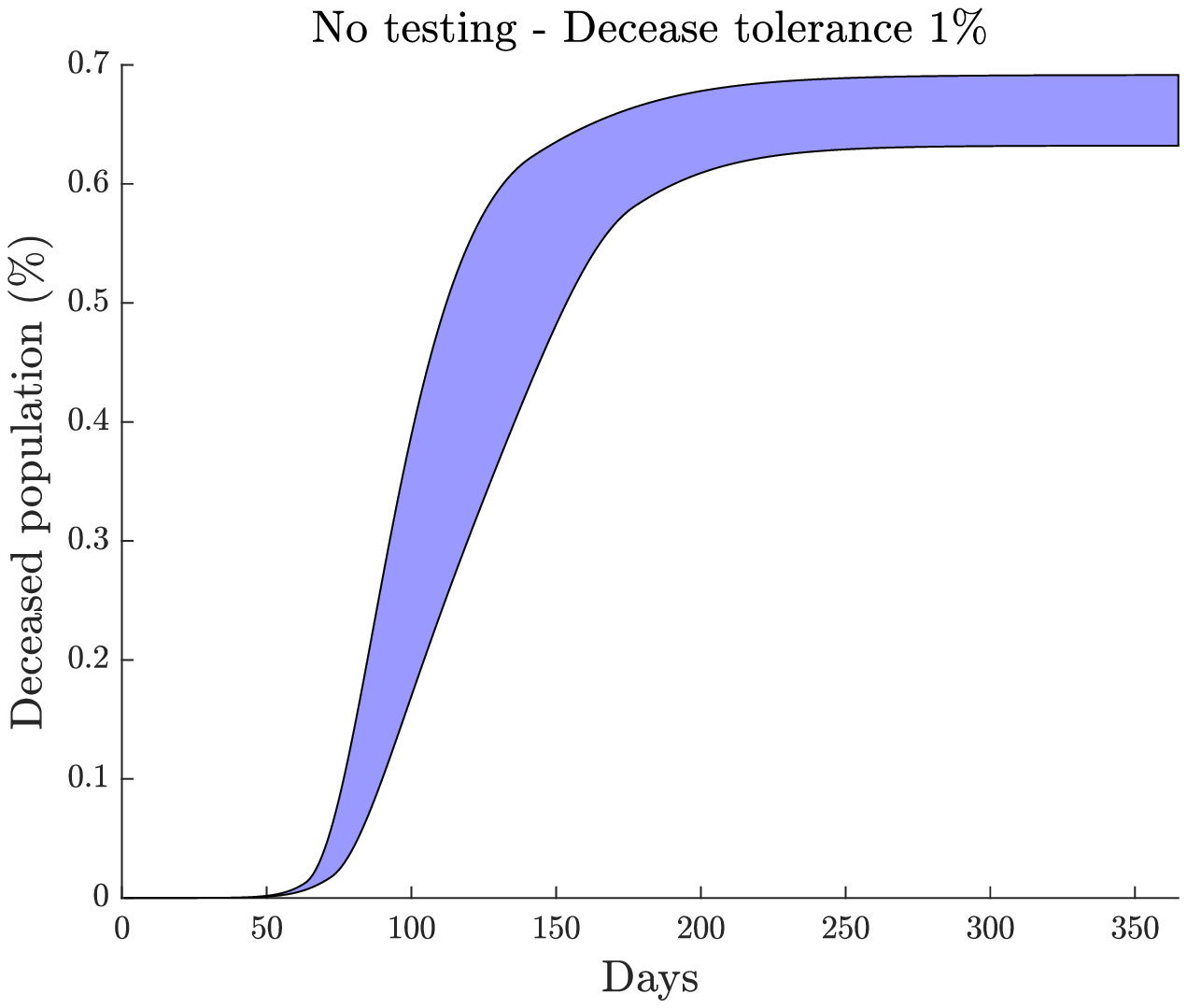}
\end{subfigure} 
\begin{subfigure}{0.5\columnwidth}
\centering
\includegraphics[scale=0.57]{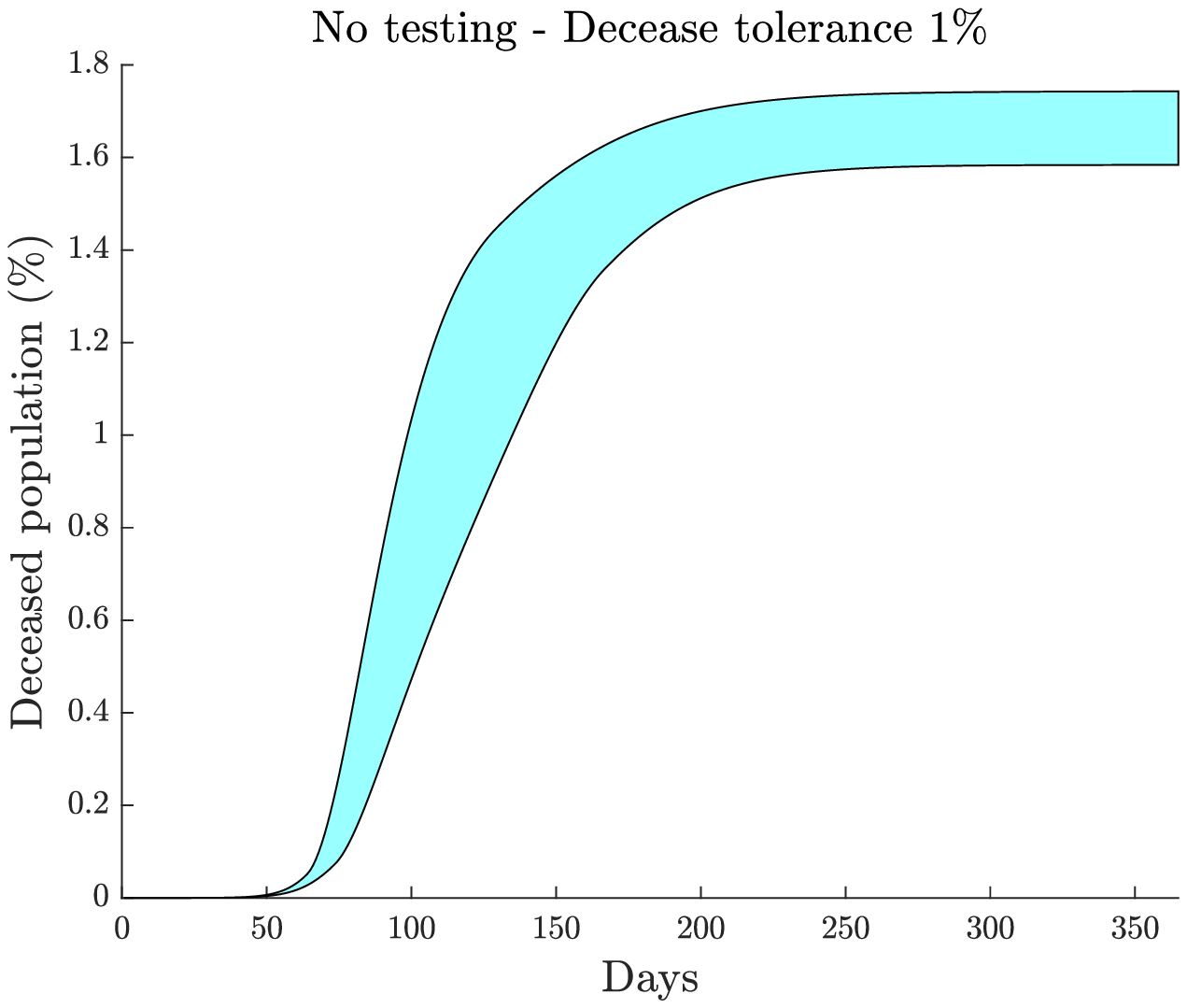}
\end{subfigure}
\vspace{-2mm} \caption[Combined effect of uncertainty in $\bar{R}_0$ and infection fatality rate on the decease rate, $1\%$ decease tolerance, no testing]{\textbf{Combined effect of uncertainty in $\bar{R}_0$ and infection fatality rate on the decease rate, $1\%$ decease tolerance, no testing.} Ranges of \ak{aggregate deceases} for $\bar{R}_0 \in [3.17, 3.38]$
when no testing is performed and a decease tolerance of $1\%$ is adopted with infection mortality rates of $0.39\%$ (left) and $1.33\%$ (right) when the optimal strategy obtained  based on $\bar{R}_0 = 3.27$ and infection mortality rate of $0.66\%$ is implemented.
} \vspace{-2mm}
\label{robust_R_9}
\end{figure}

\begin{figure}[H]
\begin{subfigure}{0.5\columnwidth}
\centering
\includegraphics[scale=0.55]{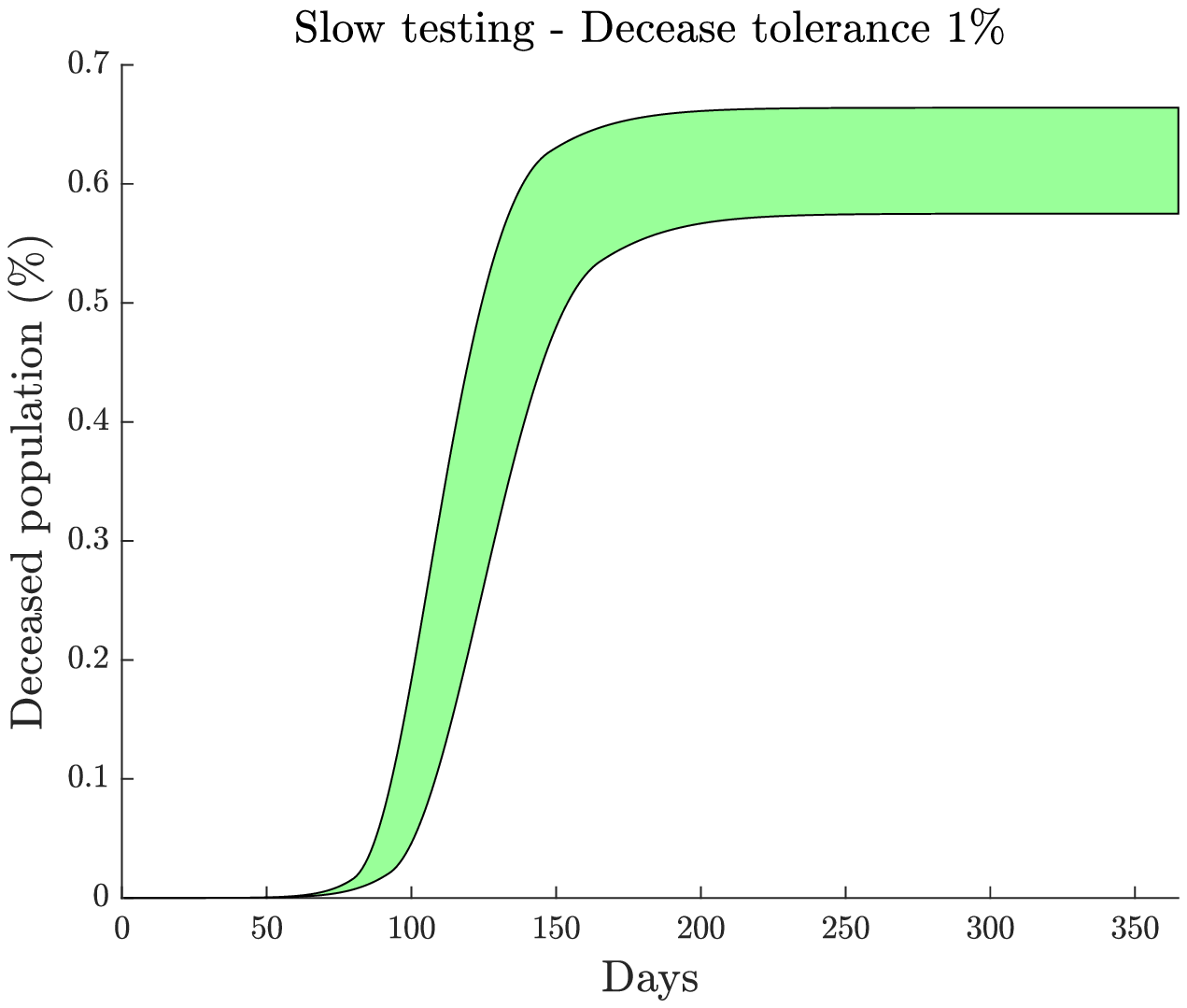}
\end{subfigure} 
\begin{subfigure}{0.5\columnwidth}
\centering
\includegraphics[scale=0.55]{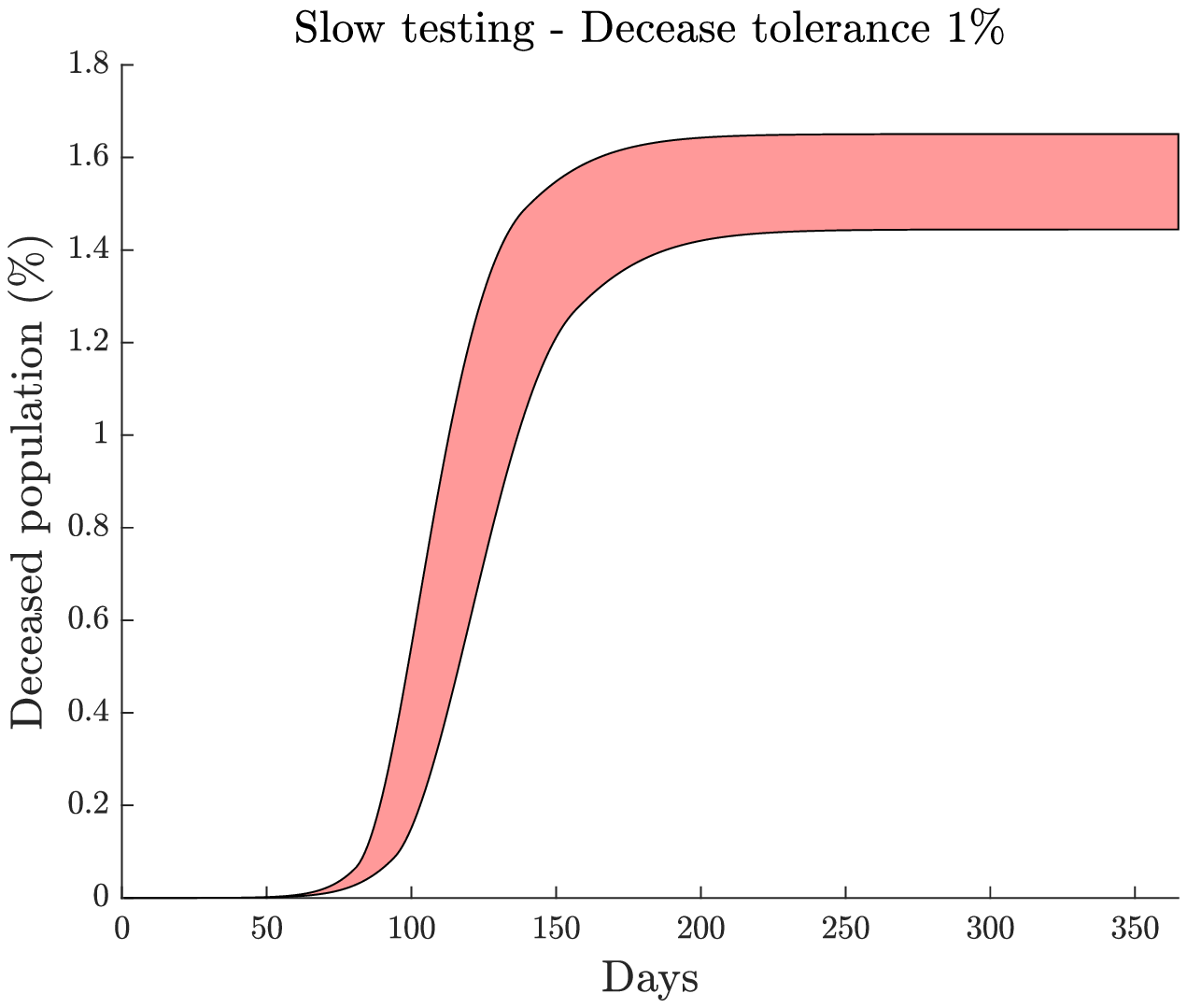}
\end{subfigure}
\vspace{-2mm} \caption[Combined effect of uncertainty in $\bar{R}_0$ and infection fatality rate on the decease rate, $1\%$ decease tolerance, slow testing]{\textbf{Combined effect of uncertainty in $\bar{R}_0$ and infection fatality rate on the decease rate, $1\%$ decease tolerance, slow testing.} Ranges of \ak{aggregate deceases} for $\bar{R}_0 \in [3.17, 3.38]$
when {a slow testing policy} and a decease tolerance of $1\%$ are adopted with infection mortality rates of $0.39\%$ (left) and $1.33\%$ (right) 
when the optimal strategy obtained  based on $\bar{R}_0 = 3.27$ and infection mortality rate of $0.66\%$ is implemented.
} \vspace{-2mm}
\end{figure}

\begin{figure}[H]
\begin{subfigure}{0.5\columnwidth}
\centering
\includegraphics[scale=0.55]{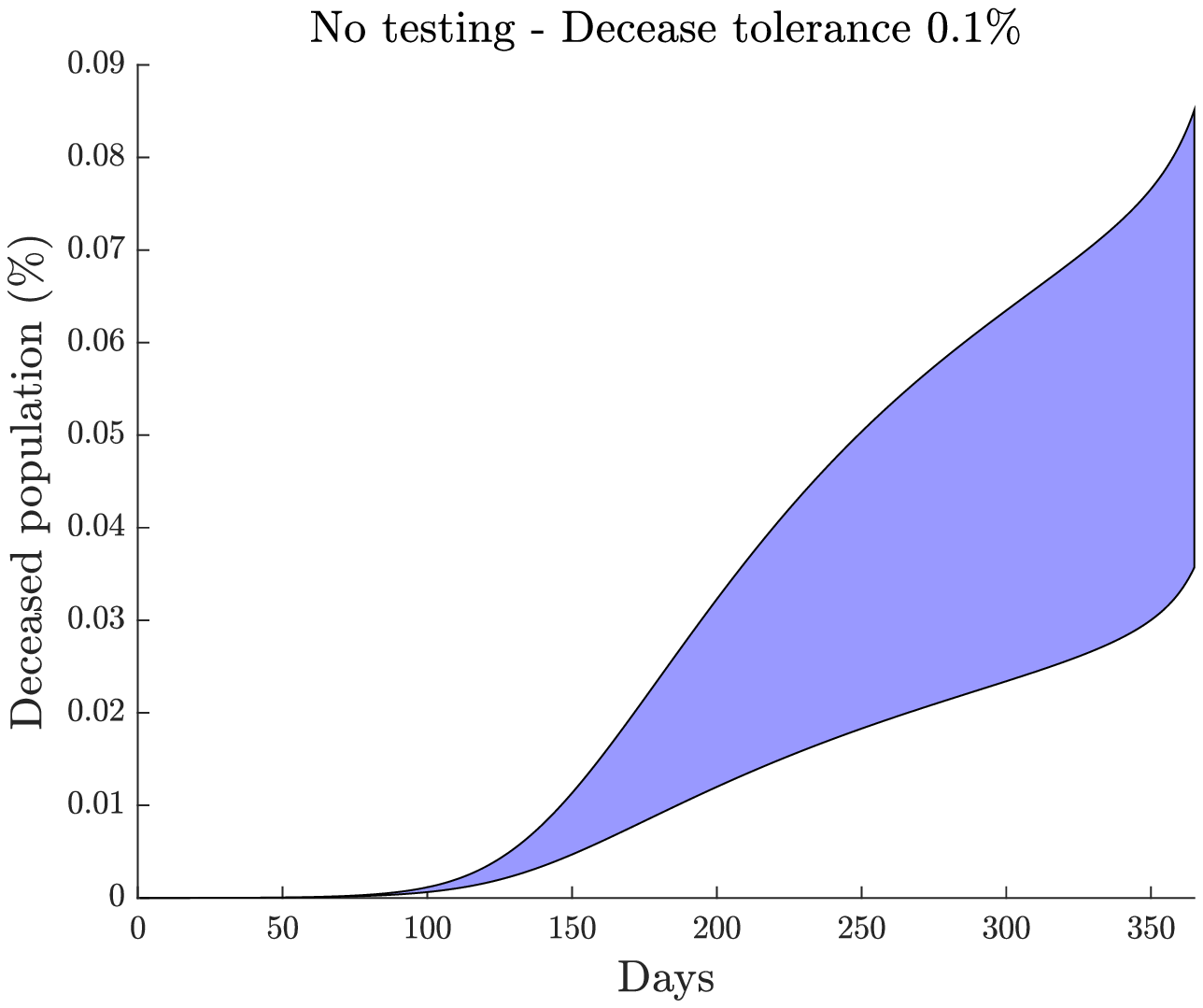}
\end{subfigure} 
\begin{subfigure}{0.5\columnwidth}
\centering
\includegraphics[scale=0.55]{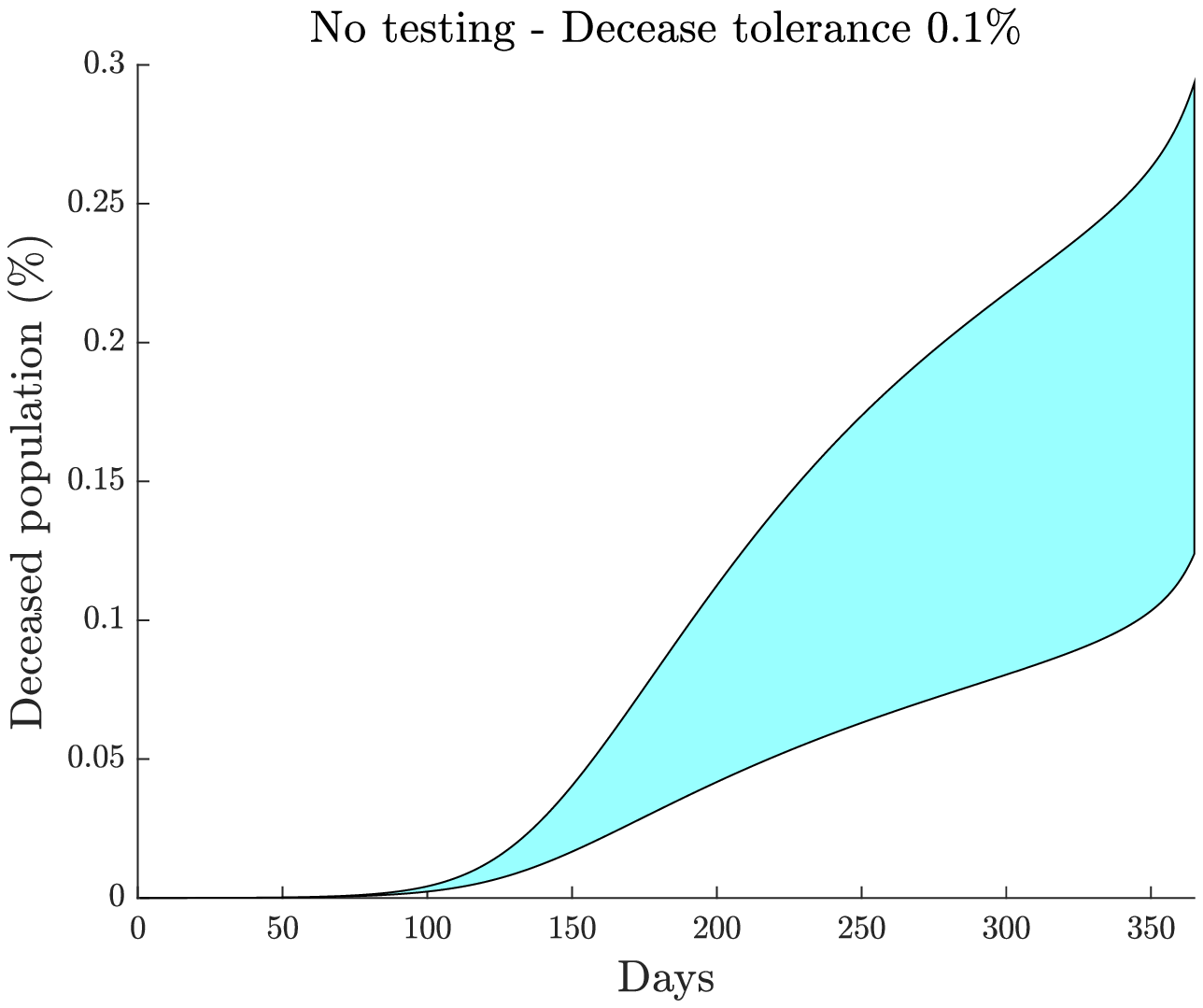}
\end{subfigure}
\vspace{-2mm} \caption[Combined effect of uncertainty in $\bar{R}_0$ and infection fatality rate on the decease rate, $0.1\%$ decease tolerance, no testing]{\textbf{Combined effect of uncertainty in $\bar{R}_0$ and infection fatality rate on the decease rate, $0.1\%$ decease tolerance, no testing.} Ranges of \ak{aggregate deceases} for $\bar{R}_0 \in [3.17, 3.38]$
when no testing is performed and a decease tolerance of $0.1\%$ is adopted with infection mortality rates of $0.39\%$ (left) and $1.33\%$ (right)
when the optimal strategy obtained  based on $\bar{R}_0 = 3.27$ and infection mortality rate of $0.66\%$ is implemented.
} \vspace{-2mm}
\end{figure}

\begin{figure}[H]
\begin{subfigure}{0.5\columnwidth}
\centering
\includegraphics[scale=0.55]{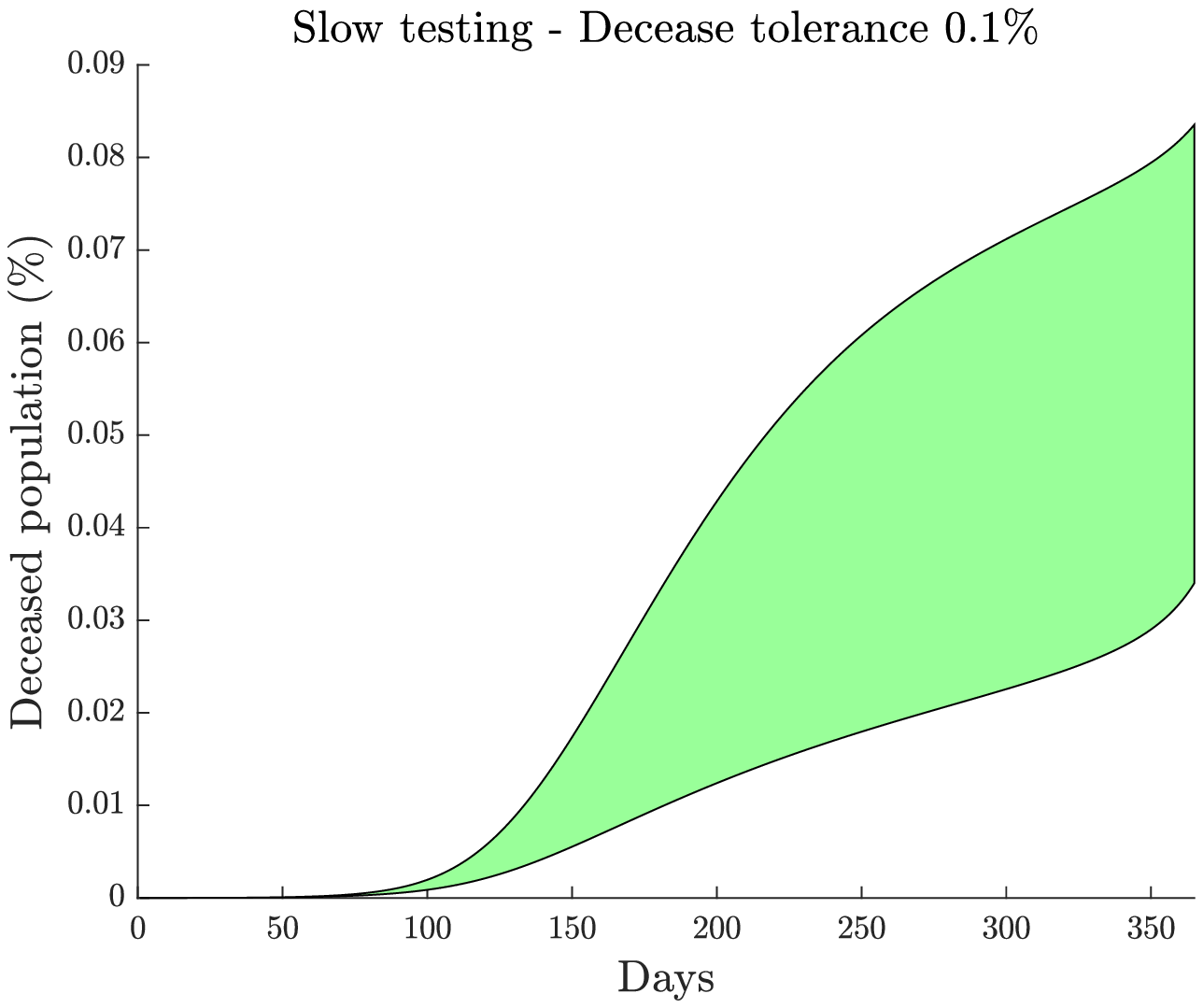}
\end{subfigure} 
\begin{subfigure}{0.5\columnwidth}
\centering
\includegraphics[scale=0.55]{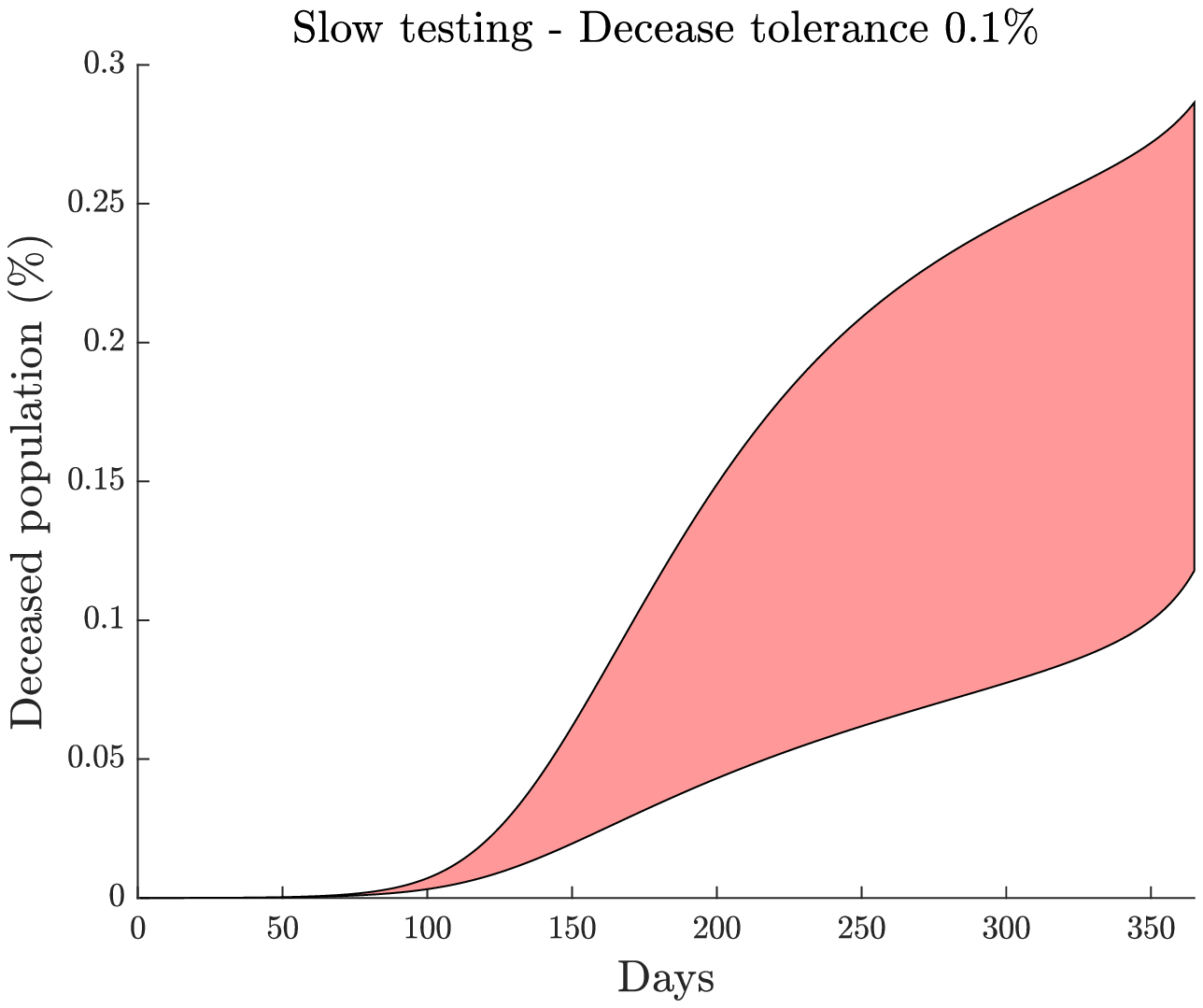}
\end{subfigure}
\vspace{-2mm} \caption[Combined effect of uncertainty in $\bar{R}_0$ and infection fatality rate on the decease rate, $0.1\%$ decease tolerance, slow testing]{\textbf{Combined effect of uncertainty in $\bar{R}_0$ and infection fatality rate on the decease rate, $0.1\%$ decease tolerance, slow testing.} Ranges of \ak{aggregate deceases} for $\bar{R}_0 \in [3.17, 3.38]$
when {a slow testing policy} and a decease tolerance of $0.1\%$ are adopted with infection mortality rates of $0.39\%$ (left) and $1.33\%$ (right)
when the optimal strategy obtained  based on $\bar{R}_0 = 3.27$ and infection mortality rate of $0.66\%$ is implemented.
}\vspace{-2mm}
\end{figure}

\begin{figure}[H]
\begin{subfigure}{0.5\columnwidth}
\centering
\includegraphics[scale=0.55]{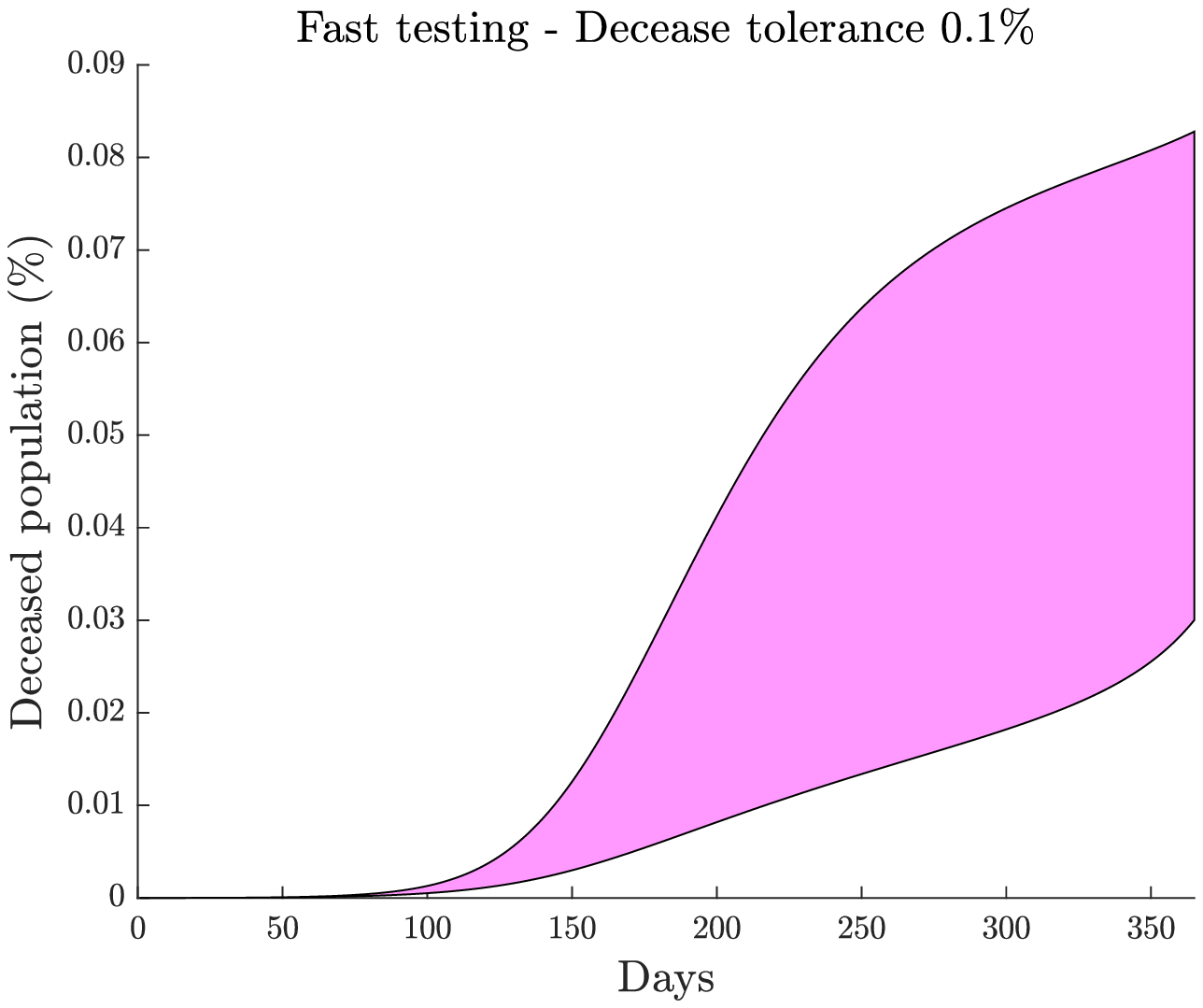}
\end{subfigure} 
\begin{subfigure}{0.5\columnwidth}
\centering
\includegraphics[scale=0.55]{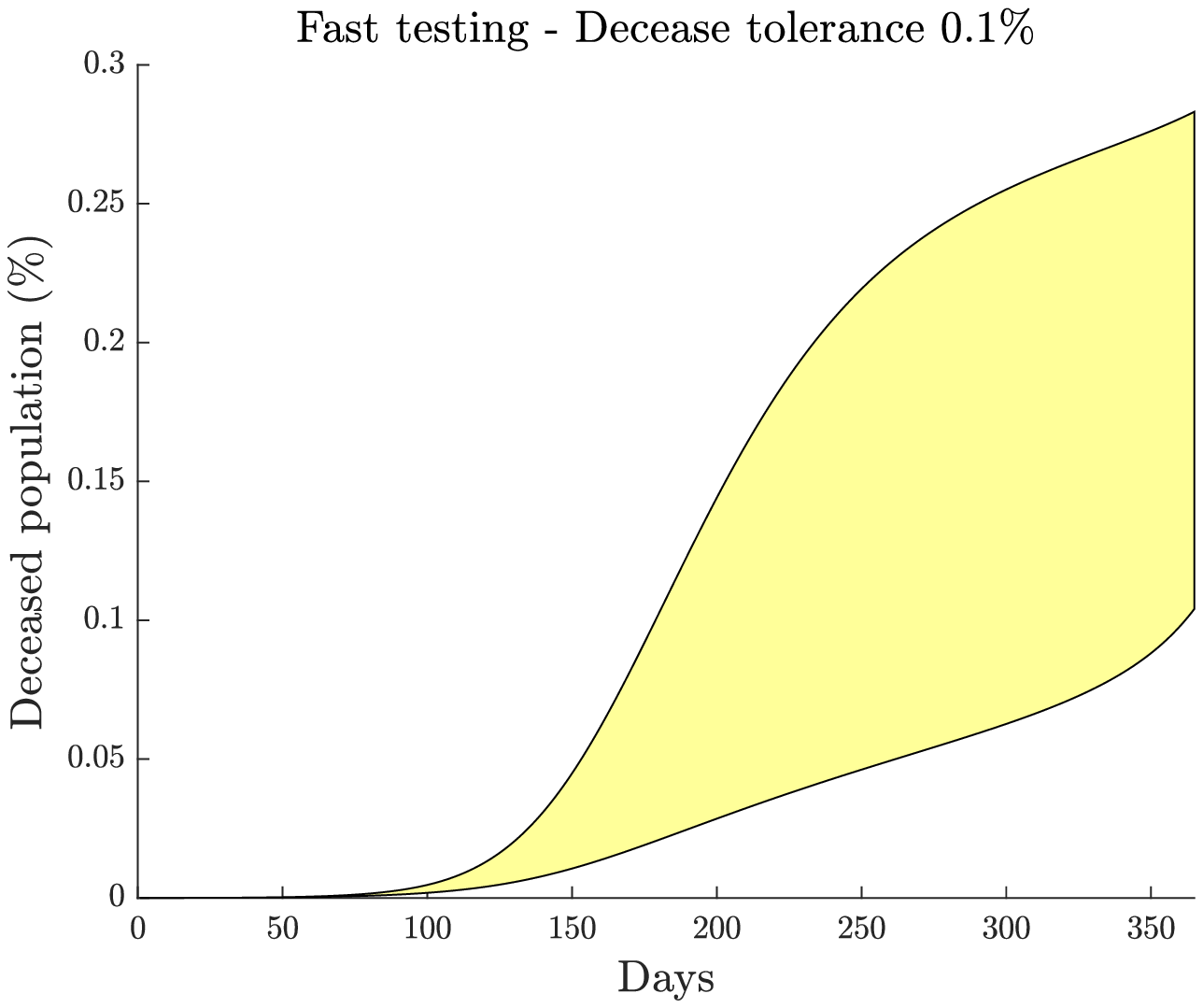}
\end{subfigure}
\vspace{-2mm} \caption[Combined effect of uncertainty in $\bar{R}_0$ and infection fatality rate on the decease rate, $0.1\%$ decease tolerance, fast testing]{\textbf{Combined effect of uncertainty in $\bar{R}_0$ and infection fatality rate on the decease rate, $0.1\%$ decease tolerance, fast testing.} Ranges of \ak{aggregate deceases} for $\bar{R}_0 \in [3.17, 3.38]$
when {a fast testing policy} and a decease tolerance of $0.1\%$ are adopted with infection mortality rates of $0.39\%$ (left) and $1.33\%$ (right)
when the optimal strategy obtained  based on $\bar{R}_0 = 3.27$ and infection mortality rate of $0.66\%$ is implemented.
}\vspace{-2mm}
\end{figure}

\begin{figure}[H]
\begin{subfigure}{0.5\columnwidth}
\centering
\includegraphics[scale=0.55]{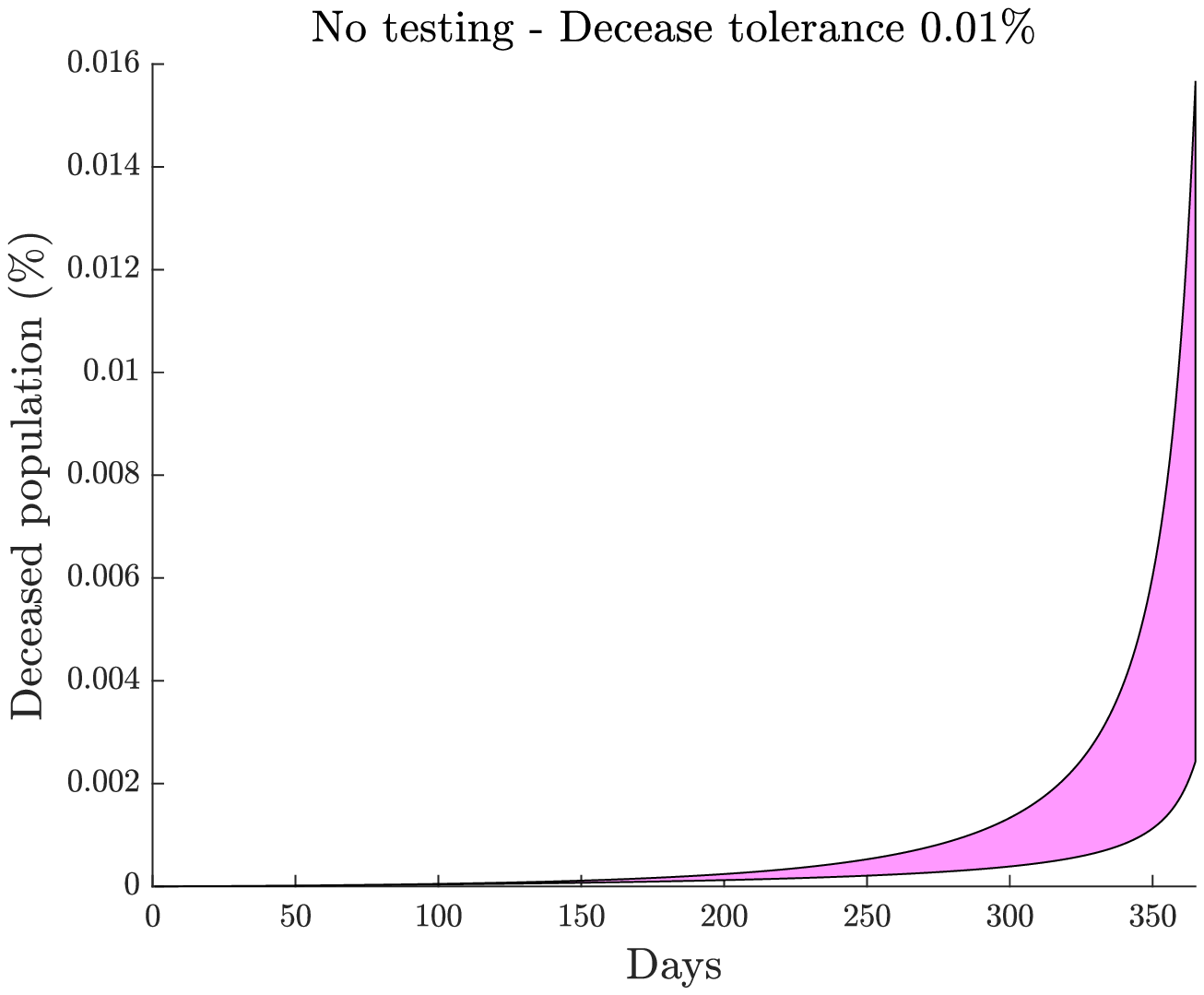}
\end{subfigure} 
\begin{subfigure}{0.5\columnwidth}
\centering
\includegraphics[scale=0.55]{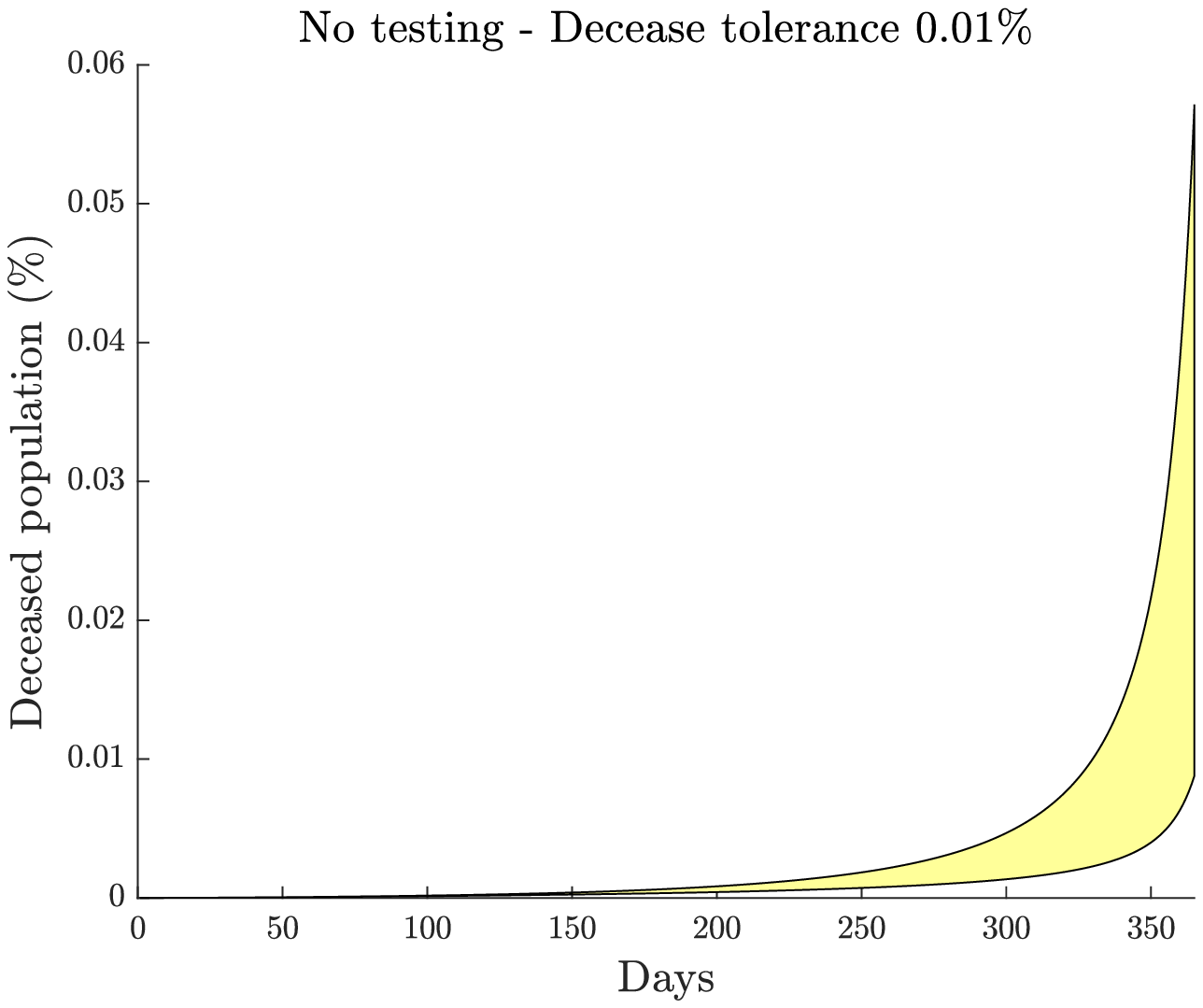}
\end{subfigure}
\vspace{-2mm} \caption[Combined effect of uncertainty in $\bar{R}_0$ and infection fatality rate on the decease rate, $0.01\%$ decease tolerance, no testing]{\textbf{Combined effect of uncertainty in $\bar{R}_0$ and infection fatality rate on the decease rate, $0.01\%$ decease tolerance, no testing.} Ranges of \ak{aggregate deceases} for $\bar{R}_0 \in [3.17, 3.38]$
when no testing is performed and a decease tolerance of $0.01\%$ is adopted with infection mortality rates of $0.39\%$ (left) and $1.33\%$ (right)
when the optimal strategy obtained  based on $\bar{R}_0 = 3.27$ and infection mortality rate of $0.66\%$ is implemented.
} \vspace{-2mm}
\end{figure}

\begin{figure}[H]
\begin{subfigure}{0.5\columnwidth}
\centering
\includegraphics[scale=0.55]{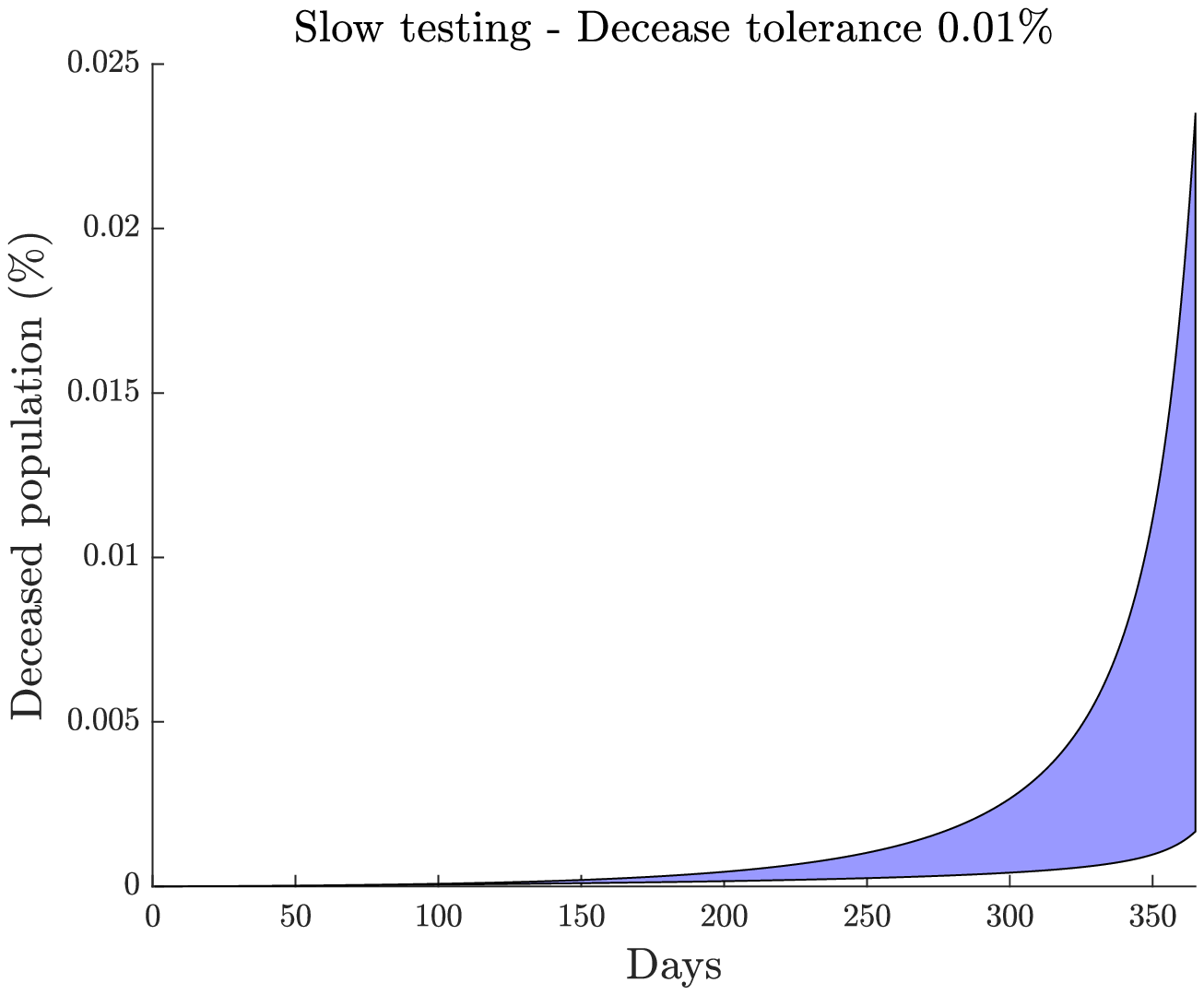}
\end{subfigure} 
\begin{subfigure}{0.5\columnwidth}
\centering
\includegraphics[scale=0.55]{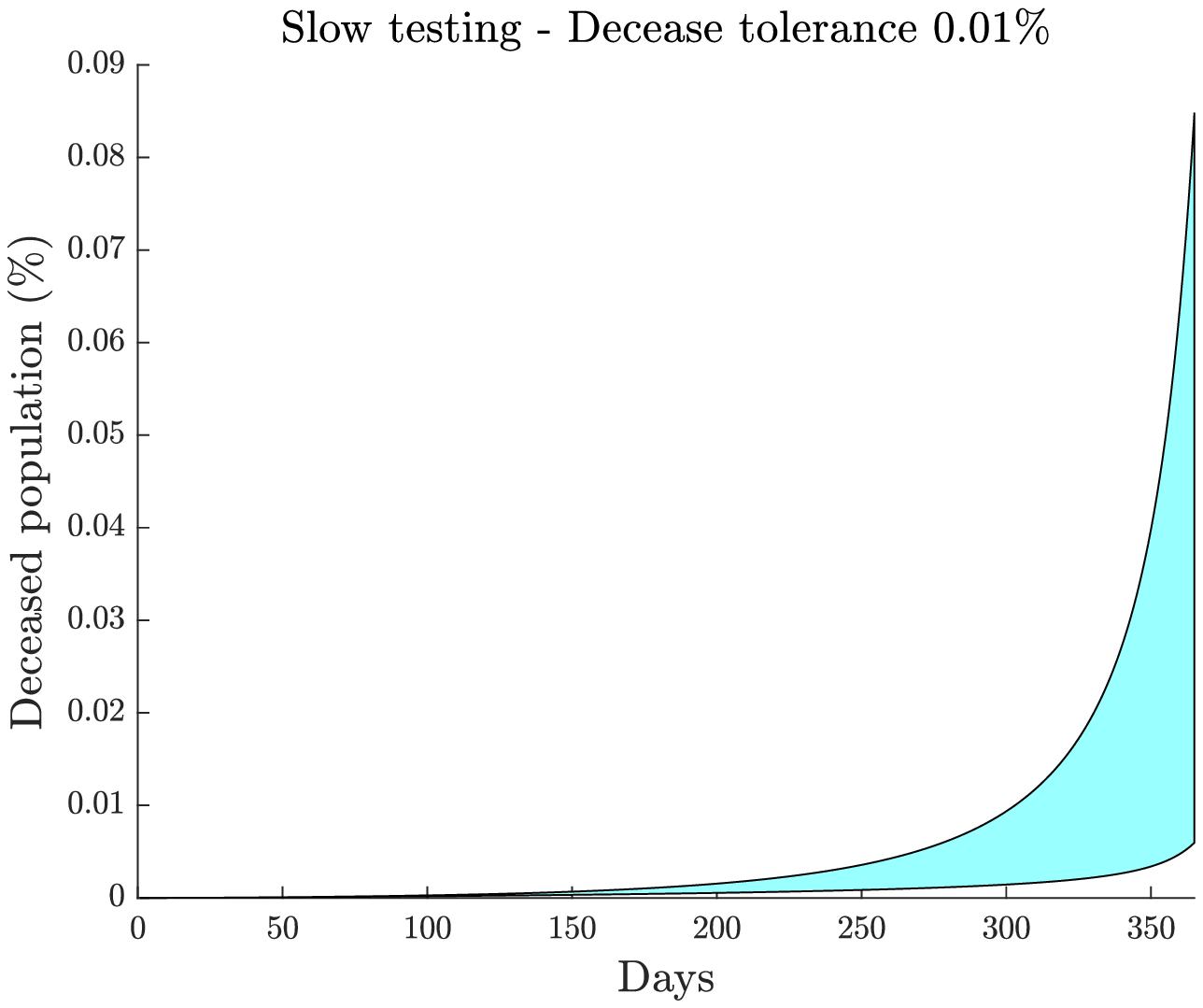}
\end{subfigure}
\vspace{-2mm} \caption[Combined effect of uncertainty in $\bar{R}_0$ and infection fatality rate on the decease rate, $0.01\%$ decease tolerance, slow testing]{\textbf{Combined effect of uncertainty in $\bar{R}_0$ and infection fatality rate on the decease rate, $0.01\%$ decease tolerance, slow testing.} Ranges of \ak{aggregate deceases} for $\bar{R}_0 \in [3.17, 3.38]$
when {a slow testing policy} and a decease tolerance of $0.01\%$ are adopted with infection mortality rates of $0.39\%$ (left) and $1.33\%$ (right)
when the optimal strategy obtained  based on $\bar{R}_0 = 3.27$ and infection mortality rate of $0.66\%$ is implemented.
}\vspace{-2mm}
\end{figure}

\begin{figure}[H]
\begin{subfigure}{0.5\columnwidth}
\centering
\includegraphics[scale=0.57]{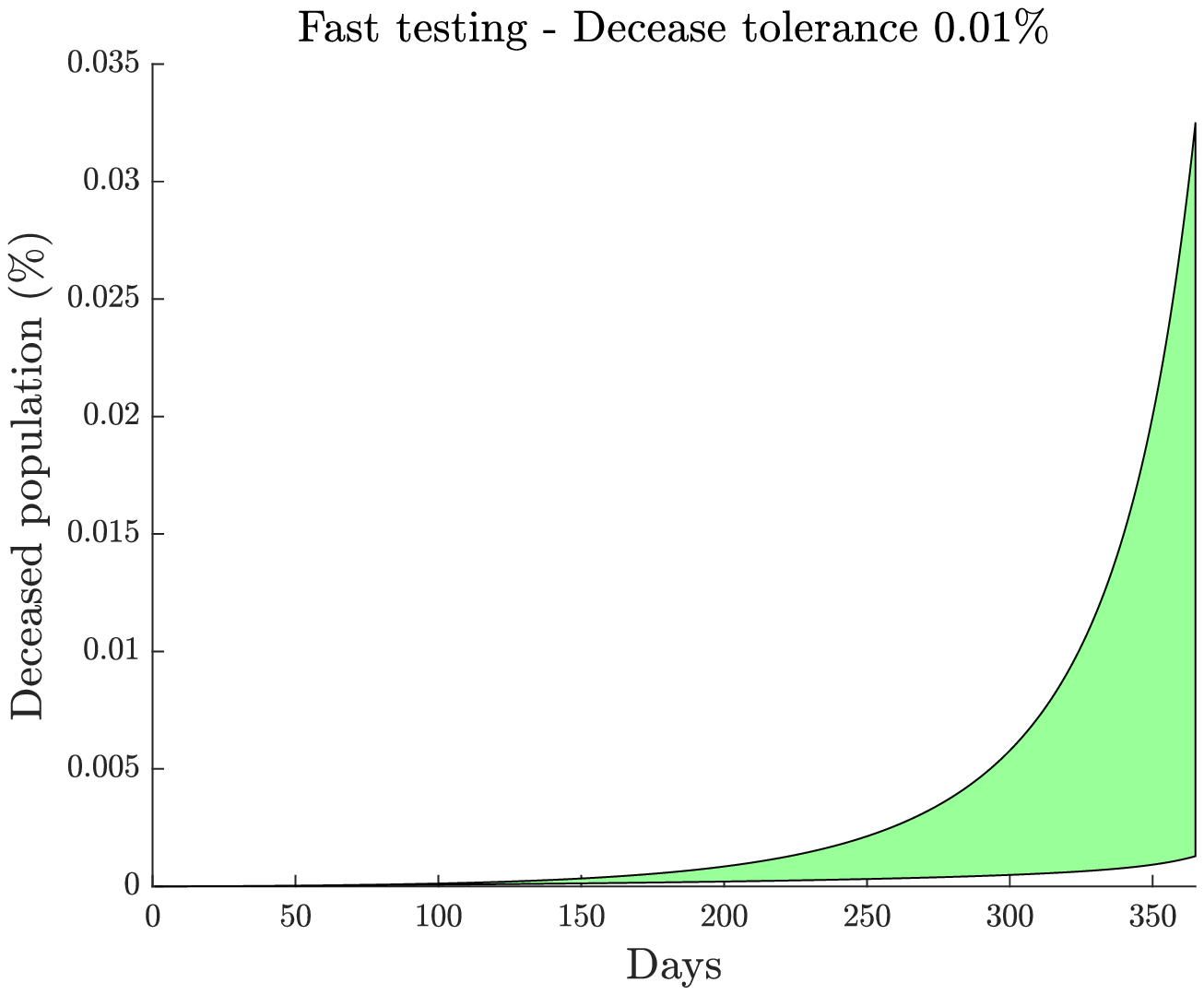}
\end{subfigure} 
\begin{subfigure}{0.5\columnwidth}
\centering
\includegraphics[scale=0.57]{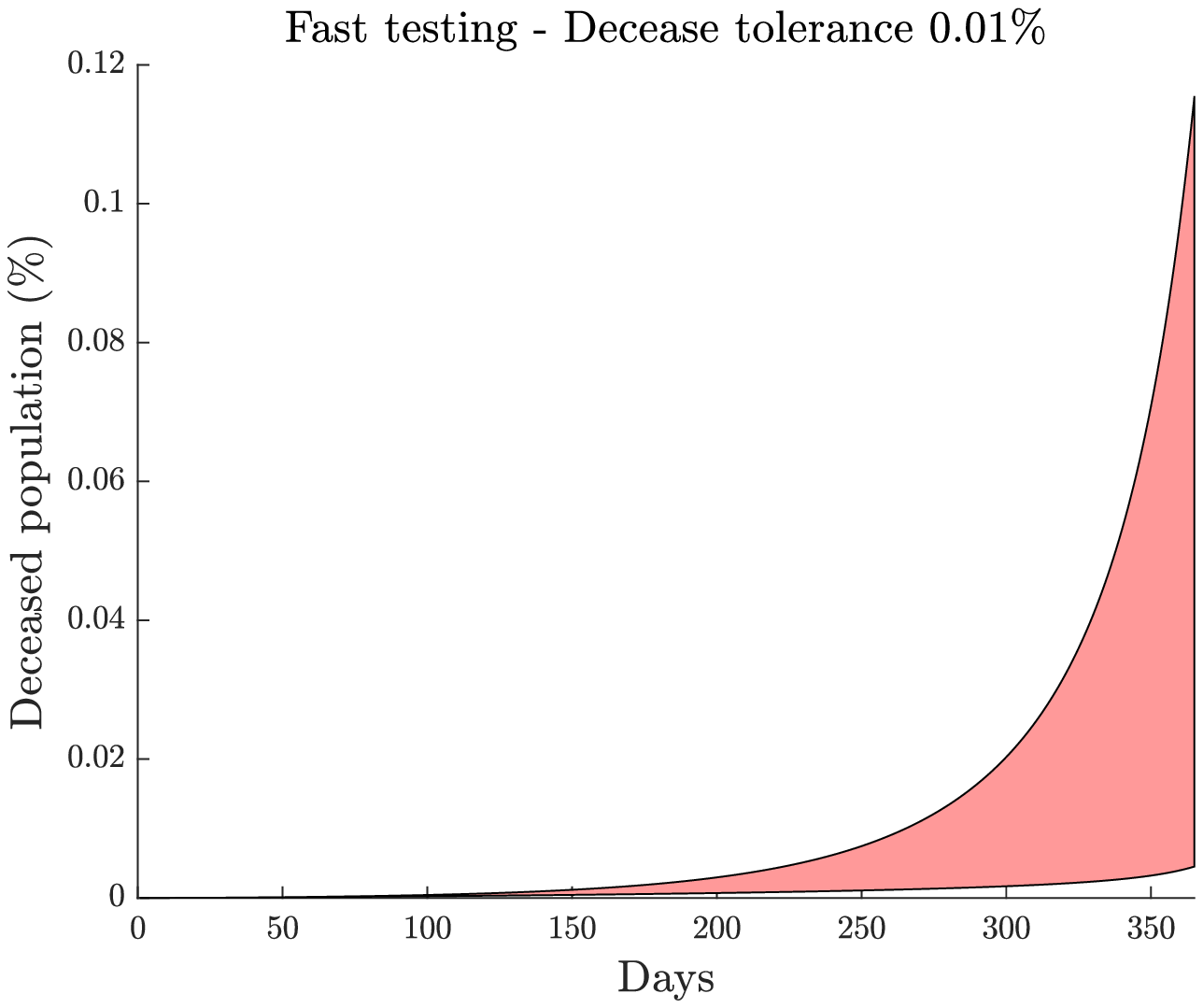}
\end{subfigure}
\vspace{-2mm} \caption[Combined effect of uncertainty in $\bar{R}_0$ and infection fatality rate on the decease rate, $0.01\%$ decease tolerance, fast testing]{\textbf{Combined effect of uncertainty in $\bar{R}_0$ and infection fatality rate on the decease rate, $0.01\%$ decease tolerance, fast testing.} Ranges of \ak{aggregate deceases} for $\bar{R}_0 \in [3.17, 3.38]$
when {a fast testing policy} and a decease tolerance of $0.01\%$ are adopted with infection mortality rates of $0.39\%$ (left) and $1.33\%$ (right)
when the optimal strategy obtained  based on $\bar{R}_0 = 3.27$ and infection mortality rate of $0.66\%$ is implemented.
}
\vspace{-2mm}
\label{robust_R_16}
\end{figure}



%
%



\subsubsection*{Acknowledgements}

This work was funded by the European Union’s Horizon 2020 research and innovation program under grant agreement 739551 (KIOS CoE) and from the Republic of Cyprus through the Directorate General for European Programs, Coordination, and Development.

%
%

\end{document}